\documentclass[11pt]{article}

\usepackage{amstext,amssymb,amsmath,amsbsy}
\usepackage{hyperref}
\usepackage{amscd}
\usepackage{amsfonts}
\usepackage{indentfirst}
\usepackage{verbatim}
\usepackage{amsmath}
\usepackage{amsthm}
\usepackage{enumerate}
\usepackage{graphicx}
\usepackage{color, soul}
\usepackage[OT1]{fontenc}
\usepackage[latin1]{inputenc}
\usepackage[english]{babel}
\usepackage{amssymb}
\usepackage{subfig}
\usepackage{algorithm}
\usepackage{algpseudocode}
\usepackage{mathtools}

\newcommand{\R}{\mathbb{R}}
\newcommand{\N}{\mathbb{N}}
\newcommand{\ds}{\displaystyle}

\newcommand{\x}{{\bf x}}
\newcommand{\p}{{\bf p}}

\newcommand{\Div}{{\rm div}}

\newcommand*{\QEDB}
{\null\nobreak\hfill\ensuremath{\square}}

\newtheorem{Theorem}{Theorem}[section]

\newtheorem{remark}{Remark}[section]

\newtheorem{Assumption}{Assumption}[section]
\newtheorem*{Assumption*}{Assumption}
\newtheorem{Definition}{Definition}[section]

\newtheorem*{problem*}{Problem}
\numberwithin{equation}{section}

\begin{document}

\title{Numerical viscosity  solutions to Hamilton-Jacobi equations via a Carleman estimate and the  convexification method}

\author{Michael Klibanov\thanks{Department of Mathematics and Statistics, University of North Carolina at
Charlotte, Charlotte, NC, 28223, USA, \texttt{mklibanv@uncc.edu}.} \and Loc H. Nguyen\thanks{Department of Mathematics and Statistics, University of North Carolina at
Charlotte, Charlotte, NC, 28223, USA, \texttt{loc.nguyen@uncc.edu}.} \and Hung V. Tran\thanks{Department of Mathematics, University of Wisconsin Madison, Madison, WI, 53706, USA, \texttt{hung@math.wisc.edu} (corresponding
author).}} 


\date{}
\maketitle
\begin{abstract}
	We propose a globally convergent numerical method, called the convexification, to numerically compute the viscosity solution to first-order Hamilton-Jacobi equations through the vanishing viscosity process where the viscosity parameter is a fixed small number. 
By convexification, we mean that we employ a suitable Carleman weight function to convexify the cost functional defined directly from the form of the Hamilton-Jacobi equation under consideration.
The strict convexity of this functional is rigorously proved using a new Carleman estimate.
We also prove that the unique minimizer of the this strictly convex functional can be reached by the gradient descent method. Moreover,  we show that the minimizer well approximates the viscosity solution of the Hamilton-Jacobi equation as the noise contained in the boundary data tends to zero.
Some interesting numerical illustrations are presented. 
\end{abstract}



\section{Introduction}

The aim of this paper is to compute viscosity solutions to a large class of Hamilton-Jacobi equations possibly involving  nonconvex Hamiltonians. 
The key ingredient for us to reach this achievement is the use of a new Carleman estimate and the convexification method. 
This method is only applicable when 
a  viscosity term is added to the Hamilton-Jacobi equation under consideration. The idea of adding the viscosity term and passing to the limit to obtain viscosity solutions is due to the seminal works \cite{CrandallLions83,CrandallEvansLions84}.
Let $R > 0$ and $\Omega = (-R, R)^d$ where $d \geq 1$ is the spatial dimension.
Let $F: \overline \Omega \times \R \times \R^d \to \R$ and $f: \overline \Omega \to \R$ be functions of the class $C^2$.
In this paper, we propose a globally convergent numerical method to solve the Hamilton-Jacobi equation
\begin{equation}
	F(\x, u(\x), \nabla u(\x)) = 0 \quad \mbox{for all } \x \in \Omega
	\label{HJ}
\end{equation}
with the Dirichlet boundary condition
\begin{equation}
	u(\x) = f(\x) \quad  \mbox{for all }  \x \in \partial \Omega.
	\label{dir}
\end{equation}
The smoothness condition imposing on $F$ and $f$ is for the simplicity. We need it when analytically establishing the convergence of the proposed method. 
	However, this technical condition can be relaxed in our numerical study. 
In this paper, we are interested in computing the viscosity solution to \eqref{HJ}--\eqref{dir}.
We only deal with the case that the Dirichlet boundary condition holds in the classical sense in this paper.
We  refer the readers to \cite{CrandallLions83, CrandallEvansLions84, Lions, Barles, BCD,Tran19} and the references therein for the theory of viscosity solutions to  \eqref{HJ}--\eqref{dir}.
It is worth mentioning that a number of different  extremely efficient and fast numerical approaches and techniques (many of which are of high orders) have been developed for Hamilton-Jacobi equations.
For finite difference monotone and consistent schemes of first-order equations and applications, see \cite{CL-rate, Sou1, BS-num, Sethian, OsFe} for details and recent developments.
If $F=F(\x,s,\p)$ is convex in $\p$ and satisfies some appropriate conditions, it is possible to construct some semi-Lagrangian approximations by the discretization of the Dynamical Programming Principle associated to the problem, see \cite{FaFe1, FaFe2} and the references therein.
For a non-exhaustive list of results along these directions, see \cite{OsherSethian, OsherShu, Tsitsiklis, Abgrall, Abgrall2, BrysonLevy, SethianVladimirsky, TsaiChengOsherZhao, KaoOsherQian, QianZhangZhao, CamilliCDGomes, CockburnMerevQian, CagnettiGomesTran, ObermanSalvador, LiQian, GallistlSprekelerSuli}.
Another approach to solve  \eqref{HJ}--\eqref{dir} is based on optimization \cite{HornBrooks, Szeliski, LeclercBobick, DanielDurou}. 
However, due to the nonlinearity of the function $F$ in \eqref{HJ}, the least squares cost functional is nonconvex and might have multiple local minima and ravines.
Hence, the methods based on optimization only provide reliable numerical solutions if good initial guesses of the true solutions are given. 

\smallskip

Unlike the mentioned optimization approach,  we propose to use the convexification method, which
does not rely on the above assumptions. 
This method is globally convergent in the sense that 
\begin{enumerate}
	\item it delivers good approximation of the true solution without knowing any advance knowledge of the true solution even when the given data is noisy;
	\item the claim in \# 1 above is rigorously proved and numerically verified.
\end{enumerate} 
To the best of our knowledge, our method is new in the context of viscosity solutions to Hamilton-Jacobi equations.
It has two major advantages in solving \eqref{HJ}--\eqref{dir} numerically.
Firstly, it works for some quite general $F(\x,u,\nabla u)$, which might be nonconvex in $\nabla u$, and it does not require a lot of specific structures on $F$.
 In particular, $F$ might be dependent on $u$ and $\nabla u$ in a rather complicated way.
Secondly, it is quite stable and robust even with some noise level on the boundary data, which occurs naturally in applications.

The main idea of the convexification method is to employ a suitable Carleman
weight function to convexify the mismatch functional derived from the given
boundary value problem. Several versions of the convexification method have
been developed since it was first introduced in \cite%
{KlibanovIoussoupova:SMA1995} for a coefficient inverse problem for a hyperbolic equation. We cite here \cite%
{Klibanov:ip2015,Klibanov:sjma1997, Klibanov:IPI2019,Klibanov:nw1997,KlibanovNik:ra2017, KhoaKlibanovLoc:SIAMImaging2020, VoKlibanovNguyen:IP2020, SmirnovKlibanovNguyen:IPI2020, Klibanov:ip2020} and references therein
for some important works in this area and their real-world applications in bio-medical imaging, non-destructed testing, 
travel time tomography, identifying anti-personnel explosive devices buried under the ground, {\it etc.} 
The crucial mathematical ingredients  that
guarantee the strict convexity of this functional are the
Carleman estimates. 
The original idea of applying Carleman estimates to prove the uniqueness for a large class of important nonlinear mathematical problems was first published in \cite{BukhgeimKlibanov:smd1981}. 
It was discovered later in \cite{KlibanovIoussoupova:SMA1995, KlibanovLiBook}, that the idea of \cite%
{BukhgeimKlibanov:smd1981} can be successfully modified to develop globally
convergent numerical methods for coefficient inverse problems using the
convexification.

In this paper, it is the first time we use a Carleman weight function to numerically solve Hamilton-Jacobi equations.
One of the strengths of the convexification method is that it does not require the convexity of $F(\x, s, \p)$ with respect to $\p$. Still, it has a drawback. The theory of the convexification method requires
 an additional information about $u_z$ on a part of $\partial \Omega$. In this paper, that part is 
\begin{equation}
		\Gamma^+ = \big\{
			\x = (x_1, x_2, \dots, x_{d - 1}, z = R): |x_i| \leq R, 1 \leq i \leq d - 1
		\big\}
		\subset \partial \Omega.
	\end{equation}
For $\x\in \R^d$, write $\x = (x_1, x_2, \dots, x_{d - 1}, z)$.
More precisely, we impose the following condition.
\begin{Assumption}
	Write $\nabla u = (u_{x_1}, \dots, u_{x_{d - 1}}, u_z)$. We assume that  $u_z$ on $\Gamma^+$ is known.
	\label{assumption 1}
\end{Assumption}

In general, the additional knowledge of $u_z$ on $\Gamma^+$ in Assumption \ref{assumption 1} makes the problem of computing solutions to Hamilton-Jacobi equations with the Dirichlet data on $\partial \Omega$ and the Neumann data on $\Gamma^+$  over-determined.
However, in many real-world circumstances, we are able to compute $u_z$ on $\Gamma^+$ from the knowledge of $u$ on $\Gamma^+$ without further measurement.
We provide here a classical example arising from  the traveling time problem. Denote by $c(\x)$, $\x \in \overline \Omega$, the velocity of the light at the point $\x$. 
Here, $c \in C(\overline \Omega, (0,\infty))$ is a given function.
Let $u(\x)$ be the minimal time for the light to travel from $\partial \Omega$ to $\x \in \Omega$.
This function is governed by the boundary value problem for the eikonal equation
	\begin{equation}
		\left\{
			\begin{array}{rcll}
				c(\x)^2|\nabla u(\x)|^2 &=& 1 &\x \in \Omega,\\
				u(\x) &=& 0 &\x \in \partial \Omega.
			\end{array}
		\right.
		\label{eik}
	\end{equation}
	Since $u = 0$ on $\partial \Omega$, in particular, $u = 0$ on $\Gamma^+.$
	Hence $u_{x_i} = 0$ on $\Gamma^+$ for all $1 \leq i \leq d - 1$. 
	The function $u_z$ on $\Gamma^+$ is given by 
	\[
		u_z(\x) = -\frac{1}{c(\x)} \quad \mbox{for all } \x \in \Gamma^+.
	\]
	Above, the case  $u_z(\x) = \frac{1}{c(\x)}$ is negligible since in reality, $u \geq 0$ on $\overline \Omega$ and therefore, its partial derivative with respect to $z$ on $\Gamma^+$ is non-positive.

\begin{remark}[Reducing Assumption \ref{assumption 1}]
We have the following points.
\begin{itemize}
\item[1.]	It follows from the example above that in general, since $u$, and; therefore, $u_{x_1},$ $\dots,$ $u_{x_{d - 1}}$ on $\Gamma^+$ are known, to verify Assumption \ref{assumption 1},  we simply solve  the equation 
\begin{equation}
		F(\x, u, u_{x_1}, \dots, u_{x_{d - 1}}, u_z) = 0
		\label{1.5}
\end{equation}
	for $u_z$. 
	So, a condition on $F$ such that Assumption \ref{assumption 1} holds true is that \eqref{1.5} is uniquely solvable for $u_z$. For example,
	\[
		\frac{\partial }{\partial z}F(\x, u, u_{x_1}, \dots, u_{x_{d - 1}}, u_z) \not = 0 \quad \mbox{for all } \x \in \Gamma^+.
	\]	
	
\item[2.] In many important Hamilton-Jacobi equations; for e.g., the eikonal equation $F(x, u, \nabla u) = c(\x)^2|\nabla u|^2 - 1$ in \eqref{eik} or $F(\x, u, \nabla u) = |u_{x_1}| - |u_{z}| - g$ (the case $d = 2$), for some function $g$, equation  \eqref{1.5} only provides $|u_z|$ rather than $u_z$ on $\Gamma^+$. In this case, the sign of $u_z$ on $\Gamma^+$ is required.

\item[3.]We provide here an example in which the sign of $u_{z}$ is known.
Let $d=3$ and $\mathbf{x}_{0}$ be a point in $\mathbb{R}^{3}\setminus \Omega 
$. In travel time tomography, the Hamilton-Jacobi equation that describes
the the travel time of light traveling  from the source $\mathbf{x}_{0}$ to
a point $\mathbf{x}\in \mathbb{R}^{3},$ has the form $c^{2}(\mathbf{x}%
)|\nabla _{\mathbf{x}}u(\mathbf{x},\mathbf{x}_{0})|^{2}=1$ for all $\mathbf{x%
}\in \mathbb{R}^{3}\setminus \{\mathbf{x}_{0}\}$ with $u(\mathbf{x}_{0},%
\mathbf{x}_{0})=0.$ Here $c(\mathbf{x})$ is the speed of the light at $%
\mathbf{x}$. It was proved in \cite[Lemma 4.1]{Klibanov:IPI2019} that if $c(%
\mathbf{x})$, $\mathbf{x}\in \Omega $ is an increasing function with respect
to $z$ and the source $x_{0}\in \left\{ z<-R\right\} ,$
then the function $u(x,x_{0})$ is strictly increasing in the $z-$%
direction for $x\in \Omega $\ implying $u_{z}>0$\
on $\Gamma ^{+}.$ This result can be extended to all dimensions
by repeating the proof in \cite[Lemma 4.1]{Klibanov:IPI2019}.
	
\item[4.] In the case when Assumption \ref{assumption 1} cannot be verified or even when it might not hold true, for e.g., 
\[
	F(\x, u, \nabla u) = 10 u + \min\{|\nabla u|, ||\nabla u| - 8| + 6\} - g
\] or 
\[ F(\x, u, \nabla u) =u+ |\nabla u| - V\cdot \nabla u 
\]
for some function $g$ and vector valued function $V$,  the convexification method still provides good numerical solutions, see Test 4 and Test 5 in Section \ref{sec:num}. However, the rigorous theorem that guarantees the efficiency of this method is missing in this paper.
\end{itemize}
\label{rm 1}
\end{remark}

\medskip

Let us give a brief description of the main results in the paper.
We consider the vanishing viscosity process (equations \eqref{2.1} and \eqref{4.1}) and aim at computing $u^{\epsilon_0}$ for $\epsilon_0>0$ sufficiently small, which is a good approximation of $u$, the viscosity solution to  \eqref{HJ}--\eqref{dir}.
The convexification is developed to compute this $u^{\epsilon_0}$.
Firstly, we obtain a new Carleman inequality in Theorem \ref{thm Car est}: 
For $\beta > 1$, $r > R + 1$, and $b > R + r$, we can find two numbers $\lambda_0 = \lambda_0(\beta, r, R, b, d)>0$, $C = C(r, R, b, d)>0$ such that for all $\lambda > \lambda_0$
and  for all $u \in C^2(\overline \Omega)$ with $u = 0$ on $\partial \Omega$ and $u_z = 0$ on $\Gamma^+$, we have
 \begin{multline*}
 \int_\Omega e^{2\lambda  (\frac{z + r}{b})^\beta} |\Delta u|^2 d\x
\geq C\lambda^3 \beta^2 (\beta - 1) b^{-3\beta} (-R + r)^{2\beta}  \int_{\Omega} e^{2\lambda  (\frac{z + r}{b})^\beta} |u|^2 d\x
\\
+ C \lambda (\beta - 1)b^{-\beta}\int_\Omega e^{2\lambda (\frac{z + r}{b})^\beta}|\nabla u|^2\,d\x.
\end{multline*}
We then use this Carleman estimate to show in Theorem \ref{thm convex} that the functional
\begin{equation}
	J_{\lambda, \beta, \eta}(u) = \int_{\Omega} e^{2\lambda (\frac{z + r}{b})^\beta} \big|-\epsilon_0 \Delta  u + F(\x,  u, \nabla  u)\big|^2 \,d\x 
	+ \eta \|u\|^2_{H^p(\Omega)}
	\label{Jdef}
\end{equation}
is strictly convex for $u\in H \cap \overline{B(M)}$.
Here, $p > \lceil d/2 \rceil + 2$ is such that $H^p \hookrightarrow C^2(\overline \Omega)$, and 
\[
H = \{u \in H^p(\Omega): u|_{\partial \Omega} = f \mbox{ and } u_z|_{\Gamma^+} = g\}, \quad
{B(M)} = \{u \in H^p(\Omega): \|u\|_{H^p(\Omega)} < M\}.
\]
Then, we use a gradient descent method (Theorem \ref{thm 4.2}) to compute the minimizer $u_{\rm min}$ of this functional. 
Assuming that $u_{\rm min} \in H \cap B(M/3)$, we can start the gradient descent method at $u^{(0)} \in H \cap B(M/3)$ and iterate
\begin{equation*}
	u^{(k)} = u^{(k - 1)} -  \kappa DJ_{\lambda, \beta, \eta}(u^{(k-1)}) \quad \text{ for } k \in \N.
\end{equation*}
Here, $ \kappa \in (0,  \kappa_0)$ where $ \kappa_0 \in (0, 1)$ depends only on $\lambda,$ $\beta$, $R,$ $r,$ $b$, $d$, $M$ and $\epsilon_0$.
We are able to obtain
\begin{equation*}
	\big\|u^{(k)} - u_{\rm min}\big\|_{H^p(\Omega)} \leq \theta^{k/2} \big\|u^{(0)} - u_{\rm min} \big\|_{H^p(\Omega)}  \quad \text{ for } k \in \N,
\end{equation*}
for some $\theta \in (0, 1)$ depending only on $ \kappa,$ $\lambda,$ $\beta$, $R,$ $r,$ $b$, $d$, $M$, $F$ and $\epsilon_0$.
Then, in Theorem \ref{thm 4.3}, we show that, even if there is a noise of size $\delta>0$, we still have a nice bound
\begin{equation*}
	\|u_{\min}^\delta - u^{\epsilon_0}\|^2_{H^1(\Omega)} \leq C(\eta \|u^{\epsilon_0}\|^2_{H^p(\Omega)} + \delta^2).
\end{equation*}
Here, $u^{\delta}_{\min}$ is the minimizer of $J_{\lambda, \beta, \eta}(u)$ with noisy data $u|_{\partial \Omega} = f^\delta$, $u_z|_{\Gamma^+} = g^\delta$.
Here, by saying that $\delta$ is the noise level,  there exists an ``error" function $\mathcal E$ satisfying
$\|\mathcal E\|_{H^p(\Omega)} \leq \delta$, $\mathcal E|_{\partial \Omega} = f^\delta - f$, $\mathcal E_z|_{\Gamma^+} = g^\delta - g$.
Combining Theorem \ref{thm 4.2} and Theorem \ref{thm 4.3}, we have for each $k \geq 1$,
\[
	\|u^{(k)} - u^{\epsilon_0}\|_{H^1(\Omega)} \leq C(\sqrt{\eta} \|u^{\epsilon_0}\|_{H^p(\Omega)} + \delta) + 
	\theta^{k/2} \|u^{(0)} - u^\delta_{\min}\|_{H^p(\Omega)}.
\]
This inequality shows the stability of our method with respect to noise. If $\theta^{k}$ and $\eta$ are $O(\delta^2)$ as $\delta$ tends to $0$, then the convergence rate is Lipschitz. 
Finally, in Section \ref{sec:num}, we implement the convexification method based on the finite difference method
and obtain interesting numerical results in two dimensions.

\medskip

We now address a bit further some state of the art  numerical methods in solving \eqref{HJ}--\eqref{dir} in the literature.
If $F$ is generically convex in $\nabla u$, there have been extremely powerful approaches to compute the solutions such as monotone numerical Hamiltonian based finite difference methods (see \cite{OsherSethian,Sethian, SethianVladimirsky, TsaiChengOsherZhao, QianZhangZhao} and the references therein).
When $F$ is nonconvex in $\nabla u$, the Lax--Friedrichs schemes (\cite{OsherShu, Abgrall, OsFe})  and the Lax--Friedrichs sweeping algorithm (\cite{KaoOsherQian,LiQian}), in which numerical viscosity terms appear naturally, are very efficient and accurate. 
Moreover, all the mentioned methods have very quick running times with not too many iterations.
Similar to the Lax--Friedrichs schemes, the addition of a viscosity term is natural in our approach as we deal with general $F$, which is possibly nonconvex in $\nabla u$.

\medskip

The paper is organized as follows.
In Section \ref{sec:prelim}, we give some preliminaries about viscosity solutions to Hamilton-Jacobi equations, which are rather well-known in the literature.
We state and prove a Carleman estimate in Theorem \ref{thm Car est} in Section \ref{sec:Carleman}.
Section \ref{sec conv} is devoted to the theoretical results of the convexification (Theorems  \ref{thm convex}--\ref{thm 4.3}), which is our main focus in this current paper.
Then, in Section \ref{sec:num}, we implement the convexification method based on the finite difference method
and obtain interesting numerical results in two dimensions.

\section{Some preliminaries about viscosity solutions to Hamilton-Jacobi equations} \label{sec:prelim}

\begin{Definition}[Viscosity solutions of \eqref{HJ}--\eqref{dir}]
Let $u \in C(\overline \Omega)$.
\begin{itemize}
\item[(a)] We say that $u$ is a viscosity subsolution to \eqref{HJ}--\eqref{dir} if for any test function $\varphi \in C^1(\overline \Omega)$ such that $u-\varphi$ has a strict maximum at $\x_0 \in \overline \Omega$, then
\[
F(\x_0, u(\x_0), \nabla \varphi(\x_0)) \leq 0 \quad \text{ if } \x_0 \in\Omega,
\]
or
\[
\min\left\{F(\x_0, u(\x_0),  \nabla\varphi(\x_0)), u(\x_0) - f(\x_0) \right\} \leq 0 \quad \text{ if } \x_0 \in\partial\Omega.
\]

\item[(b)] We say that $u$ is a viscosity supersolution to \eqref{HJ}--\eqref{dir} if for any test function $\varphi \in C^1(\overline \Omega)$ such that $u-\varphi$ has a strict minimum at $\x_0 \in \overline \Omega$, then
\[
F(\x_0, u(\x_0),  \nabla \varphi(\x_0)) \geq 0 \quad \text{ if } \x_0 \in\Omega,
\]
or
\[
\max\left\{F(\x_0, u(\x_0),  \nabla \varphi(\x_0)), u(\x_0) - f(\x_0) \right\} \geq 0 \quad \text{ if } \x_0 \in\partial\Omega.
\]

\item[(c)] We say that $u$ is a viscosity solution to \eqref{HJ}--\eqref{dir} if it is both its viscosity subsolution and its viscosity supersolution.
\end{itemize}
\end{Definition}
It is worth noting that the Dirichlet boundary condition might not hold in the classical sense for viscosity solutions to \eqref{HJ}--\eqref{dir} (see \cite[Appendix E]{Tran19} for example).
In the following, we impose some compatibility conditions to make sure that $u=f$ on $\partial \Omega$ classically.

\medskip

We write $F = F(\x, s, \p)$. Let $\nabla_{\x} F$ and $\nabla_{\p} F$  the gradient vector of $F$ with respect to the first and the third variables, respectively, and $\partial_s F=F_s$ the partial derivative of $F$ with respect to the second variable.
Here is one set of conditions that often occurs in the literature.

\begin{Assumption}
1. There exists $\alpha>0$ such that
\[
\alpha \leq F_s(\x,s,\p) \leq \frac{1}{\alpha} \quad \text{ for all } (\x,s,\p) \in \overline \Omega \times \R \times \R^d.
\]
2. $F$ is coercive in $\p$ in the sense that
\[
\lim_{|\p| \to \infty} \inf_{(\x,s) \in \overline \Omega \times \R} \left ( F(\x,s,\p)^2 + \partial_s F(\x,s,\p) |\p|^2 + \nabla_{\x}F(\x,s,\p)\cdot \p \right) = +\infty.
\]
3. There exists $\phi \in C^2(\overline \Omega)$ such that $\phi = f$ on $\partial \Omega$, and
\[
F(\x,\phi(\x),\nabla \phi(\x)) < 0 \quad \text{ in } \Omega.
\]
	\label{assumption 2}
\end{Assumption}

\begin{Theorem}
Suppose that Assumption \ref{assumption 2} holds.
Then, problem \eqref{HJ}--\eqref{dir} has a unique viscosity solution, denoted by $u$.
	Moreover, for $\epsilon > 0,$ the following problem 
	\begin{equation}
		\left\{
			\begin{array}{rcll}
				-\epsilon \Delta u^{\epsilon} + F(\x, u^\epsilon, \nabla u^\epsilon) &=& 0 &\x \in \Omega,\\
				u^{\epsilon} &=& f &\x \in \partial\Omega.
			\end{array}
		\right.
		\label{2.1}
	\end{equation}
	has a unique solution in $u^\epsilon \in C^2(\overline \Omega)$
and
	\begin{equation}
		\lim_{\epsilon \to 0^+} \|u^\epsilon - u\|_{L^\infty(\Omega)} = 0.
		\label{vanishing epsilon}
	\end{equation}
	\label{Theorem vis}
\end{Theorem}

Theorem \ref{Theorem vis} holds true under other appropriate sets of assumptions too.
We here just give one prototypical set of conditions, Assumption \ref{assumption 2}, for demonstration.
See  \cite{CrandallLions83, CrandallEvansLions84, Lions, Tran19} for its proof.
For other set of appropriate conditions on $F$, see references listed in Section \ref{sec:num} in each numerical test.

\begin{remark}
In fact, under Assumption \ref{assumption 2}, we are able to quantify the rate of convergence of $u^\epsilon$ to $u$ in $L^\infty(\Omega)$. There exists $C>0$ independent of $\epsilon \in (0,1)$ such that
\[
 \|u^\epsilon - u\|_{L^\infty(\Omega)}  \leq C \epsilon^{1/2}.
\]
See \cite{CL-rate, CagnettiGomesTran-2, Tran19} for more details.
\end{remark}

\section{A Carleman estimate} \label{sec:Carleman}

We prove a Carleman estimate that plays an important role in our proof for the convergence of the convexification method. 
In the proof of the Carleman estimate, we will need the notation
\begin{equation}
		\Gamma^- = \big\{
			\x = (x_1, x_2, \dots, x_{d - 1}, z = -R): |x_i| \leq R, 1 \leq i \leq d - 1
		\big\}
		\subset \partial \Omega.
\end{equation}

\begin{Theorem}[Carleman estimate]
For $\beta > 1$, $r > R + 1$, and $b > R + r$, we can find a number $\lambda_0 = \lambda_0(\beta, r, R, b, d)$ such that for all $\lambda > \lambda_0$
and 
 for all $u \in C^2(\overline \Omega)$ with $u = 0$  on $\partial \Omega$ and $u_z = 0$ on $\Gamma^+$, we have
 \begin{multline}
 \int_\Omega e^{2\lambda  (\frac{z + r}{b})^\beta} |\Delta u|^2 d\x
\geq C\lambda^3 \beta^2 (\beta - 1) b^{-3\beta} (-R + r)^{2\beta}  \int_{\Omega} e^{2\lambda  (\frac{z + r}{b})^\beta} |u|^2 d\x
\\
+ C \lambda (\beta - 1)b^{-\beta}\int_\Omega e^{2\lambda (\frac{z + r}{b})^\beta}|\nabla u|^2\,d\x.
\label{Car est}
\end{multline}
where $C = C(r, R, b, d)$ is a constant. 
Here, $\lambda_0$ and $C$ depend only on the listed parameters.
\label{thm Car est}
\end{Theorem}

\proof

In the proof below, $C_1,$ $C_2$ and $C_3$ are constants depending only on $r$, $R$, $b$, and $d$. We split the proof into several steps.

\noindent{\it Step 1.}
Define the function
\begin{equation}
	v = e^{\lambda (\frac{z + r}{b})^\beta} u.
\end{equation}
For $\x = (x_1, \dots, x_{d- 1}, z)$, $i = 1, \dots, d - 1$, we have
\begin{align*}
	u &= e^{-\lambda (\frac{z + r}{b})^\beta} v,
	&u_{x_i} &=e^{-\lambda (\frac{z + r}{b})^\beta} v_{x_i},  \\
	u_{x_i x_i} &=e^{-\lambda (\frac{z + r}{b})^\beta} v_{x_i x_i},   
	&u_z &=  e^{-\lambda (\frac{z + r}{b})^\beta}\Big(-\lambda \beta b^{-\beta} (z + r)^{\beta-1} v + v_z\Big),\\
\end{align*}
and
\[
	u_{zz} = e^{-\lambda (\frac{z + r}{b})^\beta}\Big(\Big( \lambda^2 \beta^2 b^{-2\beta} (z + r)^{2\beta - 2} -\lambda \beta (\beta-1) b^{-\beta}(z + r)^{\beta - 2}\Big) v
	- 2\lambda \beta b^{-\beta} (z + r)^{\beta - 1} v_z + v_{zz} \Big).
\]
Therefore, using the inequality $(a - b + c)^2 \geq -2ba - 2bc$, we have
\begin{align}
	&e^{2\lambda  (\frac{z + r}{b})^\beta} |\Delta u|^2 \nonumber\\
	=\ & \Big[ \big(\lambda^2 \beta^2 b^{-2\beta} (z + r)^{2\beta - 2} -\lambda \beta (\beta-1) b^{-\beta}(z + r)^{\beta - 2}\big) v
	- 2\lambda \beta b^{-\beta} (z + r)^{\beta - 1} v_z + \Delta v \Big]^2	\nonumber
	\\
	\geq \ &- 4\lambda^2 \beta^2 b^{-2\beta} (z + r)^{2\beta - 3}  \big(\lambda \beta b^{-\beta} (z + r)^{\beta } -(\beta-1)\big) v_z v
	-4\lambda \beta b^{-\beta} (z + r)^{\beta - 1} v_z \Delta v.
	\label{3.3}
\end{align}
Dividing both sides of \eqref{3.3} by $2\lambda \beta b^{-\beta}(z + r)^{\beta - 1}$ and integrating the resulting equation, we have
\begin{equation}
	\int_{\Omega}\frac{e^{2\lambda  (\frac{z + r}{b})^\beta} |\Delta u|^2}{2\lambda \beta b^{-\beta} (z + r)^{\beta - 1}} \,d\x
	\geq I_1 + I_2,
	\label{step1}
\end{equation}
where
\begin{equation}
	I_1 = -2\lambda \beta\int_{\Omega}  b^{-\beta}(z + r)^{\beta - 2}   \big(\lambda \beta b^{-\beta} (z + r)^{\beta } -(\beta-1)\big) v_z v \,d\x
	\label{I1}
\end{equation}
and
\begin{equation}
	I_2 = -2\int_{\Omega}  v_z \Delta v\,d\x.
	\label{I2}
\end{equation}

\noindent{\it Step 2.} In this step, we estimate $I_2$.
Since $u = 0$ on $\partial \Omega$, $v = 0$ on $\partial \Omega$.
Since $u = u_z = 0$ on $\Gamma^+$, $v_z = 0$ on $\Gamma^+$.
 Hence, $v_z = 0$ on $\partial \Omega \setminus \Gamma^-$.
Using integration by parts, we have
\begin{align*}
	I_2 &= -2 \int_{\Omega} \Div( v_z \nabla v)\,d\x
	+ 2 \int_{\Omega} \nabla v_z\cdot \nabla v\,d\x 
	= -2\int_{\Gamma^+}   |v_z|^2\,d\sigma
	+ 2\int_{\Gamma^-}  |v_z|^2\,d\sigma
	+ \int_{\Omega}  (|\nabla v|^2)_z \,d\x 
	\\
	&
	=
	 2\int_{\Gamma^-}  |v_z|^2\,d\sigma
	+ \int_{\Gamma^+}  |\nabla v|^2 \,d\sigma
	- \int_{\Gamma^-}  |\nabla v|^2 \,d\sigma.
\end{align*}
Since $v_{x_i} = 0$ on $\Gamma^+ \cup \Gamma^-$, $i = 1, \dots, d - 1$, $|\nabla v| = |v_z|$ on this set.
Using the fact that $v_z = 0$ on $\partial \Omega \setminus \Gamma^-$ again, we have
\begin{equation}
	I_2= \int_{\Gamma^+}   |v_z|^2\,d\sigma
	+ \int_{\Gamma^-}  |v_z|^2\,d\sigma 
	\geq 
	0.
	\label{3.6}
\end{equation}

\noindent{\it Step 3.} We estimate $I_1.$ It follows from \eqref{I1} that
\begin{align*}
	I_1 &= -\lambda \beta\int_{\Omega}  b^{-\beta}(z + r)^{\beta - 2}   \big(\lambda \beta b^{-\beta} (z + r)^{\beta } - (\beta-1)\big) (|v|^2)_z \,d\x
	\\
	& =-\lambda \beta\int_{\Omega} \big[ b^{-\beta}(z + r)^{\beta - 2}   \big(\lambda \beta b^{-\beta} (z + r)^{\beta } -(\beta-1)\big) |v|^2\big]_z \,d\x
	\\
	&\hspace{6cm}
	+ \lambda \beta\int_{\Omega} \big[ b^{-\beta}(z + r)^{\beta - 2}   \big(\lambda \beta b^{-\beta} (z + r)^{\beta } -(\beta-1)\big)\big]_z |v|^2 \,d\x
	\\
	&= \lambda \beta \int_{\Omega} \big[2\lambda  \beta (\beta - 1) b^{-2\beta}(z + r)^{2\beta - 3}
	- (\beta-1)(\beta - 2) b^{-\beta}(z + r)^{\beta - 3}
	 \big]|v|^2\,d\x.
\end{align*}

Now, for fixed $\beta > 1$, we can find $\lambda_0 \geq 1$, only depending on $\beta$, $r$ and $R$, sufficiently large such that for all $\lambda > \lambda_0$, 
\begin{equation}
	I_1 \geq C_1\lambda^2 \beta^2 (\beta - 1) b^{-2\beta}   \int_{\Omega} (z + r)^{2\beta} |v|^2 d\x.
	\label{3.7}
\end{equation}

\noindent{\it Step 3}. 
Combining \eqref{step1}, \eqref{3.6} and \eqref{3.7}, we have
\begin{equation}
	\int_\Omega \frac{e^{2\lambda  (\frac{z + r}{b})^\beta} |\Delta u|^2}{2\lambda \beta b^{-\beta} (z + r)^{\beta - 1}} d\x\geq 
 C_1\lambda^2 \beta^2 (\beta - 1) b^{-2\beta}   \int_{\Omega} (z + r)^{2\beta} |v|^2 d\x.
	 \label{3.8}
\end{equation}

Recall that $v = e^{\lambda  (\frac{z + r}{b})^\beta} u$. 
It follows from \eqref{3.8} that
\begin{equation*}
	\int_\Omega \frac{e^{2\lambda  (\frac{z + r}{b})^\beta} |\Delta u|^2}{2\lambda \beta b^{-\beta} (z + r)^{\beta - 1}} d\x
	\geq 
 C_1\lambda^2 \beta^2 (\beta - 1) b^{-2\beta}   \int_{\Omega} (z + r)^{2\beta} e^{2\lambda  (\frac{z + r}{b})^\beta} |u|^2 d\x.
\end{equation*}
or equivalently,
\begin{equation}
	\int_\Omega \frac{e^{2\lambda  (\frac{z + r}{b})^\beta} |\Delta u|^2}{ \beta b^{-\beta} (z + r)^{\beta - 1}} d\x
	\geq 
 C_1\lambda^3 \beta^2 (\beta - 1) b^{-2\beta}   \int_{\Omega} (z + r)^{2\beta} e^{2\lambda  (\frac{z + r}{b})^\beta} |u|^2 d\x.
	 \label{3.9}
\end{equation}

\noindent {\it Step 4.}
We  estimate the term $\ds \int_{\Omega}  e^{2\lambda (\frac{z + r}{b})^\beta} u\Delta u \,d\x$.
Using integration by parts, since $u = 0$ on $\partial \Omega,$ we have
\begin{align*}
	-\int_{\Omega}  e^{2\lambda (\frac{z + r}{b})^\beta} u\Delta u \,d\x
	&=  -\int_{\Omega} \Div\big[  e^{2\lambda (\frac{z + r}{b})^\beta} u\nabla u\big] \,d\x
	+ \int_{\Omega} \nabla(e^{2\lambda (\frac{z + r}{b})^\beta} u) \cdot \nabla u \,d\x
	\\
	&\hspace{-3cm}=  \int_{\Omega} e^{2\lambda (\frac{z + r}{b})^\beta}\big[2\lambda \beta b^{-\beta}(z + r)^{\beta-1}u + u_z\big]u_z\,d\x
	 + \sum_{i = 1}^{d - 1}\int_\Omega e^{2\lambda (\frac{z + r}{b})^\beta}|u_{x_i}|^2\,d\x
	 \\
	&\hspace{-3cm}
	= \lambda \beta  \int_{\Omega} e^{2\lambda (\frac{z + r}{b})^\beta} b^{-\beta} (z + r)^{\beta-1} (|u|^2)_z\,d\x
	 + \int_\Omega e^{2\lambda (\frac{z + r}{b})^\beta}|\nabla u|^2\,d\x
	 \\
	&\hspace{-3cm}
	= \lambda \beta  \int_{\Omega} \Big[e^{2\lambda (\frac{z + r}{b})^\beta} b^{-\beta}(z + r)^{\beta-1} |u|^2\Big]_z\,d\x
	-\lambda \beta  \int_{\Omega} \Big[e^{2\lambda (\frac{z + r}{b})^\beta} b^{-\beta} (z + r)^{\beta-1}\Big]_z |u|^2\,d\x
	\\
	&\hspace{6cm}
	 + \int_\Omega e^{2\lambda (\frac{z + r}{b})^\beta b^{-\beta}}|\nabla u|^2\,d\x.
\end{align*}
Due to the fact that $u = 0$ on $\partial \Omega$, the first integral in the last row above vanishes.
We have
\begin{multline*}
	-\int_{\Omega}  e^{2\lambda (\frac{z + r}{b})^\beta} u\Delta u \,d\x = 
-\lambda \beta  \int_{\Omega} e^{2\lambda (\frac{z + r}{b})^\beta} \Big(\lambda \beta b^{-2\beta} (z + r)^{2\beta - 2}  + (\beta - 1)b^{-\beta} (z + r)^{-2+\beta}\Big)  |u|^2\,d\x
\\
+  \int_\Omega e^{2\lambda (\frac{z + r}{b})^\beta b^{-\beta}}|\nabla u|^2\,d\x.
\end{multline*}
Hence, for $\lambda > \lambda_0$,
\begin{equation}
	-\int_{\Omega}  e^{2\lambda (\frac{z + r}{b})^\beta} u\Delta u \,d\x
	\geq
	 -C_2\lambda^2 \beta^2 b^{-2\beta}  \int_{\Omega} 	
		(z + r)^{2\beta}  e^{2\lambda (\frac{z + r}{b})^\beta} 
	 |u|^2\,d\x
	 + \int_\Omega e^{2\lambda (\frac{z + r}{b})^\beta}|\nabla u|^2\,d\x.
	 \label{3.12}
\end{equation}

\noindent{\it Step 5.} We complete the proof in this step. Using the inequality $x^2 + y^2 \geq -2xy$, we have
\begin{multline}
	\lambda^{5/2} \beta b^{-\beta} \int_{\Omega} (z + r)^{\beta-1} e^{2\lambda (\frac{z + r}{b})^\beta} |u|^2 d\x
	+ \int_\Omega \frac{e^{2\lambda  (\frac{z + r}{b})^\beta} |\Delta u|^2}{4\lambda^{1/2} \beta b^{-\beta} (z + r)^{\beta - 1}} d\x
	\\
	\geq -\lambda\int_{\Omega}  e^{2\lambda (\frac{z + r}{b})^\beta} u\Delta u \,d\x.
	\label{3.1222}
\end{multline}

Combining \eqref{3.12} and \eqref{3.1222}, we have
\begin{multline*}
	\int_\Omega \frac{e^{2\lambda  (\frac{z + r}{b})^\beta} |\Delta u|^2}{4\lambda^{1/2} \beta b^{-\beta} (z + r)^{\beta - 1}} d\x
	\geq
	- \lambda^{5/2} \beta b^{-\beta} \int_{\Omega} (z + r)^{\beta-1} e^{2\lambda (\frac{z + r}{b})^\beta} |u|^2 d\x
	\\
	  -C_2\lambda^3 \beta^2 b^{-2\beta}  \int_{\Omega} 	
		(z + r)^{2\beta}  e^{2\lambda (\frac{z + r}{b})^\beta} 
	 |u|^2\,d\x
	 +  \lambda  \int_\Omega e^{2\lambda (\frac{z + r}{b})^\beta}|\nabla u|^2\,d\x,
\end{multline*}
which implies
\begin{multline}
	\int_\Omega \frac{e^{2\lambda  (\frac{z + r}{b})^\beta} |\Delta u|^2}{4\lambda^{1/2} \beta b^{-\beta} (z + r)^{\beta - 1}} d\x
	\geq
	  -C_3\lambda^3 \beta^2 b^{-2\beta}  \int_{\Omega} 	
		(z + r)^{2\beta}  e^{2\lambda (\frac{z + r}{b})^\beta} 
	 |u|^2\,d\x
	 \\
	 +  \lambda  \int_\Omega e^{2\lambda (\frac{z + r}{b})^\beta}|\nabla u|^2\,d\x.
	 \label{3.1414}
\end{multline}
Multiply $\frac{C_1(\beta - 1)}{2C_3}$ to both sides of \eqref{3.1414} where $C_1$ is the constant in \eqref{3.9}, we yield
\begin{multline}
	\frac{C_1(\beta - 1)}{2C_3}\int_\Omega \frac{e^{2\lambda  (\frac{z + r}{b})^\beta} |\Delta u|^2}{4\lambda^{1/2} \beta b^{-\beta} (z + r)^{\beta - 1}} d\x
	\geq
	  - \frac{C_1\lambda^3 \beta^2(\beta - 1) b^{-2\beta} }{2}   \int_{\Omega} 	
		(z + r)^{2\beta}  e^{2\lambda (\frac{z + r}{b})^\beta} 
	 |u|^2\,d\x
	 \\
	 +  \frac{C_1  \lambda (\beta - 1)}{2C_3}  \int_\Omega e^{2\lambda (\frac{z + r}{b})^\beta}|\nabla u|^2\,d\x.
	 \label{3.13}
\end{multline}

Adding \eqref{3.9} and \eqref{3.13},
we have
\begin{multline}
\Big(1 + \frac{C_1(\beta - 1)}{2C_3\lambda^{1/2}}\Big) \int_\Omega \frac{e^{2\lambda  (\frac{z + r}{b})^\beta} |\Delta u|^2}{ \beta b^{-\beta} (z + r)^{\beta - 1}} d\x
\geq \frac{C_1\lambda^3 \beta^2(\beta - 1) b^{-2\beta} }{2}   \int_{\Omega} 	
		(z + r)^{2\beta}  e^{2\lambda (\frac{z + r}{b})^\beta} 
	 |u|^2\,d\x
	 \\
	 +  \frac{C_1  \lambda (\beta - 1)}{2C_3}  \int_\Omega e^{2\lambda (\frac{z + r}{b})^\beta}|\nabla u|^2\,d\x.
\label{3.14}
\end{multline}
The desired Carleman estimate \eqref{Car est} follows from  \eqref{3.14}.
\QEDB

\begin{remark}
	In our previous publications, see e.g. \cite{NguyenLiKlibanov:IPI2019, LeNguyen:2020}, the number $\beta$ must be large. The reason for \eqref{Car est} holds true when $\beta > 1$ is our trick of dividing both sides of  \eqref{3.3} by $2\lambda \beta b^{-\beta}(z + r)^{\beta - 1}$  so that the corresponding $I_2$, defined in \eqref{I2}, is nonnegative.	
	This trick was used in \cite[Theorem 3.1]{LeNguyen:2020}. However, since the Carleman weight in  \cite[Theorem 3.1]{LeNguyen:2020} is different from the one in \eqref{Car est}, the parameter $\beta$ in \cite[Theorem 3.1]{LeNguyen:2020} must be large.
	In the case when the principal differential operator in \eqref{Car est} is replaced by the general elliptic operator, this trick is not applicable.   
\end{remark}

\section{The convexification method to compute the viscosity solution to Hamilton-Jacobi equations} \label{sec conv}

In this section, we propose to use the convexification method to solve \eqref{HJ} together with the boundary condition \eqref{dir} supposing that Assumption \ref{assumption 1} holds true. 
That means we know the boundary value $u = f$ on $\partial \Omega$ and the function $u_z|_{\Gamma^+}$ can be computed, say $u_z = g$ on $\Gamma^+.$
Due to Theorem \ref{Theorem vis}, it is natural to try to approximate the solution $u$ by a function $u^{\epsilon_0}$ that satisfies
\begin{equation}
	\left\{
	\begin{array}{rcll}
		-\epsilon_0 \Delta u^{\epsilon_0} + F(\x, u^{\epsilon_0}, \nabla u^{\epsilon_0}) &=& 0 &\x \in \Omega,\\
		u^{\epsilon_0} &=& f(\x) &\x \in \partial\Omega,\\
		u^{\epsilon_0}_z &=& g(\x) &\x \in \Gamma^+.
	\end{array}
	\right.
	\label{4.1}
\end{equation}
\begin{remark}
1. In general, \eqref{4.1} is over-determined.
It might have no solution; especially when the boundary data contains noise.
However, the convexification method can deliver a function that ``most fits" \eqref{4.1} and show that this function is an approximation of the true solution to \eqref{HJ}--\eqref{dir} when the given boundary data are noiseless. 
Again, due to Assumption \ref{assumption 1}, the Neumann condition imposed in \eqref{4.1} makes sense.

2. On the other hand, \eqref{4.1} is not over-determined in the  sense that
\begin{equation}
	\left\{
	\begin{array}{rcll}
		-\epsilon_0 \Delta u^{\epsilon_0} + F(\x, u^{\epsilon_0}, \nabla u^{\epsilon_0}) &=& 0 &\x \in \Omega,\\
		u^{\epsilon_0} &=& f(\x) &\x \in \partial\Omega,\\
	\end{array}
	\right.
	\label{4.2222}
\end{equation}
might not be uniquely solvable if we do not impose Assumption \ref{assumption 2} or other sets of appropriate conditions. 
For example, when $F(\x, u^{\epsilon_0}, \nabla u^{\epsilon_0}) = -u^{\epsilon_0}$, $f=0$, and $\frac{1}{\epsilon_0}$ is an eigenvalue of $-\Delta$,  \eqref{4.2222} has multiple solutions.
Therefore, imposing the additional Neumann boundary condition on $\Gamma^+$, in the general case, is necessary.
\end{remark}
Let $p > \lceil d/2 \rceil + 2$ such that $H^p(\Omega) \hookrightarrow C^2(\overline \Omega).$
Define
\[
	H = \{
		u \in H^p(\Omega): u|_{\partial \Omega} = f \mbox{ and } u_z|_{\Gamma^+} = g
	\}.
\]
Clearly, $H$ is a closed subset of $H^p(\Omega).$ 
We will also need the following set of test functions
\[
	H_0 = \{
		u \in H^p(\Omega): u|_{\partial \Omega} = 0 \mbox{ and } u_z|_{\Gamma^+} = 0
	\}.
\]
Let $M > 0$ be chosen later, and set 
\[
	{B(M)} = \{u \in H^p(\Omega): \|u\|_{H^p(\Omega)} < M\}.
\]
We assume that $H \cap \overline{B(M)} \not = \emptyset$. 
For each $\lambda > 1$, $\beta > 1$ and $\eta > 0$, introduce the functional
\begin{equation}
	J_{\lambda, \beta, \eta}(u) = \int_{\Omega} e^{2\lambda (\frac{z + r}{b})^\beta} \big|-\epsilon_0 \Delta  u + F(\x,  u, \nabla  u)\big|^2 \,d\x 
	+ \eta \|u\|^2_{H^p(\Omega)}
	\quad u \in H \cap \overline{B(M)}
	\label{Jdef}
\end{equation}
where $\epsilon_0$ is a fixed small positive number.

\medskip

The convexification theorem is to prove that for each $\beta > 1$, there is a number $\lambda_0 > 1$ such that for all $\lambda \geq \lambda_0$ and for all $\eta > 0$, the function $J_{\lambda, \beta, \eta}$ is strictly convex in $H \cap \overline{B(M)}.$ 
The word ``convexification" is suggested by the fact that the Carleman weight function $e^{2\lambda (\frac{z + r}{b})^\beta}$ convexifies this functional.
The convexification theorem is stated below.
\begin{Theorem}[The convexification theorem]
1. For all $\lambda, \beta >1$, and $\eta > 0$, for all $u \in H^p(\Omega)$, $h\in H_0$, we have	\begin{equation}
		\lim_{H_0 \in h \to 0} \frac{|J_{\lambda, \beta, \eta}(u + h) - J_{\lambda, \beta, \eta}(u) -
		DJ_{\lambda, \beta, \eta}(u)h |}{\|h\|_{H^p(\Omega)}} = 0
		\label{4.2}
	\end{equation}
	where 
	\begin{multline}
		DJ_{\lambda, \beta, \eta}(u)h = 
		 2\int_{\Omega} e^{2\lambda (\frac{z + r}{b})^\beta}\Big[-\epsilon_0 \Delta u +  F(\x, u, \nabla u)\Big]
		  \Big[-\epsilon_0 \Delta h + \partial_s F(\x, u, \nabla u) h
		  \\
		   + \nabla_{\p} F(\x, u, \nabla u) \cdot \nabla h\Big]\,d\x 
	+ 2\eta\langle u, h\rangle_{H^p(\Omega)}.
	\label{DJ}
	\end{multline}

2.	Let $M$ be an arbitrarily large number. 
For each $\beta> 1$,
  $\lambda >\lambda_0 = \lambda_0(\epsilon_0, M, b, d, r, F, \beta) > 1$, $\eta > 0$, $u \in H\cap \overline{B(M)}$ and $v \in H \cap \overline{B(M)}$, we have
 \begin{equation}
	J_{\lambda, \beta, \eta}(v) - J_{\lambda, \beta, \eta}(u) - DJ_{\lambda, \beta, \eta}(u)h \geq C \|v - u\|^2_{H^1(\Omega)} + \eta \|v - u\|^2_{H^p(\Omega)},
	\label{convex}
\end{equation}
and
\begin{equation}
	\langle DJ_{\lambda, \beta, \eta}(v) - DJ_{\lambda, \beta, \eta}(u),(v - u)\rangle_{H^p(\Omega)}  \geq C \|v - u\|^2_{H^1(\Omega)} + \eta \|v - u\|^2_{H^p(\Omega)}.
	\label{convex1}
\end{equation}
Here, the constant $C$ depends only on $\lambda,$ $\beta$, $R,$ $r,$ $b$, $d$, $M$, $F$ and $\epsilon_0$.
As a result, the functional $J_{\lambda, \beta, \eta}$ has a unique minimizer in $\overline{B(M)}$.
\label{thm convex}
\end{Theorem}

The key point for us to successfully establish the inequalities \eqref{convex} and \eqref{convex1} is the presence of the Carleman weight function $e^{2\lambda (\frac{z + r}{b})^\beta}$ and the use of the Carleman estimate \eqref{Car est}.
A direct consequence of the inequalities \eqref{convex} and \eqref{convex1} is the strict convexity of $J_{\lambda, \beta, \eta}$ in $H \cap \overline{B(M)}.$
It is worth mentioning that $C$ in  \eqref{convex}--\eqref{convex1} depends on the viscosity coefficient $\epsilon_0$, and we need $\lambda >\lambda_0 = \lambda_0(\epsilon_0, M, b, d, r, F, \beta) > 1$ so that  $J_{\lambda, \beta, \eta}$ is convex in $H \cap \overline{B(M)}.$
On the other hand,  $D^2 J_{\lambda, \beta, \eta} \geq \eta Id$ in $H \cap \overline{B(M)}$, which means that the uniform convexity of $J_{\lambda, \beta, \eta}$ only depends on $\eta$, not $\epsilon_0$.
The proof of Theorem \ref{thm convex} is similar to that of the convexification theorem in \cite{KlibanovNik:ra2017}, which is originally designed to solve highly nonlinear and severely ill-posed inverse problems. Since this is the first time the convexification method is employed in the area of numerical methods for Hamilton-Jacobi equations, we present the proof here for the reader's convenience. 

\begin{remark}
For all $u \in H^p(\Omega)$, since $DJ_{\lambda, \beta, \eta}(u)$ is a bounded linear map from $H_0$ into $\R$, by the Riesz theorem there exists uniquely the function $J'_{\lambda, \beta, \nu}(u) \in H_0$ such that
\begin{equation*}
	\langle J'_{\lambda, \beta, \nu}(u), h \rangle_{H^p(\Omega)} =  DJ_{\lambda, \beta, \eta}(u) h
\end{equation*}
for all $h \in H_0.$
\label{rem 4.2}
\end{remark}

\proof[Proof of Theorem \ref{thm convex}]
	Since $F$ is in the class $C^2(\overline \Omega \times \R \times \R^d)$, for all $u, h \in H^p(\Omega)$ we can write
	\begin{equation}
		F(\x, u + h, \nabla (u + h))  = F(\x, u, \nabla u) + \partial_s F(\x, u, \nabla u) h + \nabla_{\p} F(\x, u, \nabla u) \cdot \nabla h
		+ O(|h|^2) + O(|\nabla h|^2).
		\label{4.4444}
	\end{equation}
	Here, $O(s)$ is the quantity satisfying $|O(s)| \leq C|s|$ where $C$ is a  constant that might depend on an upper bound of  $\|u\|_{C^1(\overline \Omega)}$  and the function $F$.
	Using \eqref{Jdef} and \eqref{4.4444}, we obtain
	\begin{align*}
		J_{\lambda, \beta, \eta}(&u + h) - J_{\lambda, \beta, \eta}(u) -  2\eta \langle u, h\rangle_{H^p(\Omega)} - \eta\|h\|^2_{H^p(\Omega)}
		\\
		&= \int_{\Omega} e^{2\lambda(\frac{z + r}{b})^\beta} \Big(|-\epsilon_0 \Delta(u+h) + F(\x, u + h, \nabla (u + h))|^2 - |-\epsilon_0 \Delta u + F(\x, u, \nabla u)|^2\Big)\,d\x 
		\\
		&= \int_{\Omega} e^{2\lambda(\frac{z + r}{b})^\beta} \big(-\epsilon_0 \Delta h + F(\x, u + h, \nabla (u + h)) - F(\x, u, \nabla u)\big)\big(-\epsilon_0 \Delta h - 2 \epsilon_0 \Delta u 
		\\
		&\hspace{8cm}
		+ F(\x, u + h, \nabla (u + h)) + F(\x, u, \nabla u)\big)\,d\x				\\
		&= \int_{\Omega} e^{2\lambda(\frac{z + r}{b})^\beta} \Big[-\epsilon_0 \Delta h + \partial_s F(\x, u, \nabla u) h + \nabla_{\p} F(\x, u, \nabla u) \cdot \nabla h + O(|h|^2) + O(|\nabla h|^2)\Big]
		\\ &\hspace{1cm}  \Big[
		2\big(-\epsilon_0 \Delta u + F(\x, u, \nabla u)\big)
		-\epsilon_0 \Delta h
		+ \partial_s F(\x, u, \nabla u) h + \nabla_{\p} F(\x, u, \nabla u) \cdot \nabla h 
		\\
		&\hspace{10cm}+ O(|h|^2) + O(|\nabla h|^2)
		\Big]\,d\x.
	\end{align*}
Hence,
\begin{multline}
	J_{\lambda, \beta, \eta}(u + h) - J_{\lambda, \beta, \eta}(u) - 2\eta\langle u, h\rangle_{H^p(\Omega)}  - \eta\|h\|^2_{H^p(\Omega)}
	- 2\int_{\Omega} e^{2\lambda (\frac{z + r}{b})^\beta} \Big[ -\epsilon_0 \Delta u + F(\x, u, \nabla u) \Big]
	\\
	\Big[-\epsilon_0 \Delta h
	 + \partial_s F(\x, u, \nabla u) h 
	+ \nabla_{\p} F(\x, u, \nabla u) \cdot \nabla h\Big]\,d\x 
	\\
	=2\int_{\Omega} e^{2\lambda (\frac{z + r}{b})^\beta}  \Big[\epsilon_0^2|\Delta h|^2 +  O(|h|^2) + O(|\nabla h|^2)\Big]\,d\x.
	\label{4.4}
\end{multline}
Defining $DJ_{\lambda, \beta, \eta}$ as in \eqref{DJ} and using \eqref{4.4}, we have
\begin{multline*}
	\lim_{H_0 \ni h \to 0} \frac{|J_{\lambda, \beta, \eta}(u + h) - J_{\lambda, \beta, \eta}(u) - DJ_{\lambda, \beta, \eta}(u)h|}{\|h\|_{H^p(\Omega)}} 
	\\
	\leq 
	\lim_{H_0 \ni h \to 0} \frac{|J_{\lambda, \beta, \eta}(u + h) - J_{\lambda, \beta, \eta}(u) - DJ_{\lambda, \beta, \eta}(u)h|}{\|h\|_{H^2(\Omega)}} = 0.
\end{multline*}
We have proved part 1 of this theorem. 
We next prove part 2. 

\medskip

For any $u$ and $v$ in $H \cap \overline{B(M)}$, set $h = v - u \in H_0$. It follows from \eqref{4.4} that
\begin{align}
	J_{\lambda, \beta, \eta}(v) - J_{\lambda, \beta, \eta}(u)  -  & DJ_{\lambda, \beta, \eta}(u)(v - u)
	 = J_{\lambda, \beta, \eta}(u + h) - J_{\lambda, \beta, \eta}(u) - DJ_{\lambda, \beta, \eta}(u)h \nonumber\\
	&\geq \int_{\Omega}e^{2\lambda (\frac{z + r}{b})^\beta}\Big[ \epsilon_0^2|\Delta h|^2 -  O(|h|^2) - O(|\nabla h|^2)\Big]\,d\x + \eta \|h\|_{H^p(\Omega)}^2. \label{4.6}
\end{align}
Note that $h|_{\partial \Omega}=0$ and $h_z|_{\Gamma^+} = 0$.
Applying Theorem \ref{thm Car est}, for each $\beta > 1$, we can find $\lambda_0$ depending on $R$, $M$, $F$, $\beta,$  $b$ and $d$  such that for all $\lambda > \lambda_0$,
\begin{equation}
	\int_{\Omega}e^{2\lambda (\frac{z + r}{b})^\beta} |\Delta h|^2\,d\x \geq C\lambda \int_{\Omega}e^{2\lambda (\frac{z + r}{b})^\beta}( |\nabla h|^2 + \lambda^2 |h|^2)\,d\x.
	\label{4.7}
\end{equation}
Here, the constant $C$ is allowed to depend on $\beta.$
We now choose $\lambda$ such that $\epsilon_0^2 \lambda$ is sufficiently large. 
Combining \eqref{4.6} and \eqref{4.7}, we have 
\begin{equation*}
	J_{\lambda, \beta, \eta}(v) - J_{\lambda, \beta, \eta}(u) -  DJ_{\lambda, \beta, \eta}(u)(v - u)\geq C \|h\|^2_{H^1(\Omega)} + \eta \|h\|^2_{H^p(\Omega)}.
\end{equation*}
We have proved \eqref{convex}.
Here, $C$ depends on $\lambda,$ $\beta,$ $\epsilon_0$, $R$, $M$, $b$ and $d$. 
Interchanging the roles of $u$ and $v$ in \eqref{convex} and adding the resulting estimate to \eqref{convex}, we obtain \eqref{convex1}. 

\medskip

We next prove that if \eqref{convex1} holds true, $J_{\lambda, \beta, \eta}$ has a unique minimizer.
The existence of the minimizer is obvious. It follows from the fact that $H \cap \overline{B(M)}$ is convex in $H^p(\Omega)$ and the compact embedding of $H^p(\Omega)$ to $H^2(\Omega).$
An alternative way to obtain the existence of the minimizer is to argue similarly to the proofs of Theorem 2.1 in \cite{KlibanovNik:ra2017} or Theorem 4.1 in \cite{KlibanovLeNguyenARL:preprint2021}.
We now prove the uniqueness. 
Let $u_1$ and $u_2$ be two local minimizers of $J_{\lambda, \beta, \eta}$ on $H \cap \overline{B(M)}$. 
It is clear that, see \cite[Lemma 2]{KlibanovNik:ra2017},
\[
	 DJ_{\lambda, \beta, \eta}(u_1)(u_2 - u_1) \geq 0 \quad
	\mbox{and} \quad
	 DJ_{\lambda, \beta, \eta}(u_2)(u_1 - u_2) \geq 0.
\]
Thus,
\begin{equation}
	 (DJ_{\lambda, \beta, \eta}(u_1) - DJ_{\lambda, \beta, \eta}(u_2))(u_1 - u_2) \leq 0 
	\label{4.12}
\end{equation}
Combining \eqref{convex1} for $u_1$ and $u_2$ and \eqref{4.12}, we have
\[
	C\|u_1 - u_2\|^2_{H^1(\Omega)} + \eta \|u_1 - u_2\|^2_{H^p(\Omega)}  \leq 0.
\]
The proof is complete.
\QEDB

\begin{Theorem}[The convergence of the gradient descent method]
 Let $\lambda, \beta$ and $\eta$ as in part 2 of Theorem \ref{thm convex} and let $J'_{\lambda, \beta, \eta}$ be as in Remark \ref{rem 4.2}. 
Let $u^{(0)}$ be any function in $H \cap B(M/3)$.
For $k \in \N$, define
\begin{equation}
	u^{(k)} = u^{(k - 1)} - \kappa J_{\lambda, \beta, \eta}'(u^{(k-1)})
	\label{gdm}
\end{equation}
for all $ \kappa \in (0,  \kappa_0)$ where $ \kappa_0 \in (0, 1)$ is a number that depends only on $\lambda,$ $\beta$, $R,$ $r,$ $b$, $d$, $M$ and $\epsilon_0$.
Then, if the minimizer $u_{\rm min}$ of $J_{\lambda, \beta, \eta}$ is in $H \cap B(M/3)$ then
there is a number $\theta \in (0, 1)$ depending only on $ \kappa,$ $\lambda,$ $\beta$, $R,$ $r,$ $b$, $d$, $M$, $F$ and $\epsilon_0$ such that for all $k \geq 1$
\begin{equation}
	\big\|u^{(k)} - u_{\rm min}\big\|_{H^p(\Omega)} \leq \theta^{k/2} \big\|u^{(0)} - u_{\rm min} \big\|_{H^p(\Omega)}.
	\label{4.1313}
\end{equation}
\label{thm 4.2}
\end{Theorem}


Theorem \ref{thm 4.2} and the estimate \eqref{4.1313} guarantee that the minimizer of $J_{\lambda, \beta, \eta}$ can be found by the popular gradient descent method. 
The success is due to the hypothesis that the desired minimizer is in the interior of $B(M/3)$. 
We do not experience any difficulty due to not checking this condition in the numerical study.
However, to be more rigorous,
in general, if this condition cannot be verified, one can use the projected gradient method as in \cite{KlibanovNik:ra2017, KhoaKlibanovLoc:SIAMImaging2020} to find the minimizer. In this paper, we choose the gradient descent method because the implementation of the projection in the projected gradient method is more complicated while there are many ready-to-use packages; for e.g., the optimization toolbox of Matlab, for the gradient descent method.
In other words, Theorem \ref{thm 4.2}  significantly reduces our efforts in implementation. 
Although the proof of Theorem \ref{thm 4.2} is similar to the proofs of \cite[Theorem 6]{Klibanov:2ndSAR2021} in 1D case and \cite[Theorem 2.2]{LeNguyen:preprint2021} in higher dimensions, we briefly present the proof of Theorem \ref{thm 4.2} here for the convenience of the reader.

\proof[Proof of Theorem \ref{thm 4.2}]
Since $F$ is a function in the class $C^2(\overline \Omega \times \R \times \R^d)$, it is obvious that $\nabla_{\x}F$, $\partial_s F$, $\nabla_{\p} F$ are all Lipschitz continuous in any bounded subdomain of $\overline \Omega \times \R \times \R^{d}$. 
As a result, $DJ_{\lambda, \beta, \eta}$, see \eqref{DJ} for its definition, is Lipschitz continuous on $H \cap \overline{B(M)}$. 
Hence, there is a positive number $L$ such that
\begin{equation}
	\|J_{\lambda, \beta, \eta}'(u_1) - J_{\lambda, \beta, \eta}'(u_2)\|_{H^p(\Omega)} \leq L \|u_1 - u_2\|_{H^p(\Omega)}
	\label{4.13}
\end{equation}
for all $u_1, u_2 \in H \cap B(M).$
We claim that for any $k \geq 0$, 
\begin{equation}
	\|u^{(k)}  - u_{\rm min}\|_{H^p(\Omega)} < \frac{2M}{3}.
	\label{bbb}
\end{equation}
 This is true when $k = 0$. Assume \eqref{bbb} is true for some $k$, we will prove that \eqref{bbb} holds true for $k + 1.$ Due to \eqref{bbb}, $u^{(k)} \in H \cap B(M)$.
Since $u_{\rm min}$ is in $H \cap B(M/3)$, $J_{\lambda, \beta, \eta}'(u_{\rm min}) = 0$.
Using \eqref{convex1} and \eqref{4.13}, we have
\begin{align*}
	\big\|u^{(k+1)} &- u_{\rm min}\big\|_{H^p(\Omega)}^2 
	= \big\|u^{(k)} - u_{\rm min} -  \kappa [J_{\lambda, \beta, \eta}'(u^{(k)})	- J_{\lambda, \beta, \eta}'(u_{\rm min})] \big\|_{H^p(\Omega)}^2
	\\
	&= \big\|u^{(k )} - u_{\rm min}\big\|_{H^p(\Omega)}^2 
	- 2 \kappa \big\langle J_{\lambda, \beta, \eta}'(u^{(k)})	- J_{\lambda, \beta, \eta}'(u_{\rm min}), u^{(k )} - u_{\rm min}\big\rangle_{H^p(\Omega)}
	\\
	&\hspace{7cm}
	+  \kappa^2\big\|J_{\lambda, \beta, \eta}'(u^{(k)})	- J_{\lambda, \beta, \eta}'(u_{\rm min})\big\|_{H^p(\Omega)}^2
	\\
	&\leq (1 - 2 \kappa\eta +  \kappa^2 L^2)\big\|u^{(k)} - u_{\rm min} \big\|_{H^p(\Omega)}^2.
\end{align*}
Choosing $ \kappa \in (0,  \kappa_0)$ where $ \kappa_0 = \eta L^{-2}$, we have $\theta = 1 - 2 \kappa\eta +  \kappa^2 L^2 \in (0, 1)$ and
\begin{equation}
	\big\|u^{(k+1)} - u_{\rm min}\big\|_{H^p(\Omega)} \leq \theta^{1/2} \big\|u^{(k)} - u_{\rm min} \big\|_{H^p(\Omega)} < \big\|u^{(k)} - u_{\rm min} \big\|_{H^p(\Omega)} < \frac{2M}{3}.
	\label{4.15}
\end{equation}
We have proved \eqref{bbb}.  The estimate \eqref{4.1313}  follows from \eqref{4.15}.
\QEDB

\medskip

Consider the case when the boundary data $f = u|_{\partial \Omega}$ contains noise with the noise level $\delta > 0$.
Since  the knowledge of $g = u_z|_{\Gamma^+}$ can be computed from the knowledge of $u$ on this set, see Assumption \ref{assumption 1} and Remark \ref{rm 1}, the Neumann data $u_z|_{\Gamma^+}$ also contains noise. We assume that the noise level is still $\delta$.
Denote by $f^{\delta}$ and $g^{\delta}$ the noisy boundary data and denote the corresponding noiseless data $f^*$ and $g^*$.
Here, by saying that $\delta$ is the noise level,  there exists an ``error" function $\mathcal E$ satisfying
\[
	\left\{
		\begin{array}{ll}
			\|\mathcal E\|_{H^p(\Omega)} \leq \delta,\\
			\mathcal E|_{\partial \Omega} = f^\delta - f^*,\\
			\mathcal E_z|_{\Gamma^+} = g^\delta - g^*.
		\end{array}
	\right.
\]
The following theorem guarantees that the minimizer of $J_{\lambda, \beta, \eta}$ subject to the boundary constraints defined by the noisy data  is an approximation of the true solution to \eqref{4.1}.
\begin{Theorem}
Assume that the set
\[
	H^{\delta} = \{
		u \in H^p(\Omega): u|_{\partial \Omega} = f^\delta \mbox{ and } u_z|_{\Gamma^+} = g^\delta
	\}
\] is nonempty 
	and let $u^{\delta}_{\min}$ be the minimizer of $J_{\lambda, \beta, \eta}$ in $H^{\delta}.$
	Assume further that Problem \eqref{4.1}, in which the Dirichlet and Neumann data $f$ and $g$  are replaced by $f^*$ and $g^*$, respectively, has the unique solution $u^{\epsilon_0} \in B(M - \delta)$.
Then, for all $\lambda$ and $\beta$ such that \eqref{convex} holds true, we have
\begin{equation}
	\|u_{\min}^\delta - u^{\epsilon_0}\|^2_{H^1(\Omega)} \leq C(\eta \|u^{\epsilon_0}\|^2_{H^p(\Omega)} + \delta^2).
	\label{4.1717}
\end{equation}
\label{thm 4.3}
\end{Theorem}

\proof
	Define $u = u^{\epsilon_0} + \mathcal{E}.$ 
	We have $u \in H^{\delta} \cap B(M).$
	Since $u^\delta_{\min}$ is the minimizer of $J_{\lambda, \beta, \eta}$ in $H^{\delta}$, by \cite[Lemma 2]{KlibanovNik:ra2017},
	\[\langle J_{\lambda, \beta, \eta}'(u^\delta_{\min}), u^\delta_{\min} - u^{\epsilon_0} - \mathcal{E} \rangle_{H^p(\Omega)}\leq  0.\]
	Applying \eqref{convex} for $u$ and $u_{\min}$ gives
	 \begin{align}
	J_{\lambda, \beta, \eta}(u^{\epsilon_0} + \mathcal{E})
	& \geq J_{\lambda, \beta, \eta}(u^{\epsilon_0} + \mathcal{E}) - J_{\lambda, \beta, \eta}(u^\delta_{\rm min}) - \langle J_{\lambda, \beta, \eta}'(u^{\delta}_{\rm min}),  u^{\epsilon_0} + \mathcal{E} - u^{\delta}_{\rm min} \rangle_{H^p(\Omega)}
	 \nonumber
	\\
	&\geq C \|u^{\epsilon_0} + \mathcal{E} - u^\delta_{\min}\|^2_{H^1(\Omega)} + \eta \|u^{\epsilon_0} + \mathcal{E} - u^\delta_{\min}\|^2_{H^p(\Omega)}
	\nonumber
	\\
	&\geq C \|u^{\epsilon_0} + \mathcal{E} - u^\delta_{\min}\|^2_{H^1(\Omega)}. \label{4.17}
\end{align}
Applying \eqref{4.2}, since $u^{\epsilon_0}$ solves \eqref{4.1}, we have
\begin{align}
	J_{\lambda, \beta, \eta}(u^{\epsilon_0} + \mathcal{E}) 
	&= J_{\lambda, \beta, \eta}(u^{\epsilon_0}) + J_{\lambda, \beta, \eta}'(u^{\epsilon_0})(\mathcal E) + C\delta^2
	= \eta \|u^{\epsilon_0}\|^2_{H^p(\Omega)} + 2\eta\langle u^{\epsilon_0}, \mathcal E \rangle_{H^p(\Omega)} + C\delta^2 \nonumber
	\\
	&\leq  \eta \|u^{\epsilon_0}\|^2_{H^p(\Omega)} + \eta^2\|u^{\epsilon_0}\|_{H^p(\Omega)}^2 +  C\delta^2
	\leq  2\eta \|u^{\epsilon_0}\|^2_{H^p(\Omega)} + C\delta^2.
	\label{4.18}
\end{align}
Using the inequality $(a - b)^2 \geq \frac{1}{2}a^2 - b^2$ and combining \eqref{4.17} and \eqref{4.18}, we have
\begin{align*}
	 \frac{1}{2}\|u^{\epsilon_0} - u^\delta_{\min}\|^2_{H^1(\Omega)} -  \|\mathcal{E}\|^2_{H^1(\Omega)}
	  \leq  \|u^{\epsilon_0} + \mathcal{E} - u^\delta_{\min}\|^2_{H^1(\Omega)} 
	  \leq C(\eta \|u^{\epsilon_0}\|^2_{H^p(\Omega)} + \delta^2)
\end{align*}
The estimate \eqref{4.1717} follows.
\QEDB


\medskip

Combining Theorem \ref{thm 4.2} and Theorem \ref{thm 4.3}, we have for each $k \geq 1$,
\[
	\|u^{(k)} - u^{\epsilon_0}\|_{H^1(\Omega)} \leq C(\sqrt{\eta} \|u^{\epsilon_0}\|_{H^p(\Omega)} + \delta) + 
	\theta^{k/2} \|u^{(0)} - u^\delta_{\min}\|_{H^p(\Omega)}
\]
for some constant $C > 0$ and $\theta \in (0, 1)$ where $u^{(k)}$ is the minimizing sequence defined in Theorem \ref{thm 4.2}. This inequality shows the stability of our method with respect to noise. If $\theta^{k}$ and $\eta$ are $O(\delta^2)$, then the convergence rate is Lipschitz.

\begin{remark}
In the case when the function $F$ is such that the comparison principle for $-\epsilon_0 \Delta u + F(\x, u, \nabla u) = 0$ holds true; for e.g., $F$ is strictly increasing with respect to its second variable (Assumption \ref{assumption 2}), one can prove the stability of $u$ with respect to noise without imposing Assumption \ref{assumption 1} by the fact that $u - \delta$ and $u+\delta$ are a subsolution and supersolution of $-\epsilon_0 \Delta u + F(\x, u, \nabla u) = 0$, respectively. The main reason for us to successfully establish the stability without assuming the comparison principle is due to the presence of the Neumann data in Assumption \ref{assumption 1}.
\label{rm 4.2}
\end{remark}
In the next section, we present the numerical implementation and some interesting numerical results.


\section{Numerical study} \label{sec:num}

We implement the convexification method based on the finite difference method. For the simplicity, in this section, we only consider the case $d = 2.$
On $\overline{\Omega} = [-R, R]^2$, we arrange a uniform grid of points 
\begin{equation}
	\mathcal{G} = \big\{\x_{ij} = (x_i, z_j): x_i = -R + (i - 1) \delta_x, z_j = -R + (j - 1)\delta_z), 1 \leq i, j \leq N\big\}
	\label{5.1}
\end{equation}
where $\delta_x = \delta_z = h = 2R/(N - 1)$.
In our computation, $R = 1$ and $N = 50.$ 
The finite difference version of the objective function $J_{\lambda, \beta, \eta}(u)$ is given by 
\begin{multline}
	J_{\lambda, \beta, \eta}(u) = h^2 \sum_{i, j = 2}^{N-1} e^{2\lambda \big|\frac{z_j + r)}{b}\big|^\beta}
	\Big|
		-\epsilon_0 \Delta^h u(x_i, z_j)
		 + F((x_i, z_j), u(x_i, z_j), \partial_x^h u(x_i, z_j), \partial_z^h u(x_i, z_j))
	\Big|^2
	\\
	+ \eta h^2   \Big(\sum_{i, j = 1}^{N}|u(x_i, z_j)|^2 +\sum_{i, j = 2}^{N-1} \big( |\partial_x^h u(x_i, z_j)|^2 + |\partial_z^h u(x_i, z_j)|^2 + |\Delta^h u(x_i, z_j)|^2\big)\Big).
\end{multline}
where
\begin{align*}
	\Delta^h u(x_i, z_j) &= \frac{u(x_i, z_{j-1}) + u(x_i, z_{j + 1}) + u(x_{i-1}, z_j) + u(x_{i+1}, z_{j}) -4u(x_i, z_j)  }{h^2}\\
	\partial_x^h u(x_i, z_j) &=\frac{u(x_{i+1}, z_{j}) - u(x_{i - 1}, z_{j})}{2h}\\
	\partial_z^h u(x_i, z_j) &=  \frac{u(x_i, z_{j+1}) - u(x_i, z_{j-1})}{2h}.
\end{align*}

\begin{remark}
1. In our computation $R = 1$, $N = 50$, $\epsilon_0 = 10^{-3}$, $\lambda = 2,$ $\beta = 8$, $b = 10$, $r = 1.2$ and $\eta = 10^{-4}$. Although in the theoretical part, the parameters $\lambda$ and $\beta$ should be large,  these values are already good for the numerical part. For the simplicity, we use this set of parameters for all tests below.

2.	In theory $p > 3$ when $d = 2$. However, in numerical study, we can reduce the norm in the regularization term to $p = 2$ to simplify the implementation and to improve the speed of computation. We do not experience any difficulty with this small change.
\end{remark}

In our implementation, instead of writing the computational code for the gradient descent method, we use the optimization toolbox of Matlab, in which the gradient descent method is coded. 
More precisely, we use the command ``fmincon" to minimize the functional $J_{\lambda, \beta, \eta}$ subject to the boundary constraint in \eqref{4.1}. 
The command ``fmincon" requires an initial solution, that is the function $u^{(0)}$ in Theorem \ref{thm 4.2}.
The function $u^{(0)}$ is naturally assigned to be the zero function.
This function $u^{(0)}$ does not satisfy the Dirichlet boundary condition on $\partial \Omega$ and Neumann condition on $\Gamma^+$. However, the command ``fmincon" corrects this error automatically.


We present here five (5) numerical tests in which the given boundary data are noisy with noise level $\delta = 0\%$, $5\%$, and $10\%$, respectively.
In each test, $u_{\rm true}$ and $u_{\rm comp}$ denote the true and computed viscosity solutions, respectively.
Given functions $f,g$, the noisy versions of $f, g$ are given by
\[
	f^\delta = f(1 + \delta \cdot \mbox{rand}) \quad \text{ and } \quad
	g^\delta = g(1 + \delta \cdot \mbox{rand}),
\]
where rand is the function that generates uniformly distributed random numbers in $[-1, 1].$ 
The relative computed error is defined as
\[
	{\rm err}(\delta) = \frac{\|u_{\rm comp} - u_{\rm true}\|_{L^{\infty}(\Omega)}}{\|u_{\rm true}\|_{L^{\infty}(\Omega)}} \quad \text{ for } \delta > 0.
\]

{\it Test 1.} In this test, we find the viscosity solution to \eqref{HJ}--\eqref{dir} where 
\begin{equation}
	F(\x, s, \p) =  \frac{1}{150} s + |\p| + \frac{1}{150} (x^2 + z^2) - 2\sqrt{x^2 + z^2} \quad \mbox{for all } \x \in \Omega, s \in \R, \p \in \R^2
	\label{F1}
\end{equation} and the boundary data are given by
\begin{equation}
	u(\x) = f(\x) =  -(x^2 + z^2)\quad \mbox{for all } \x = (x, z) \in \partial\Omega
	\label{f1}
\end{equation}
and
\begin{equation}
	u_z(\x) = g(\x) =  -2z \quad \mbox{for all } \x = (x, z) \in \Gamma^+.
	\label{g1}
\end{equation}
The true solution is $u_{\rm true}(x, z) = -(x^2 + z^2)$, which is smooth in $\Omega$.
We are here in a standard setting, and the convergence of $u^\epsilon$, the solution to \eqref{2.1}, to the solution $u_{\rm true}$ to \eqref{HJ}--\eqref{dir} is guaranteed in \cite{CrandallLions83, CrandallEvansLions84, Lions, Barles, BCD,Tran19}.

\begin{figure}[h!]
		\subfloat[The true solution $u_{\rm true} = -(x^2 + z^2)$.]{\includegraphics[width=.3\textwidth]{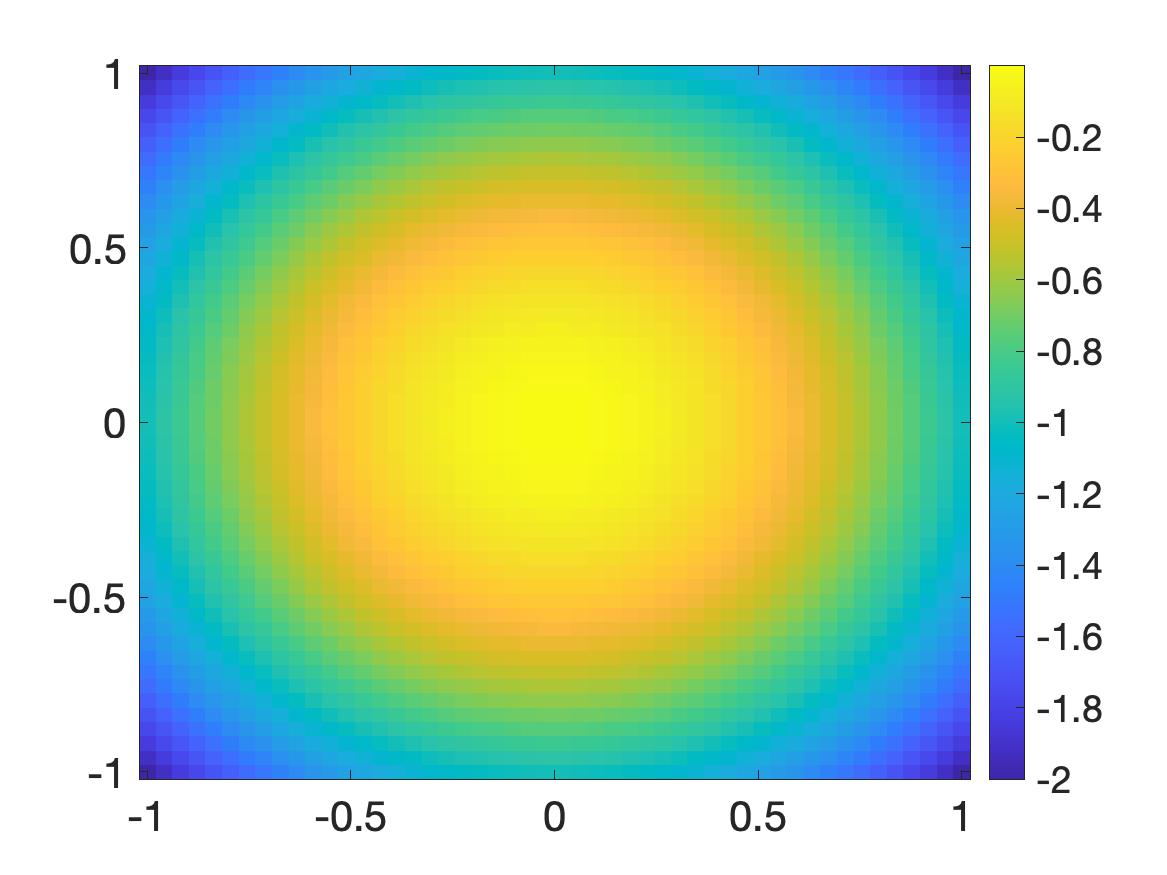}}
		
		\subfloat[The solution $u_{\rm comp}$, computed from noiseless boundary data.]{\includegraphics[width=.3\textwidth]{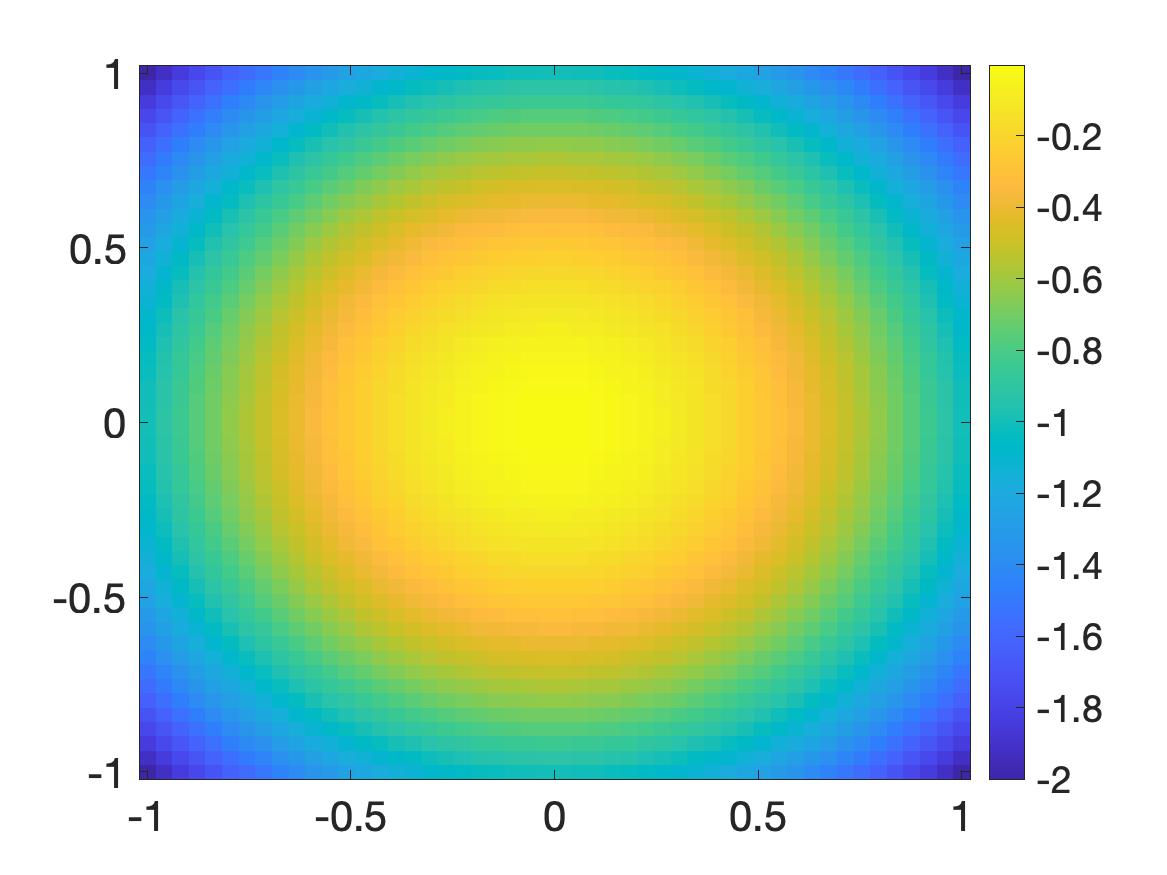}}
		\quad
		\subfloat[The solution $u_{\rm comp}$, computed from 5\% noisy boundary data.]{\includegraphics[width=.3\textwidth]{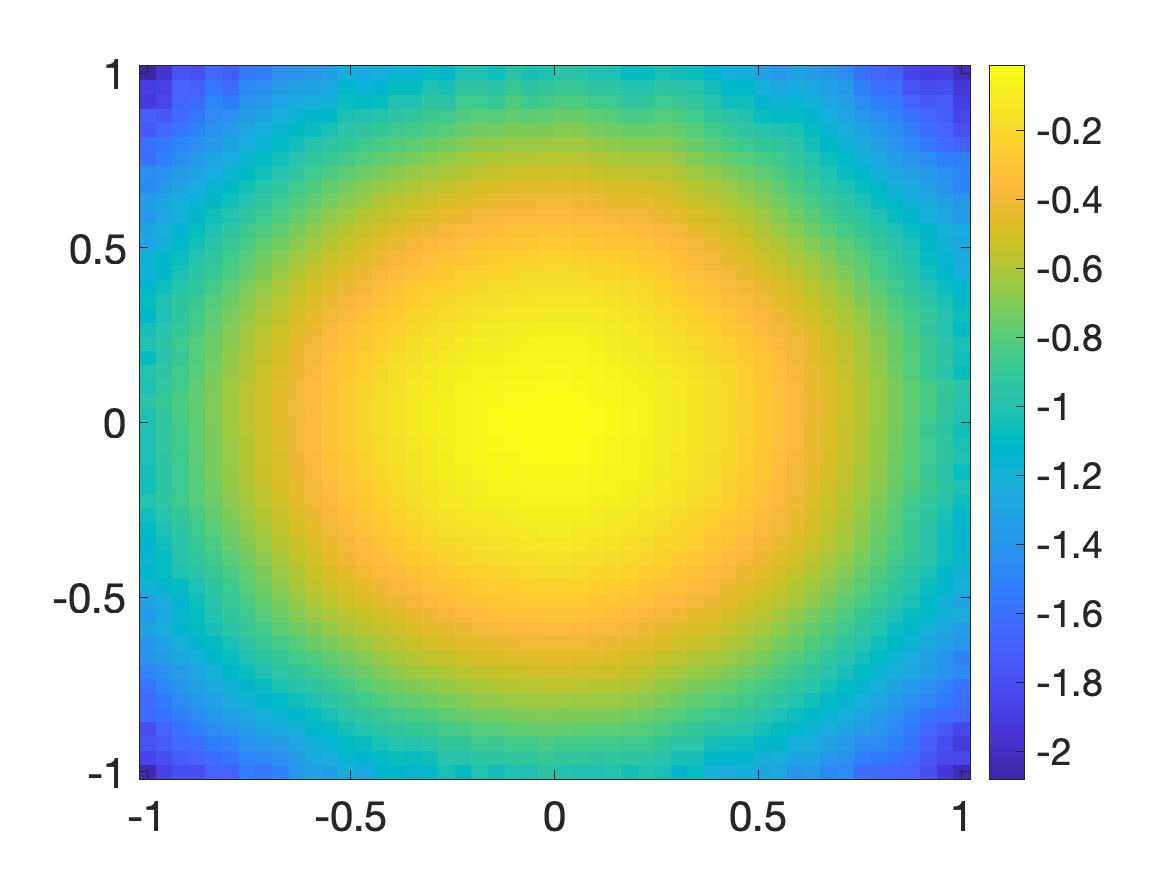}}
		\quad
		\subfloat[The solution $u_{\rm comp}$, computed from 10\% noisy boundary data.]{\includegraphics[width=.3\textwidth]{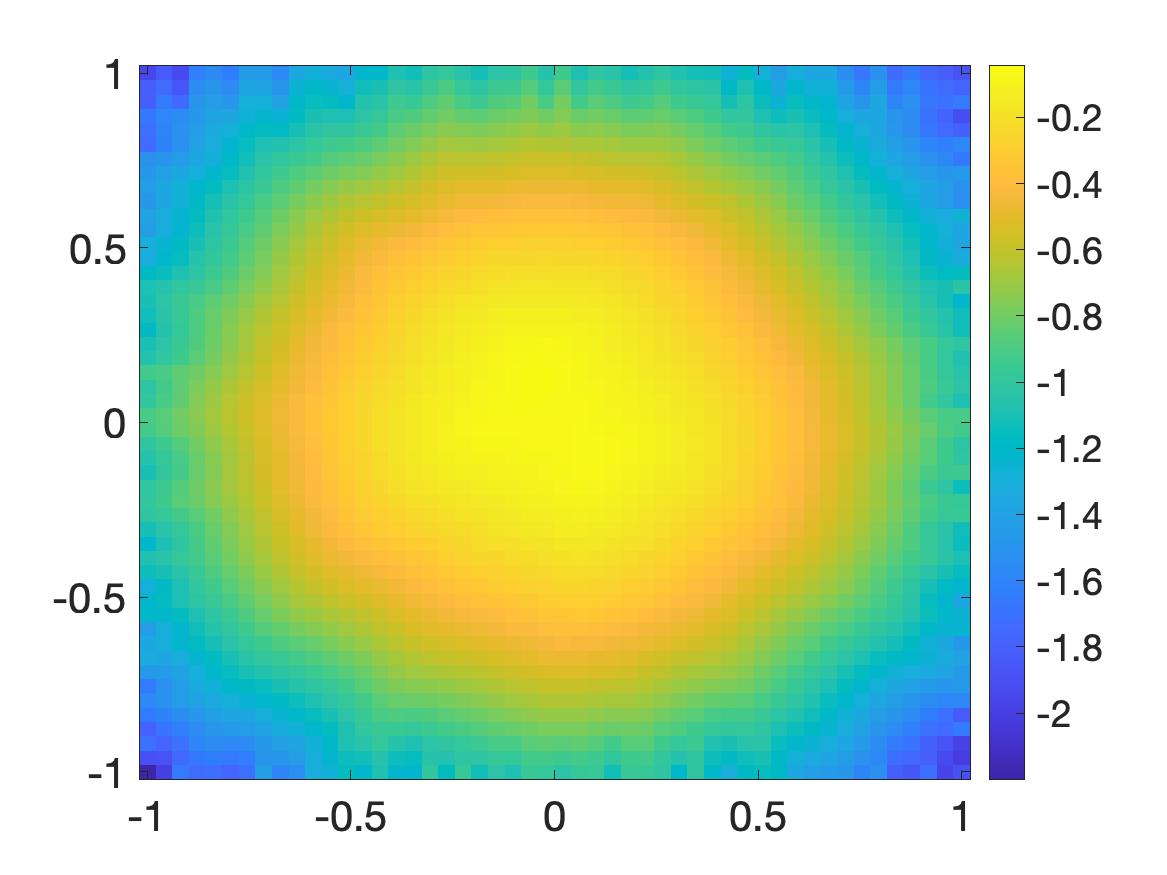}}
		
		\subfloat[{The relative error $\frac{|u_{\rm comp} - u_{\rm true}|}{\|u_{\rm true}\|_{L^{\infty}}}$, $\delta = 0\%$.}]{\includegraphics[width=.3\textwidth]{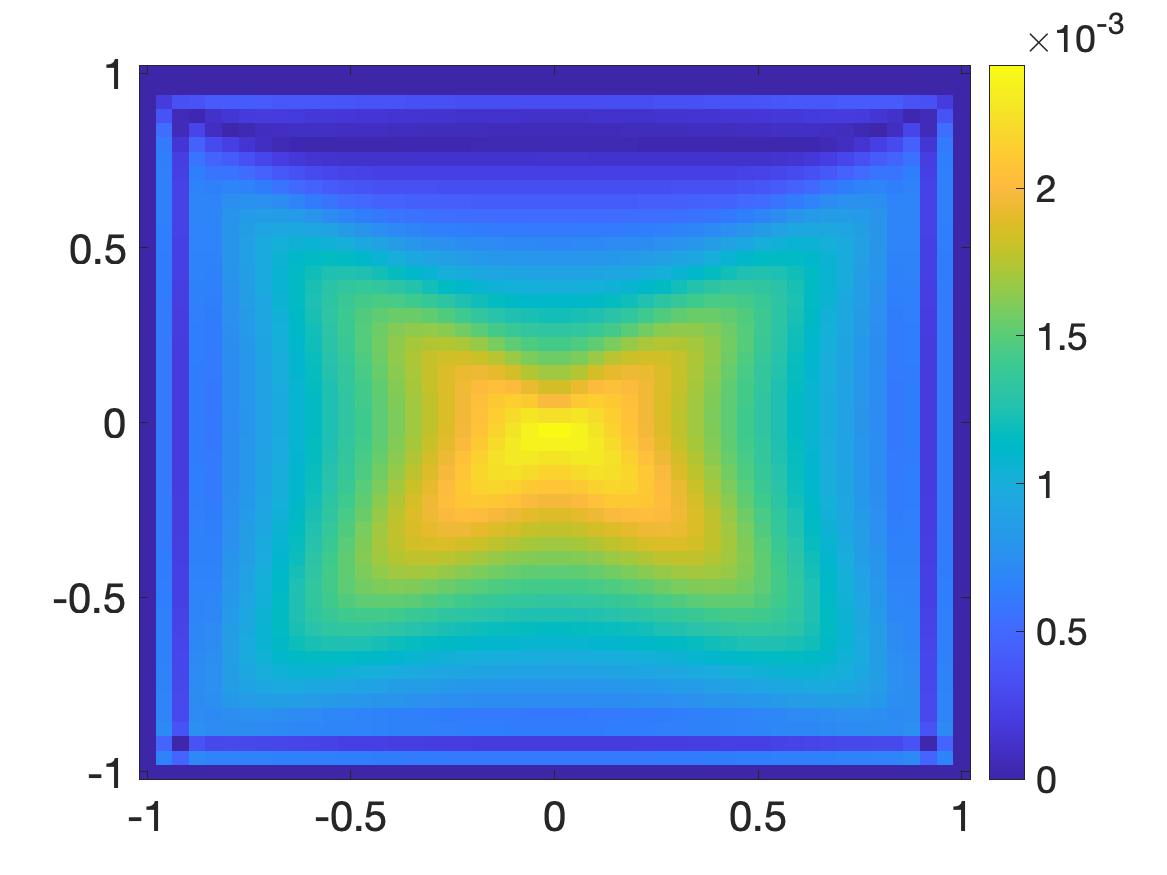}}
		\quad
		\subfloat[[{The relative error $\frac{|u_{\rm comp} - u_{\rm true}|}{\|u_{\rm true}\|_{L^{\infty}}}$, $\delta = 5\%$.}]{\includegraphics[width=.3\textwidth]{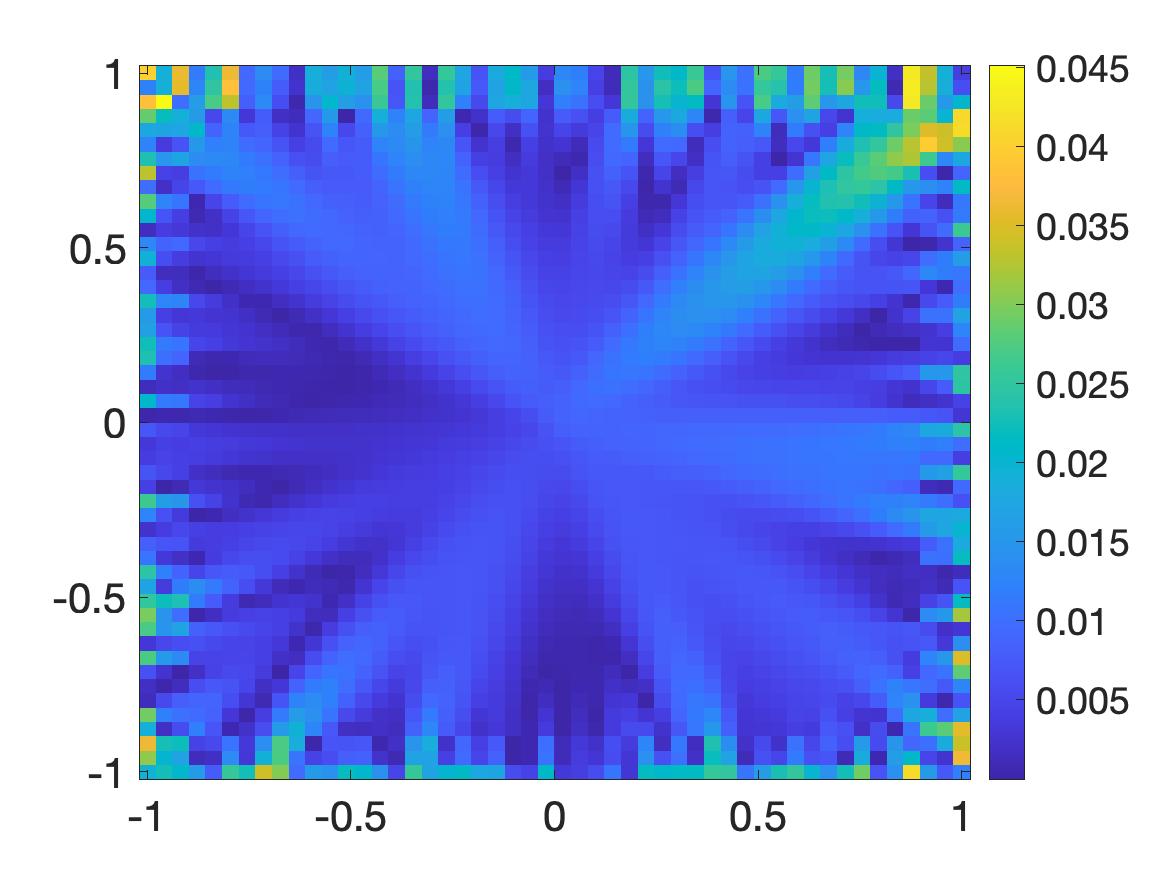}}
		\quad
		\subfloat[[{The relative error $\frac{|u_{\rm comp} - u_{\rm true}|}{\|u_{\rm true}\|_{L^{\infty}}}$, $\delta = 10\%$.}]{\includegraphics[width=.3\textwidth]{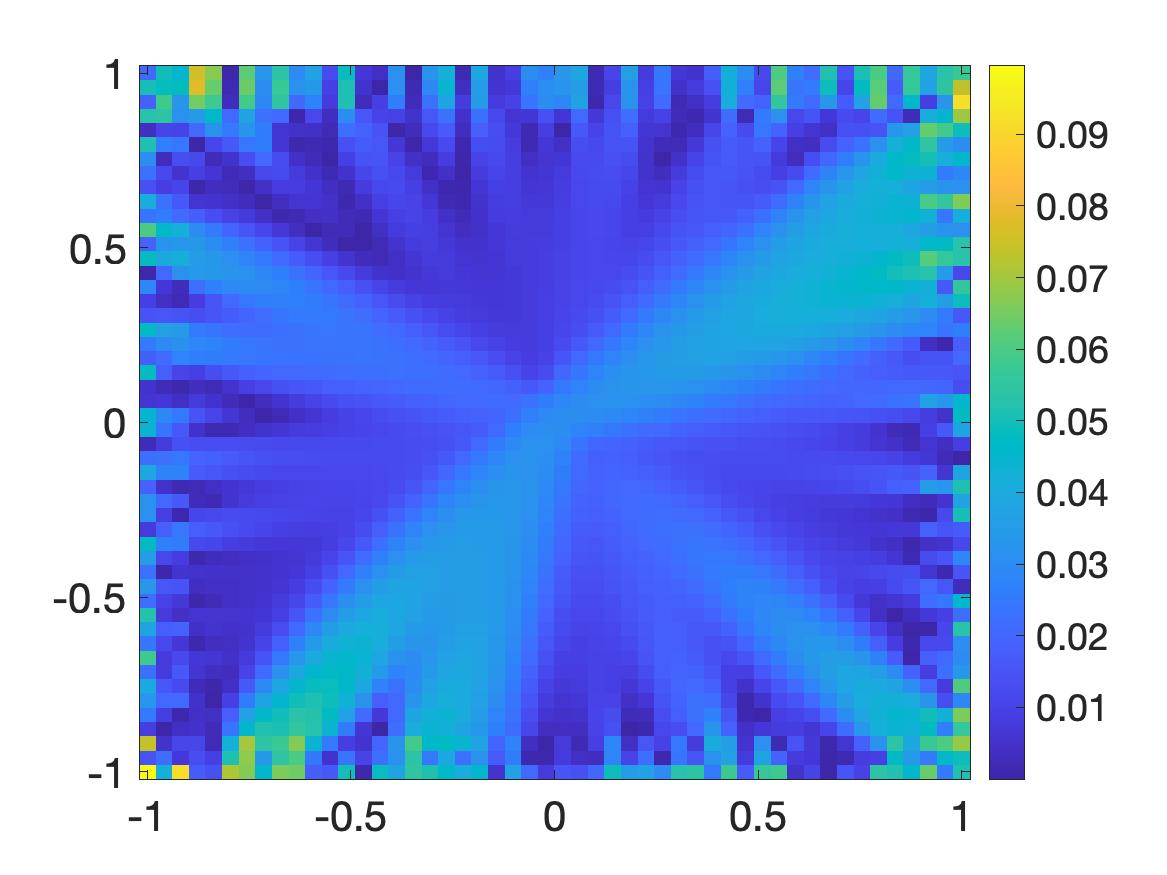}}
		
		\caption{\label{fig 1} Test 1. The true and computed viscosity solutions with $\delta$ is $0\%$,  5\%, and 10\% noisy boundary Dirichlet data on $\partial \Omega$ and Neumann data on $\Gamma^+$.  }
\end{figure}

We show in Figure \ref{fig 1} the graphs of this true solution and that of the computed ones from noiseless and noisy boundary data. 
It is evident that we successfully obtained computed viscosity solutions to \eqref{HJ}--\eqref{dir}.  By adding high level noise into the boundary data, we have numerically shown that the convexification is stable.
{It is evident from Table \ref{tab1} that the relative computed error is about the noise level, 
which clearly illustrates the Lipschitz stability in Theorem \ref{thm 4.3}.}
In this test, the true solution is smooth. 
The function $F(\x, s, \p)$ is strictly increasing with respect to $s$.

\begin{table}[h!]
\caption{\label{tab1} {Test 1. The performance of the convexification method. The computational time is the time for a Precisions Workstations T7810 with 24 cores to compute the solution $u_{\rm comp}$. In this table, 
the relative $L^\infty(\Omega)$ error is $\|u_{\rm comp} - u_{\rm true}\|_{L^\infty(\Omega)}/\|u_{\rm true}\|_{L^\infty(\Omega)}$}.} 
{
\centering
\begin{tabular}{c c c  c }
\hline\hline
Noise level &computational time&  number of iterations &  relative $ L^\infty(\Omega)$ error\\ 
\hline
0\%&23.47 minutes&279 &$0.24\%$  \\
5\%& 24.40 minutes &292&4.51\%\\
10\%& 27.35 minutes &329&9.95\%\\
\hline
\end{tabular}
}
\end{table}

\medskip

\begin{figure}[h!]
		\subfloat[The true solution $u_{\rm true} = -(|x| + |z|)$.]{\includegraphics[width=.3\textwidth]{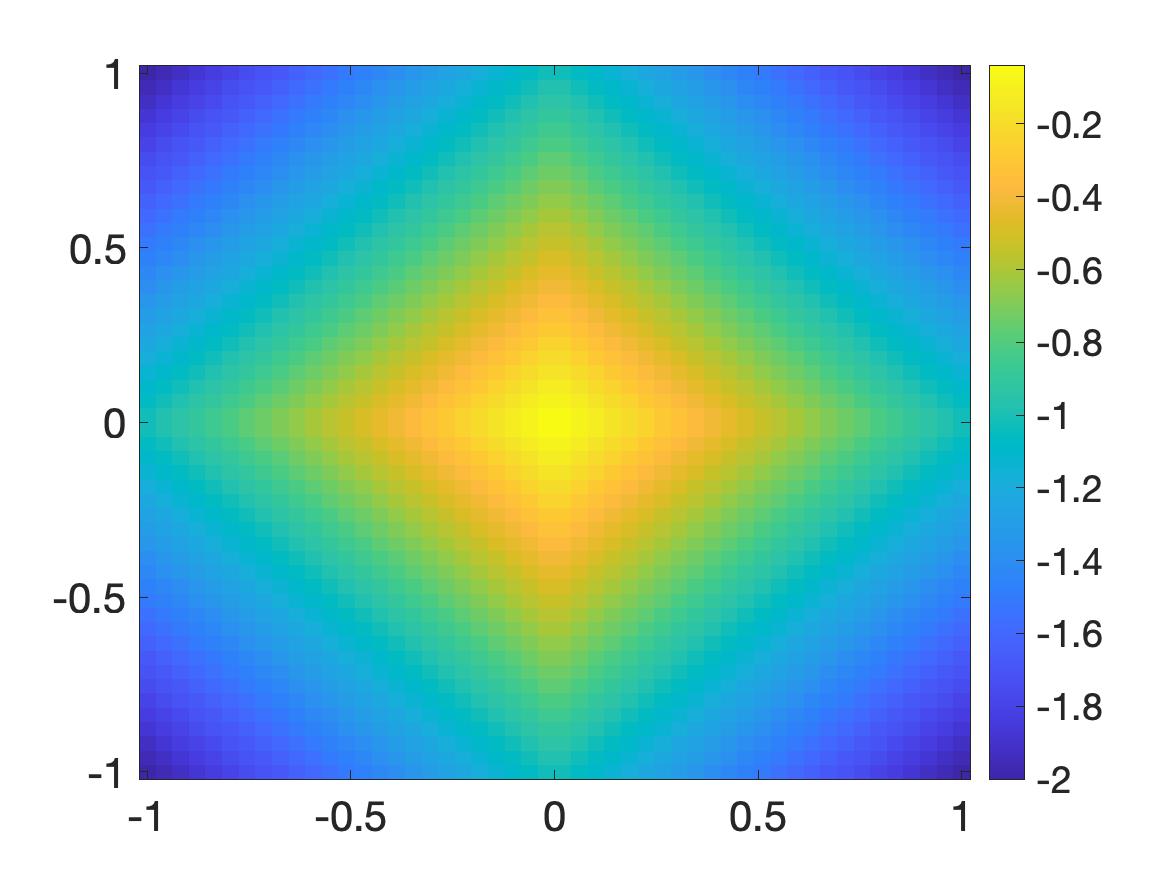}}
	
		\subfloat[The solution $u_{\rm comp}$, computed from noiseless boundary data.]{\includegraphics[width=.3\textwidth]{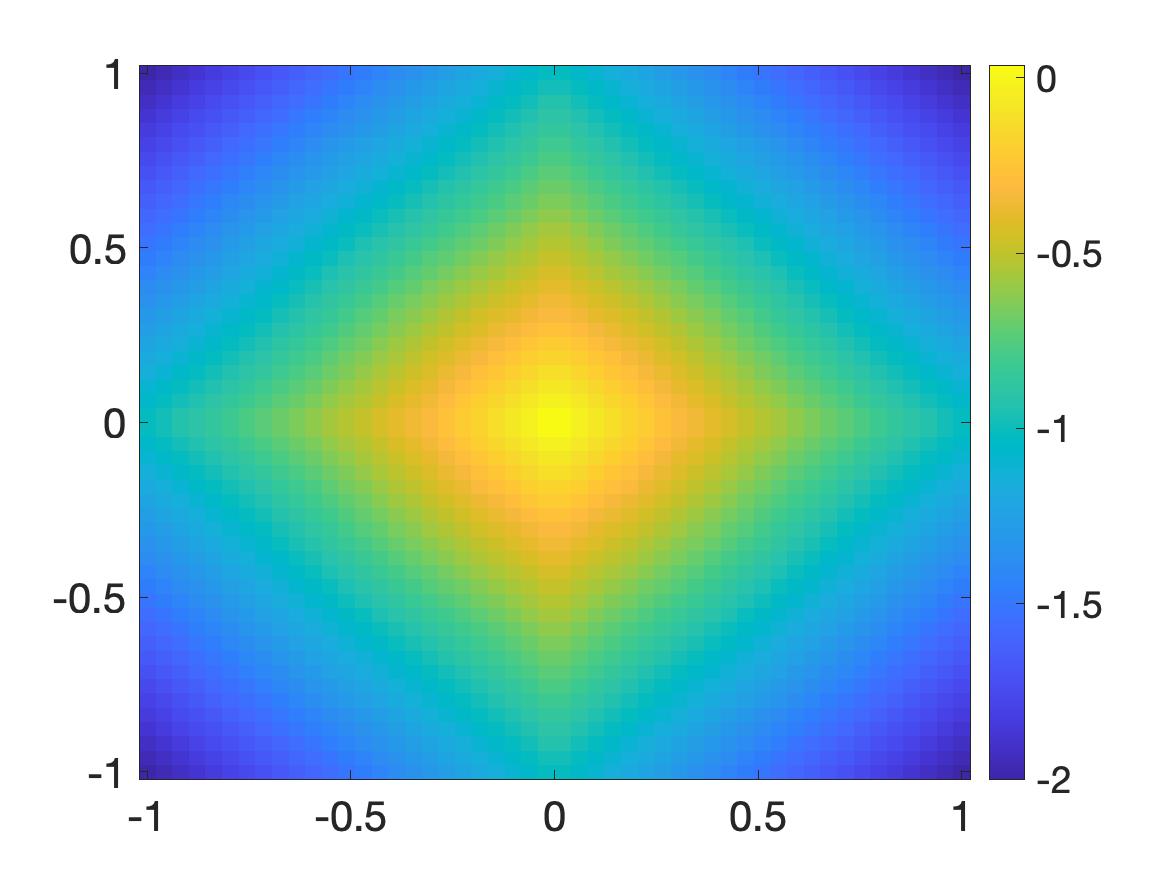}}
		\subfloat[The solution $u_{\rm comp}$, computed from 5\% noisy boundary data.]{\includegraphics[width=.3\textwidth]{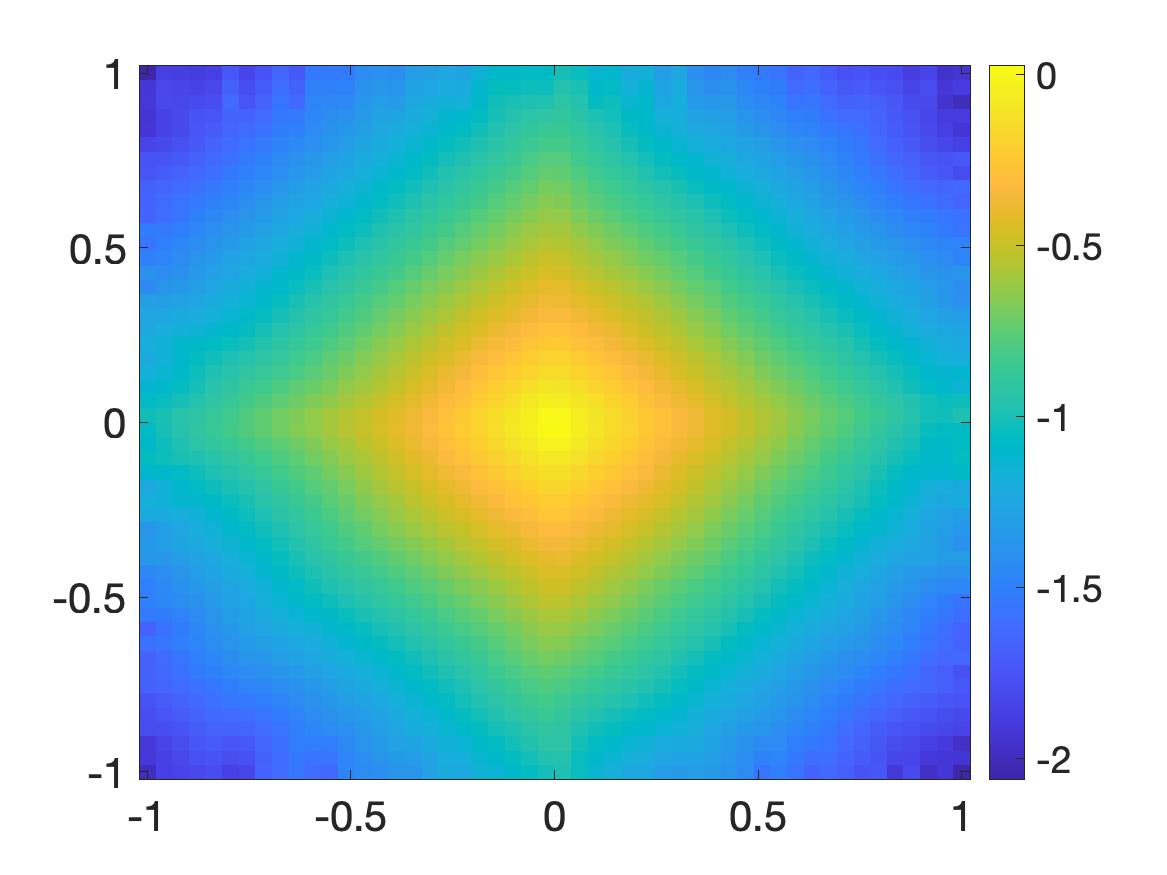}}
		\quad		
		\subfloat[The solution $u_{\rm comp}$, computed from 10\% noisy boundary data.]{\includegraphics[width=.3\textwidth]{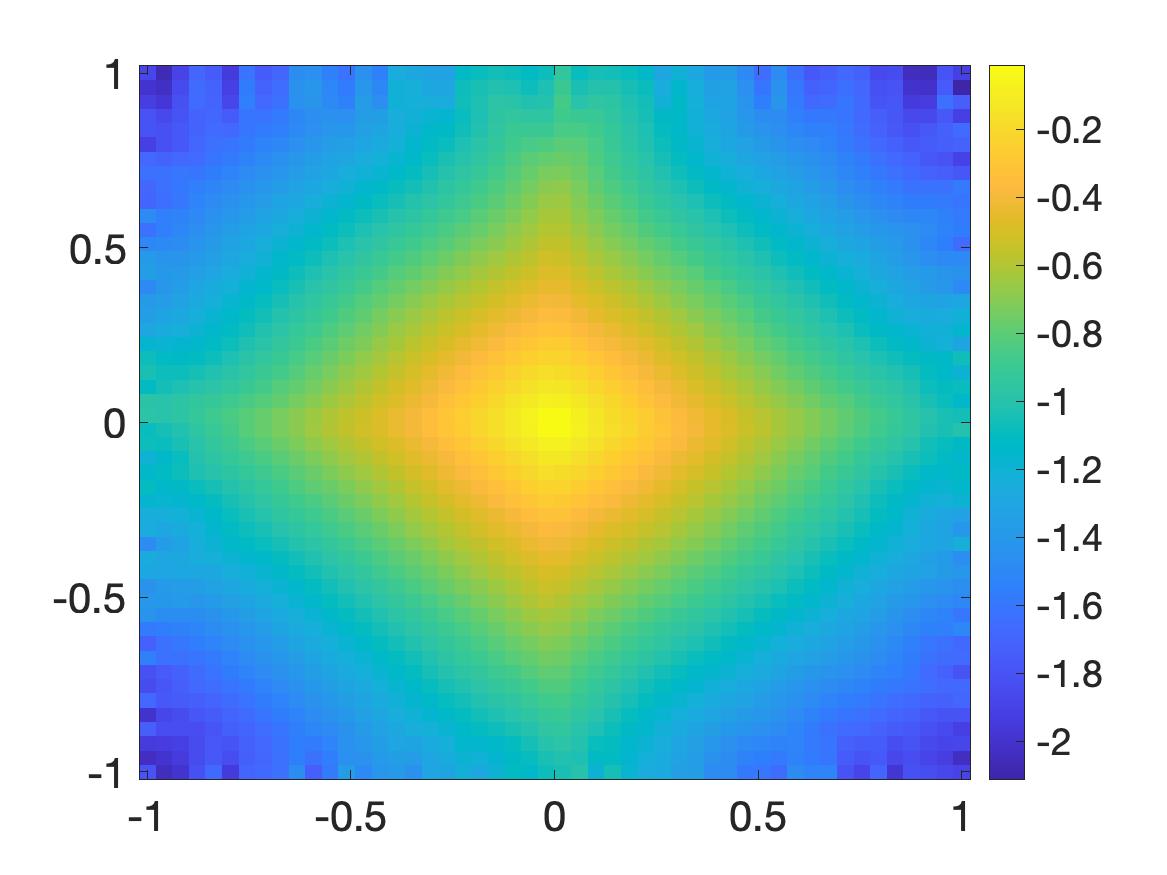}}

	\subfloat[{The relative error $\frac{|u_{\rm comp} - u_{\rm true}|}{\|u_{\rm true}\|_{L^{\infty}}}$, $\delta = 0\%$.}]{\includegraphics[width=.3\textwidth]{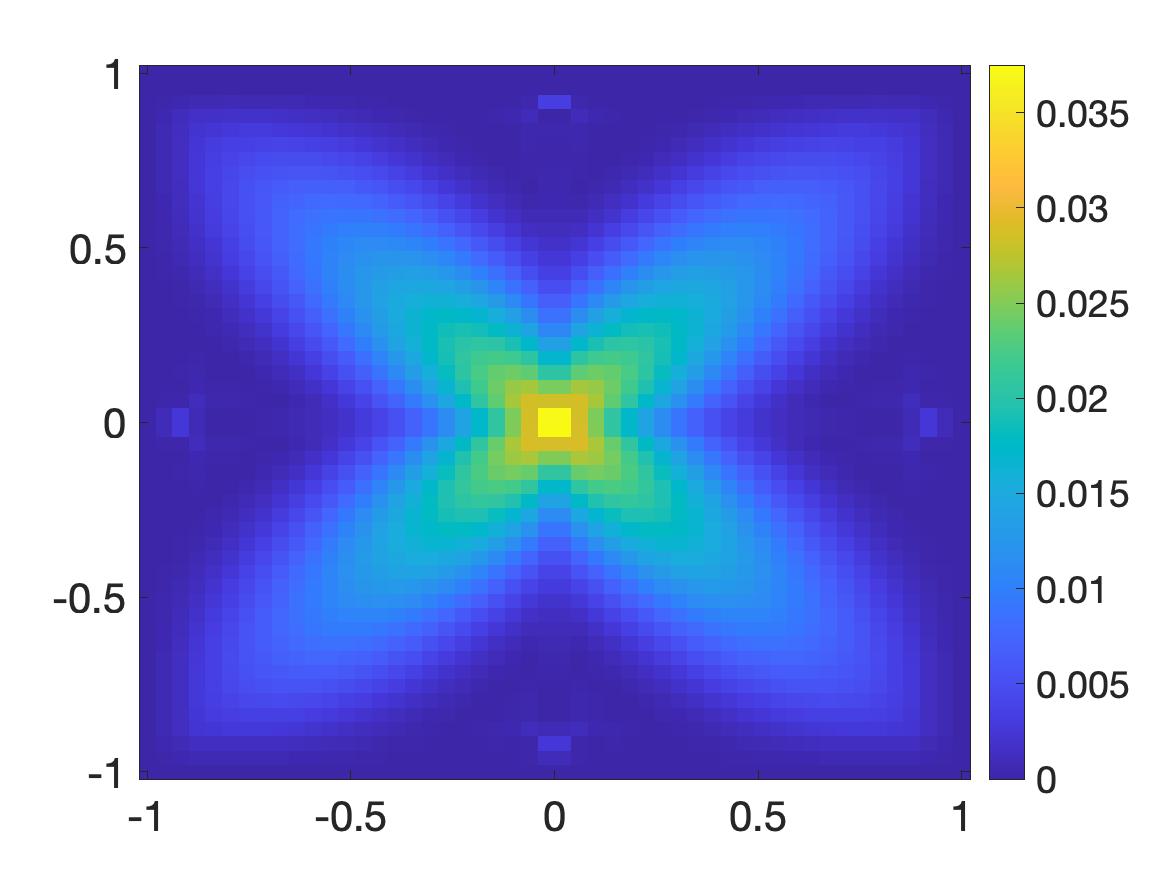}}
		\quad
		\subfloat[[{The relative error $\frac{|u_{\rm comp} - u_{\rm true}|}{\|u_{\rm true}\|_{L^{\infty}}}$, $\delta = 5\%$.}]{\includegraphics[width=.3\textwidth]{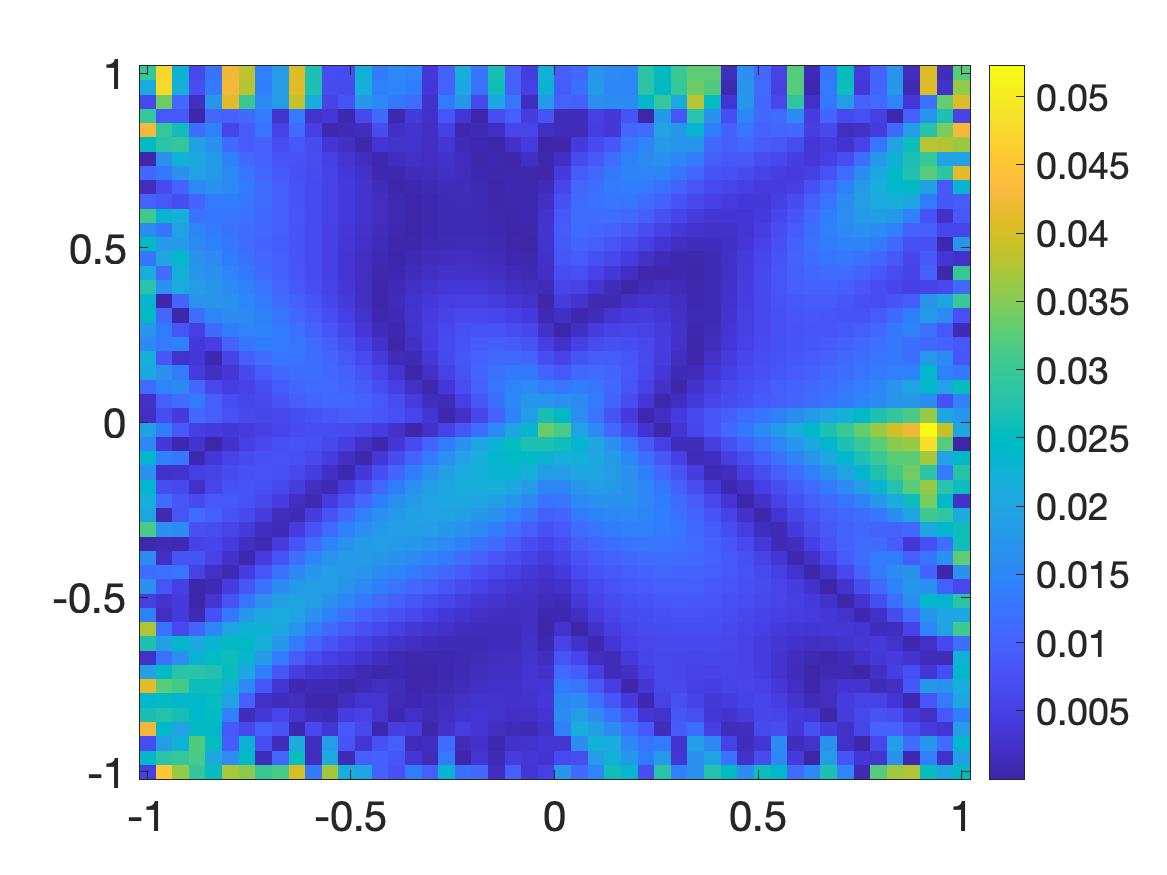}}
		\quad
		\subfloat[[{The relative error $\frac{|u_{\rm comp} - u_{\rm true}|}{\|u_{\rm true}\|_{L^{\infty}}}$, $\delta = 10\%$.}]{\includegraphics[width=.3\textwidth]{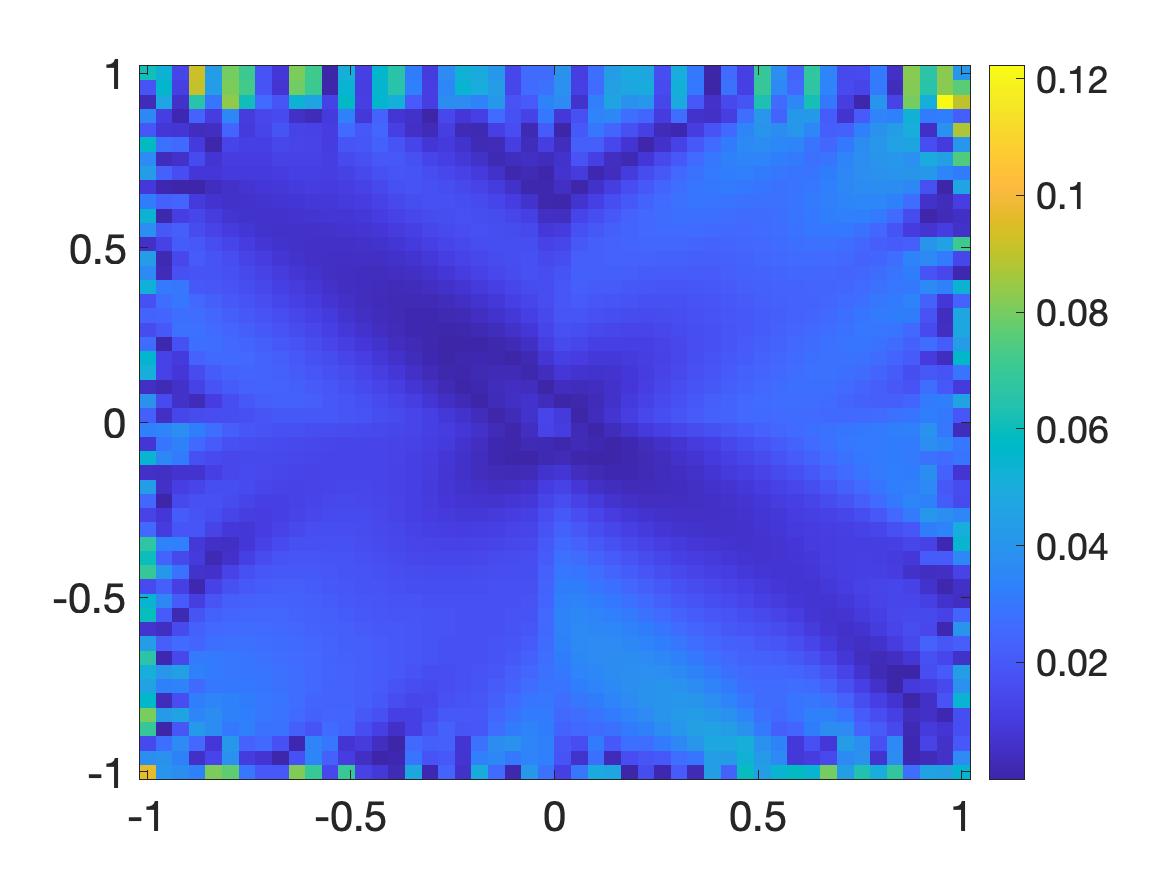}}

		\caption{\label{fig 2}  Test 2. The true and computed viscosity solutions with $\delta$ is $0\%$,  5\%, and 10\% noisy boundary Dirichlet data on $\partial \Omega$ and Neumann data on $\Gamma^+$.  }
\end{figure}

{\it Test 2.} We now solve the {\it eikonal} equation of which the function $F$ in \eqref{HJ} is not in the class $C^1$. 
In this test, the function $F$ is given by
\begin{equation}
	F(\x, s, \p) = |\p| - \sqrt{2} 
	\quad \mbox{for all } \x \in \Omega, s \in \R, \p \in \R^2
	\label{F2}
\end{equation}
and the boundary data are
\begin{equation}
	u(\x) = f(\x) = -(|x| + |z|) 
	\quad \mbox{for all } \x = (x, z) \in \partial \Omega
	\label{f2}
\end{equation}
and
\begin{equation}
	u_z(\x) = g(\x) = \left\{
		\begin{array}{rl}
			1 &z < 0,\\
			-1 &z > 0,
		\end{array}
	\right.
	\quad \mbox{for all } \x = (x, z) \in \Gamma^+.
	\label{g2}
\end{equation}
We claim that the true solution to \eqref{HJ}--\eqref{dir} is $u_{\rm true}(\x) = -(|x| + |z|)$ for all $\x = (x, z) \in \Omega.$
Intuitively, this claim holds as the graph of $u_{\rm true}$ only has corners from above and $F$ is convex in $\p$.
Let us provide a rigorous verification here.
If $x \neq 0$ and $z \neq 0$, then $u_{\rm true}$ is differentiable at $\x=(x,z)$, and
\[
\nabla u_{\rm true}(\x) =\left(-\frac{x}{|x|}, - \frac{z}{|z|} \right ) \quad \Rightarrow \quad |\nabla u_{\rm true}(\x) |=\sqrt{2}.
\]
If $x=0$ or $z=0$, then $u_{\rm true}$ is not differentiable at $\x=(x,z)$.
We can only find smooth test functions that touch $u_{\rm true}$ from above at $\x$, and we cannot find smooth test functions that touch $u_{\rm true}$ from below at $\x$.
Let $\phi$ be a smooth test function that touches $u_{\rm true}$ from above at $\x$.
Without loss of generality, we only need to consider the case $x=0$.
If $z \neq 0$, then we have that
\[
\phi_x(\x) \in [-1,1], \ \phi_z(\x)=-\frac{z}{|z|}.
\]
If $z = 0$, then we have that
\[
\phi_x(\x) \in [-1,1], \ \phi_z(\x) \in [-1,1].
\]
In both cases,
\[
|\nabla \phi(\x) | \leq \sqrt{2}.
\]
Thus, the subsolution test holds for $u_{\rm true}$ at $\x$.
We conclude that $u_{\rm true}$ is a viscosity solution  to \eqref{HJ}--\eqref{dir}.
This true solution and its computed versions $u_{\rm comp}$ from noisy boundary data are displayed in Figure \ref{fig 2}.
The convergence of $u^\epsilon$, the solution to \eqref{2.1}, to the solution $u_{\rm true}$ to \eqref{HJ}--\eqref{dir} is guaranteed in \cite{FlemingSouganidis, Tran11} with convergence rate $O(\epsilon^{1/2})$.

This test is more challenging than Test 1. 
In this test, the functions $F$ and $f$ are both not smooth. 
So, the smoothness condition for the theoretical part does not satisfy. 
However, the convexification method still provides reliable solutions even when the given boundary is  noisy with the noise level $\delta = 10\%.$
This shows the robustness of the convexification method.
{The relative errors are acceptable, see Table \ref{tab2}.}

\begin{table}[h!]
\caption{\label{tab2} {Test 2. The performance of the convexification method. The computational time is the time for a Precisions Workstations T7810 with 24 cores to compute the solution $u_{\rm comp}$. In this table, 
the relative $L^\infty(\Omega)$ error is $\|u_{\rm comp} - u_{\rm true}\|_{L^\infty(\Omega)}/\|u_{\rm true}\|_{L^\infty(\Omega)}$.}} 
\centering
{
\begin{tabular}{c c c  c }
\hline\hline
Noise level &computational time&  number of iterations &  relative $ L^\infty(\Omega)$ error\\ 
\hline
0\%&  48 minutes&581&3.75\%\\
5\%& 89 minutes&1071& 5.23\%\\
10\%& 64 minutes & 766 &12.22 \%\\
\hline
\end{tabular}
}
\end{table}

\medskip

{\it Test 3.}
We next consider a more complicated Hamilton-Jacobi equation. Unlike the previous two tests, the function $F$ in this test is not convex and; more interestingly, not coercive and not continuous. It is given by
\begin{multline}
	F(\x, s, \p) =20 s +  |p_1| - |p_2| 
	\\
	- \Big(20\big(-|x - 0.5| + e^{\sin(\pi(x^2 + z^2))}\big) + G(x, z) - \Big|2\pi z \cos(\pi(x^2 + z^2)) e^{\sin(\pi(x^2 + z^2))}\Big| \Big) 
	\label{5.9999}
\end{multline} 
for all $\x = (x, z) \in \Omega, s \in \R, \p= (p_1, p_2) \in \R^2$ where
\[
	G(x, z) = 
	\left\{
		\begin{array}{ll}
			 \big|1 + 2\pi x \cos(\pi(x^2 + z^2)) e^{\sin(\pi(x^2 + z^2))}\big|  &x < 0.5,\\
			 \big|-1 + 2\pi x \cos(\pi(x^2 + z^2)) e^{\sin(\pi(x^2 + z^2))}\big| &x > 0.5.
		\end{array}
	\right.
\]
Note that $G$ and $F$ are discontinuous at $x=0.5$ in general.
The boundary data are given by
\begin{equation}
	u(\x) = f(\x) =  - |x - 0.5| + e^{\sin(\pi(x^2 + z^2))} \quad \mbox{for all } \x = (x, z) \in \partial \Omega
\end{equation}
and
\begin{equation}
	u_z(\x) = g(\x) = 2\pi z \cos(\pi(x^2 + z^2)) e^{\sin(\pi(x^2 + z^2))} \quad \mbox{for all } \x = (x, z) \in \Gamma^+.
\end{equation}
The function $u_{\rm true}(x, z) = - |x - 0.5| + e^{\sin(\pi(x^2 + z^2))} $ is the true viscosity solution for this test.
Indeed, if $x \neq 0.5$, then $u_{\rm true}$ is differentiable at $\x=(x,z)$, and
\[
\nabla u_{\rm true}(\x) =\left(-\frac{x-0.5}{|x-0.5|}+ 2\pi x \cos(\pi(x^2 + z^2)) e^{\sin(\pi(x^2 + z^2))}, 2\pi z \cos(\pi(x^2 + z^2)) e^{\sin(\pi(x^2 + z^2))} \right ), 
\]
which gives that $F(\x,u_{\rm true}(\x),\nabla u_{\rm true}(\x) )=0$.
If $x=0.5$, then $u_{\rm true}$ is not differentiable at $\x=(x,z)$.
We can only find smooth test functions that touch $u_{\rm true}$ from above at $\x$, and we cannot find smooth test functions that touch $u_{\rm true}$ from below at $\x$.
Let $\phi$ be a smooth test function that touches $u_{\rm true}$ from above at $\x$.
Then,
\[
\phi_x(\x) \in 2\pi x \cos(\pi(x^2 + z^2)) e^{\sin(\pi(x^2 + z^2))}+[-1,1], \ \phi_z(\x)=2\pi z \cos(\pi(x^2 + z^2)) e^{\sin(\pi(x^2 + z^2))},
\]
which yields that $F_*(\x,u_{\rm true}(\x),\nabla \phi(\x) ) \leq 0$.
Here, $F_*$ is the lower semicontinuous envelope of $F$.
Therefore, the subsolution test holds for $u_{\rm true}$ at $\x$.
We imply that $u_{\rm true}$ is a viscosity solution  to \eqref{HJ}--\eqref{dir}.

\begin{figure}[h!]
		\subfloat[The true solution $u_{\rm true}(x, z) = - |x - 0.5| + e^{\sin(\pi(x^2 + z^2))} $.]{\includegraphics[width=.3\textwidth]{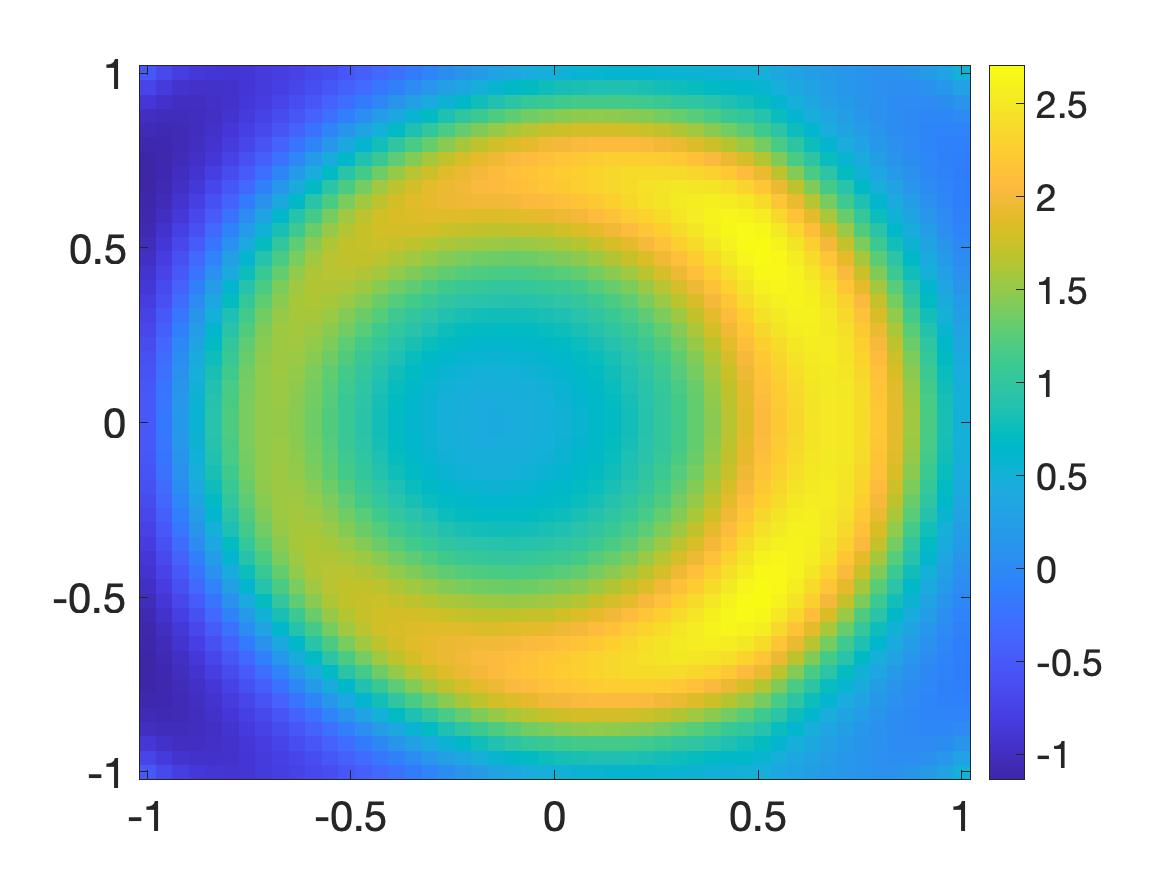}}
		
		\subfloat[The solution $u_{\rm comp}$, computed from noiseless boundary data.]{\includegraphics[width=.3\textwidth]{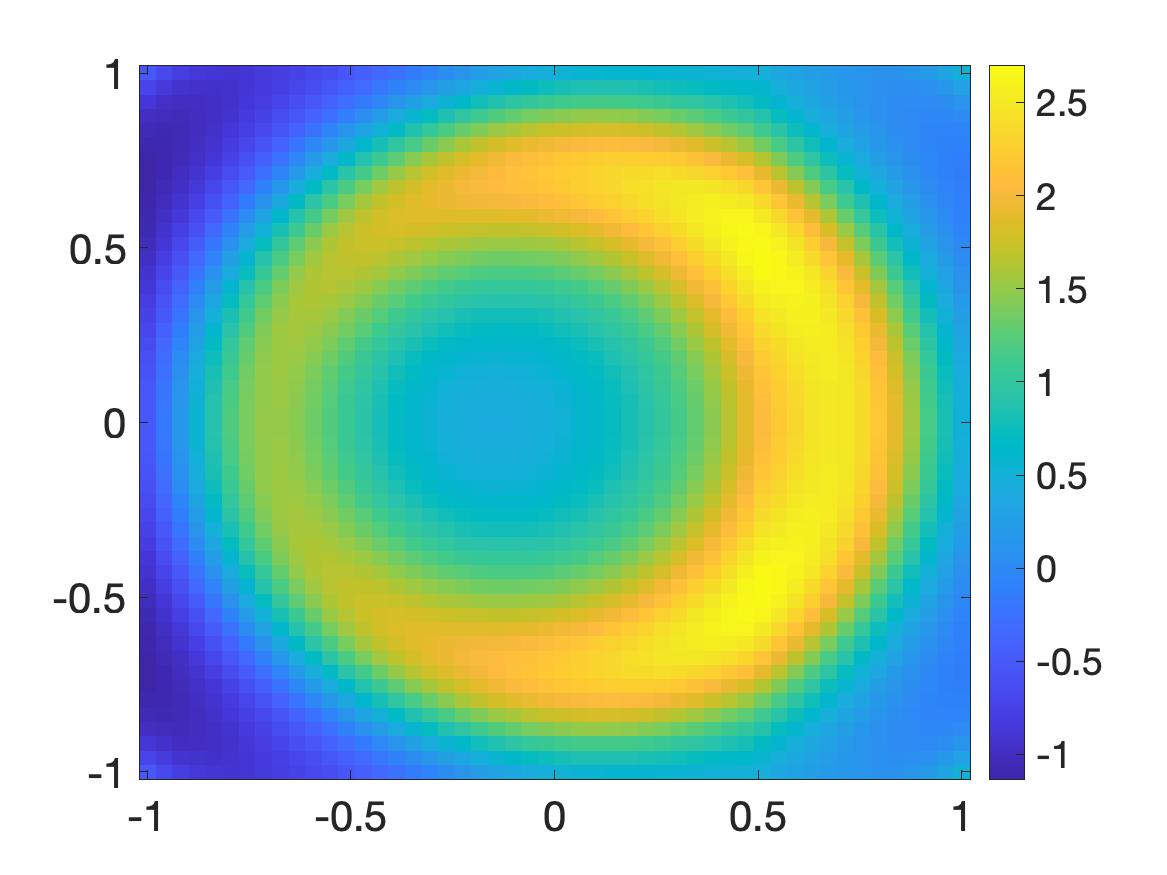}}
		\quad
		\subfloat[The solution $u_{\rm comp}$, computed from 5\% noisy boundary data.]{\includegraphics[width=.3\textwidth]{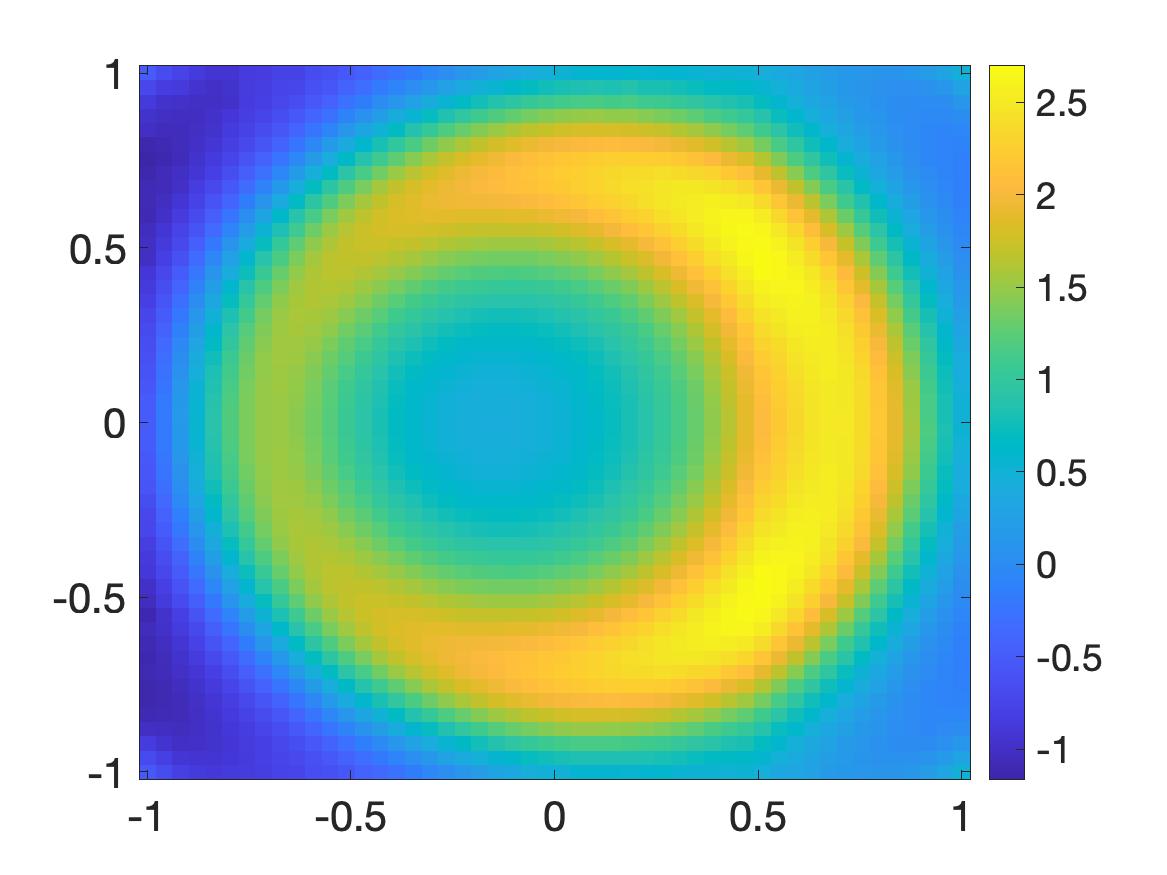}}
		\quad
		\subfloat[The solution $u_{\rm comp}$, computed from 10\% noisy boundary data.]{\includegraphics[width=.3\textwidth]{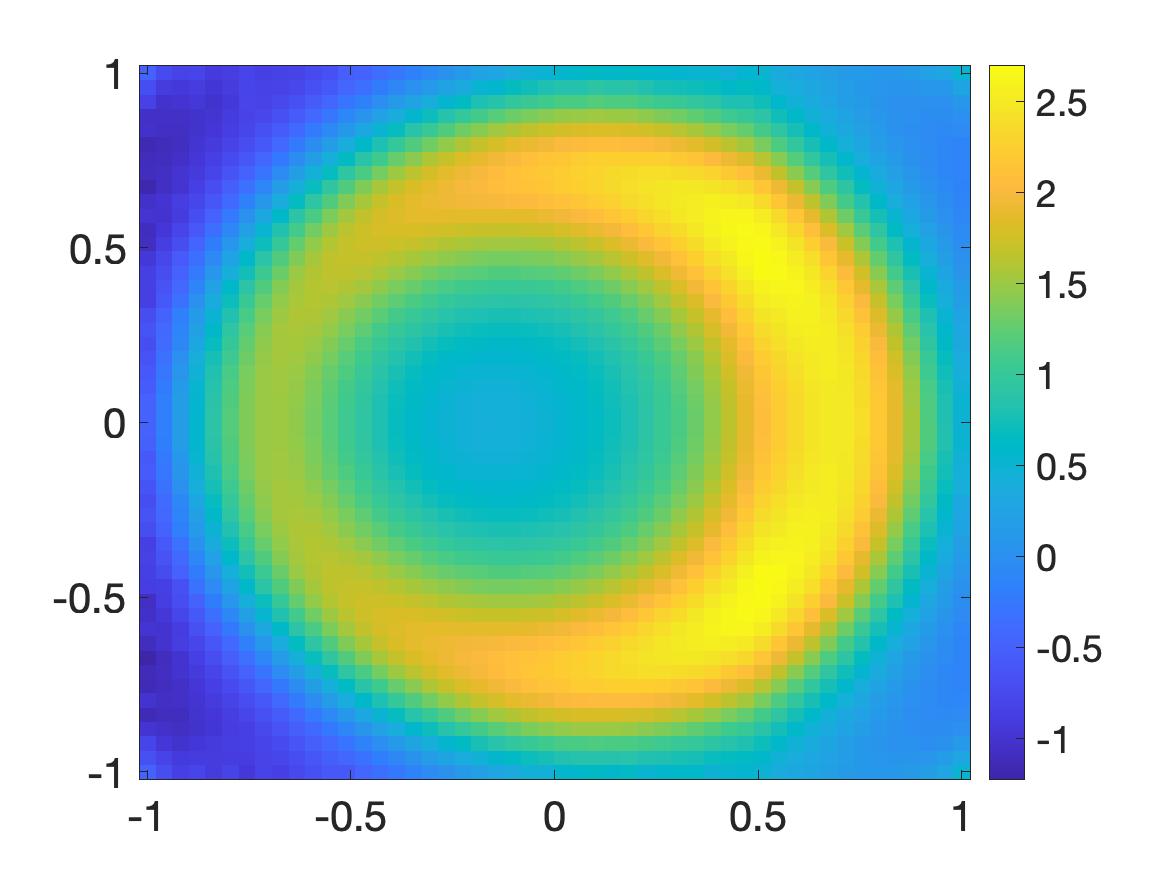}}

		\subfloat[{The relative error $\frac{|u_{\rm comp} - u_{\rm true}|}{\|u_{\rm true}\|_{L^{\infty}}}$, $\delta = 0\%$.}]{\includegraphics[width=.3\textwidth]{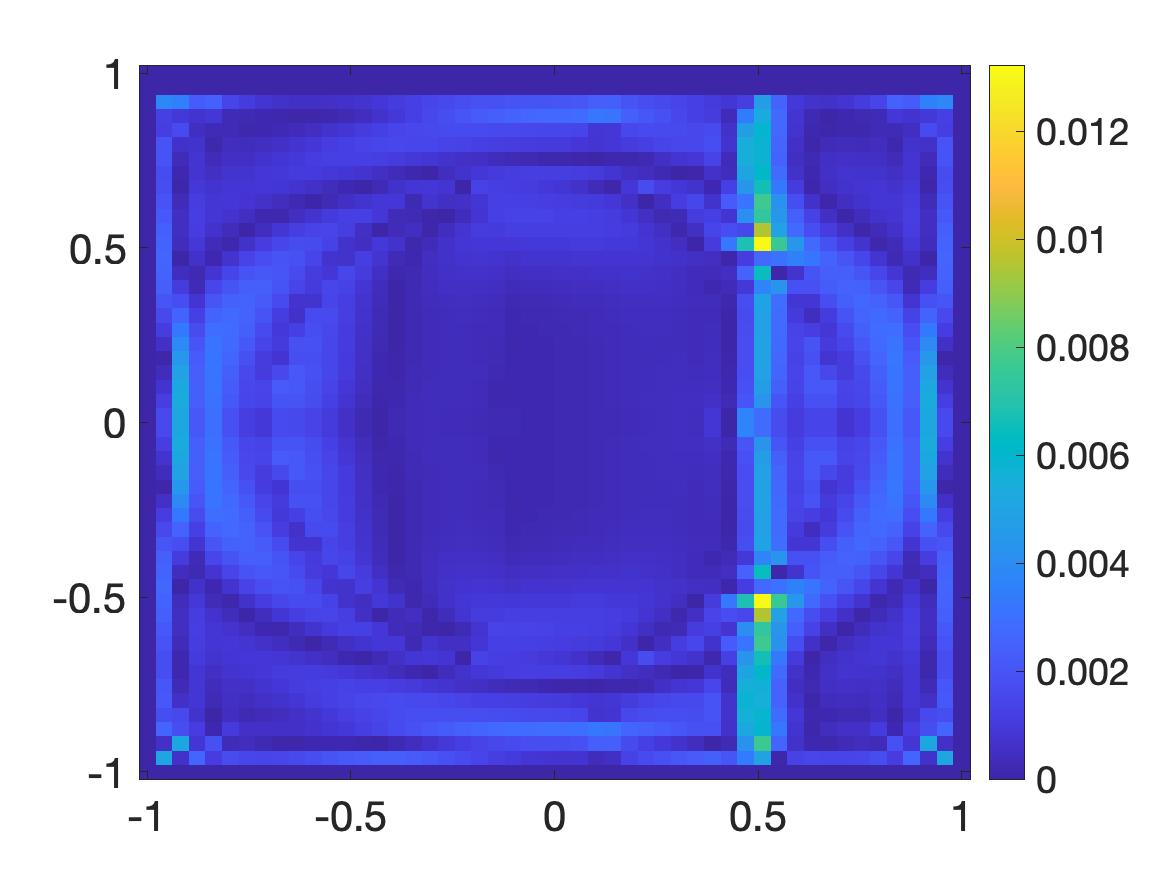}}
		\quad
		\subfloat[[{The relative error $\frac{|u_{\rm comp} - u_{\rm true}|}{\|u_{\rm true}\|_{L^{\infty}}}$, $\delta = 5\%$.}]{\includegraphics[width=.3\textwidth]{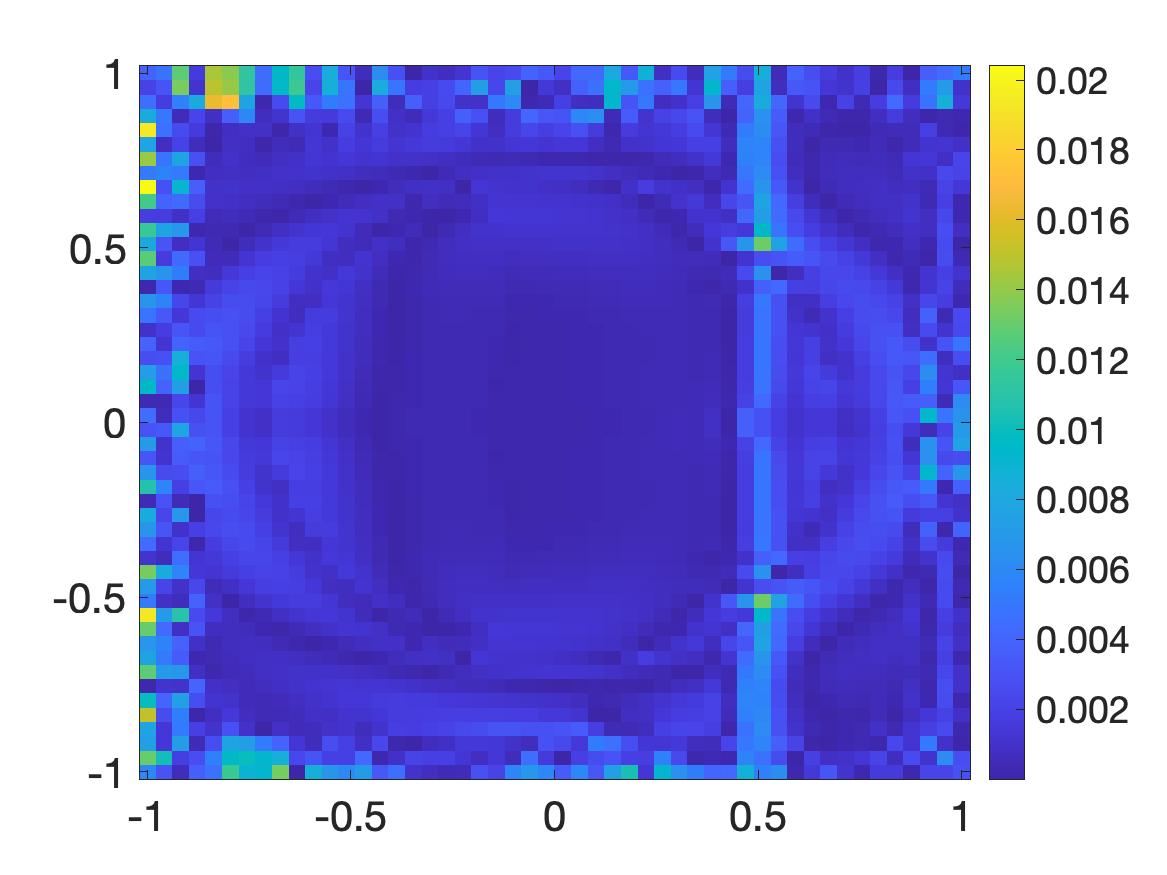}}
		\quad
		\subfloat[[{The relative error $\frac{|u_{\rm comp} - u_{\rm true}|}{\|u_{\rm true}\|_{L^{\infty}}}$, $\delta = 10\%$.}]{\includegraphics[width=.3\textwidth]{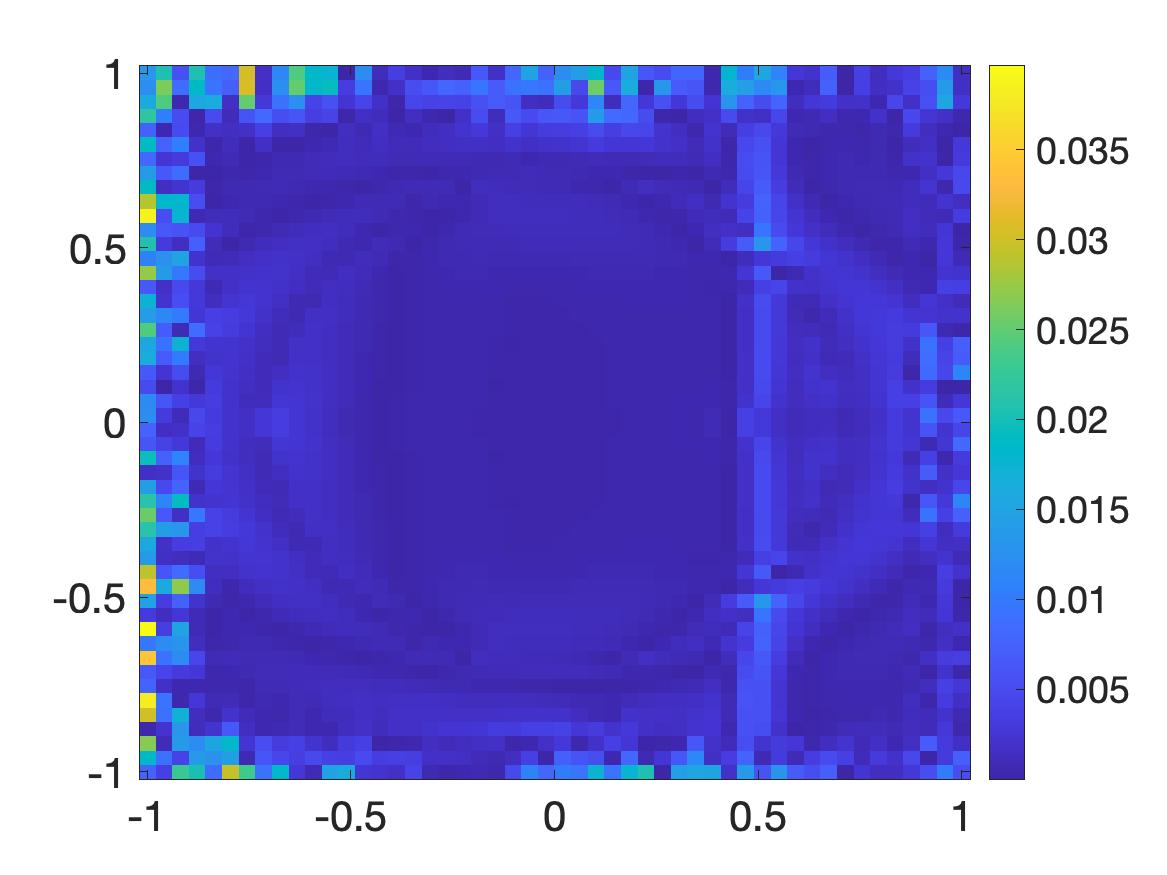}}
		\caption{\label{fig 3} Test 3. The true and computed viscosity solutions with $\delta$ is $0\%$,  5\%, and 10\% noisy boundary Dirichlet data on $\partial \Omega$ and Neumann data on $\Gamma^+$.}
\end{figure}

 The numerical results  are given in Figure \ref{fig 3}.
As mentioned, this test is interesting since the function $F$, see \eqref{5.9999}, is nonconvex, noncoercive, and discontinuous.
{Solving the Hamilton-Jacobi equation with this Hamiltonian is challenging.}
Some existing methods  might not be applicable.
In contrast, the numerical results in Figure \ref{fig 3} are  out of expectation. 
The errors of computation are small, see Table \ref{tab3},  although the solution has complicated structure.
This kind of nonconvex Hamiltonian occurs in the context of two-player zero-sum differential games (see \cite{BCD, Tran19}).

\begin{table}[h!]
\caption{\label{tab3}{Test 3. The performance of the convexification method. The computational time is the time for a Precisions Workstations T7810 with 24 cores to compute the solution $u_{\rm comp}$. In this table, 
the relative $L^\infty(\Omega)$ error is $\|u_{\rm comp} - u_{\rm true}\|_{L^\infty(\Omega)}/\|u_{\rm true}\|_{L^\infty(\Omega)}$.}} 
{
\centering
\begin{tabular}{c c c  c }
\hline\hline
Noise level &computational time&  number of iterations &  relative $ L^\infty(\Omega)$ error\\ 
\hline
0\%& 10.72 minutes & 129 & 1.32\%\\
5\%& 11.12 minutes &131 & 2.04\%\\
10\%& 8.43 minutes & 101 & 3.97\%\\
\hline
\end{tabular}
}
\end{table}

\begin{remark}
In Tests 1, 2 and 3 above, we are in the context that the knowledge of $u_z$ on $\Gamma^+$ can be computed from the knowledge of $u$ and the form of the Hamilton-Jacobi equation, see Assumption \ref{assumption 1} and Remark \ref{rm 1}.
However, if the given Hamilton-Jacobi equation is rather complicated as in Tests 4 and 5  below, solving $u_z|_{\Gamma^+}$ from $u|_{\Gamma^+}$ is impossible. 
In this case, we minimize $J_{\lambda, \beta, \eta}$ on the set $\{u \in H^p(\Omega): u(\x) = f(\x) \mbox{ for } \x \in \partial \Omega\}$. 
The convexification method still gives us out of expectation numerical results in these two tests. 
However, the proofs of the convexification theorem and the convergence of the numerical scheme for the problem with only Dirichlet boundary condition \eqref{HJ}--\eqref{dir} are extremely challenging, and they are out of the scope of this paper. 
\end{remark}

\medskip

\begin{figure}[h!]
		\subfloat[The true solution $u_{\rm true} = -|x| - 1$.]{\includegraphics[width=.3\textwidth]{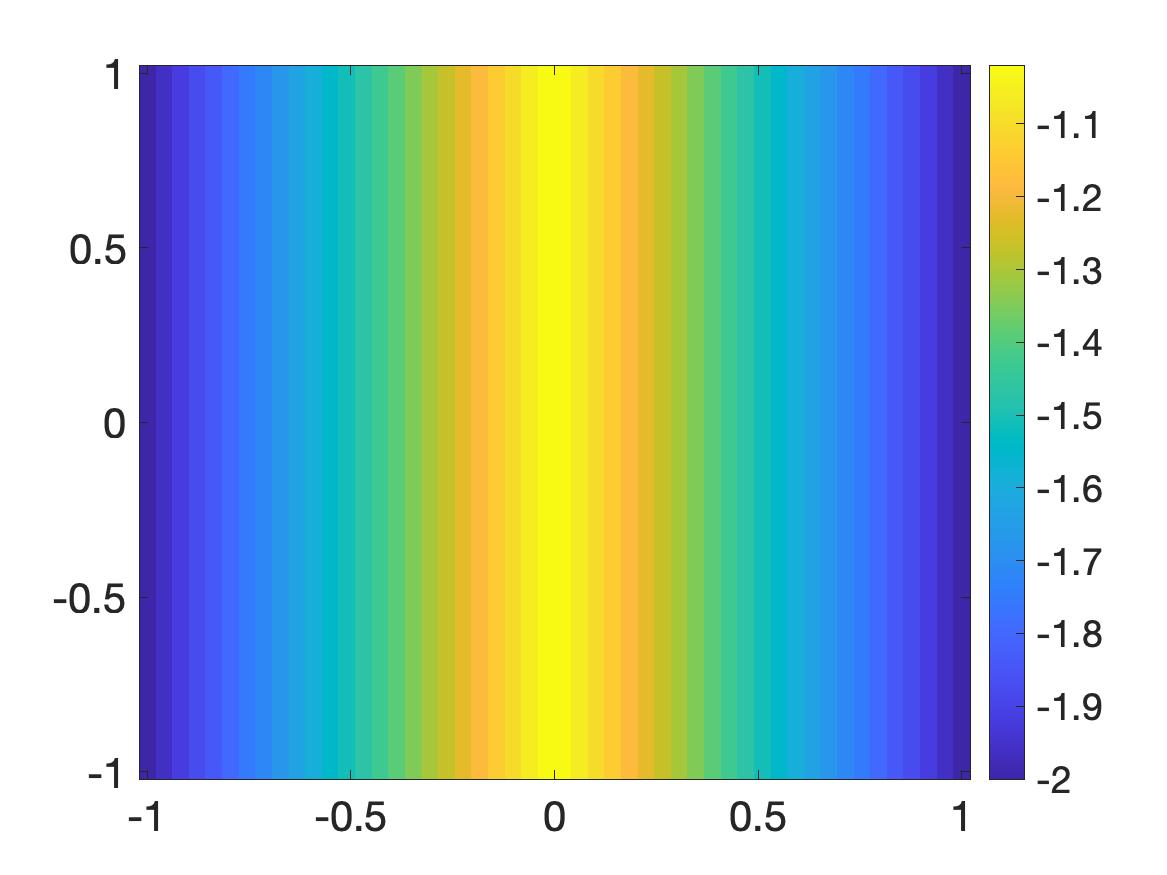}}
		
		\subfloat[The solution $u_{\rm comp}$, computed from noiseless boundary data.]{\includegraphics[width=.3\textwidth]{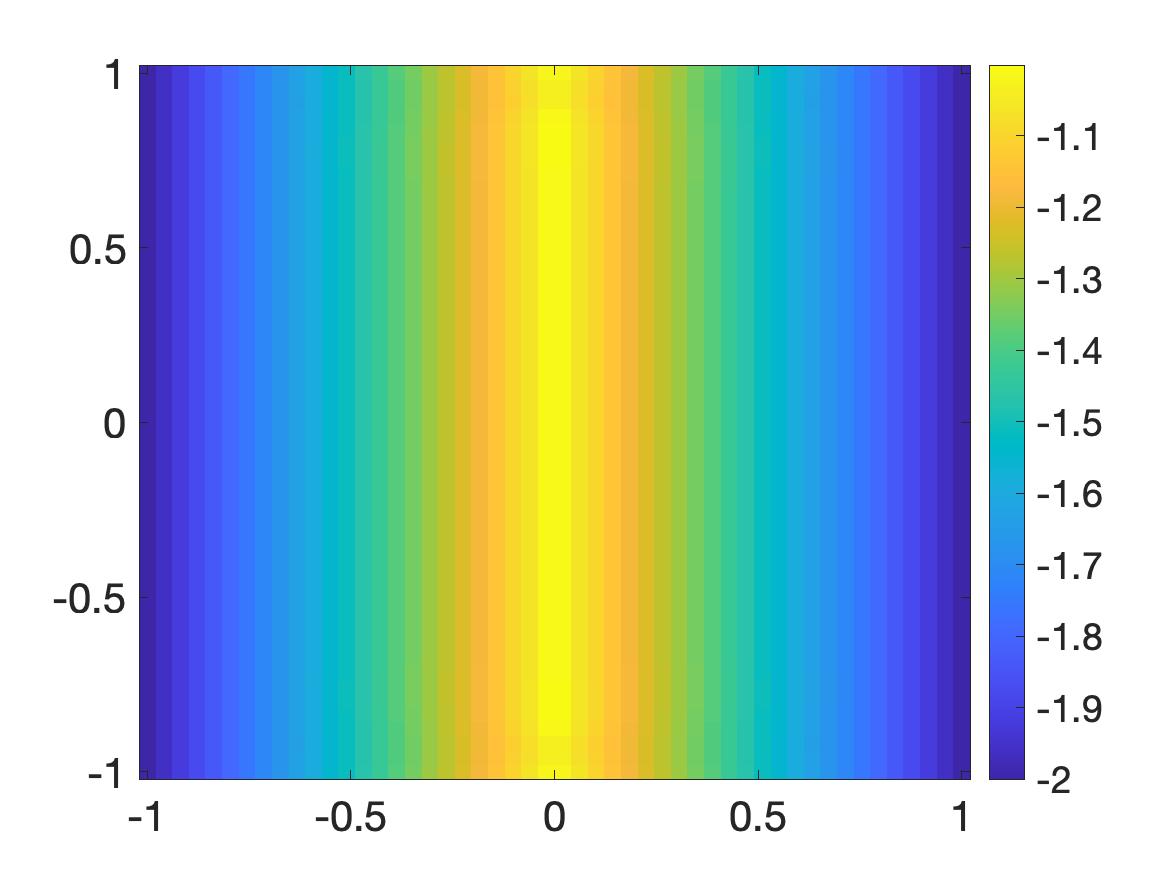}}
		\quad
		\subfloat[The solution $u_{\rm comp}$, computed from 5\% noisy boundary data.]{\includegraphics[width=.3\textwidth]{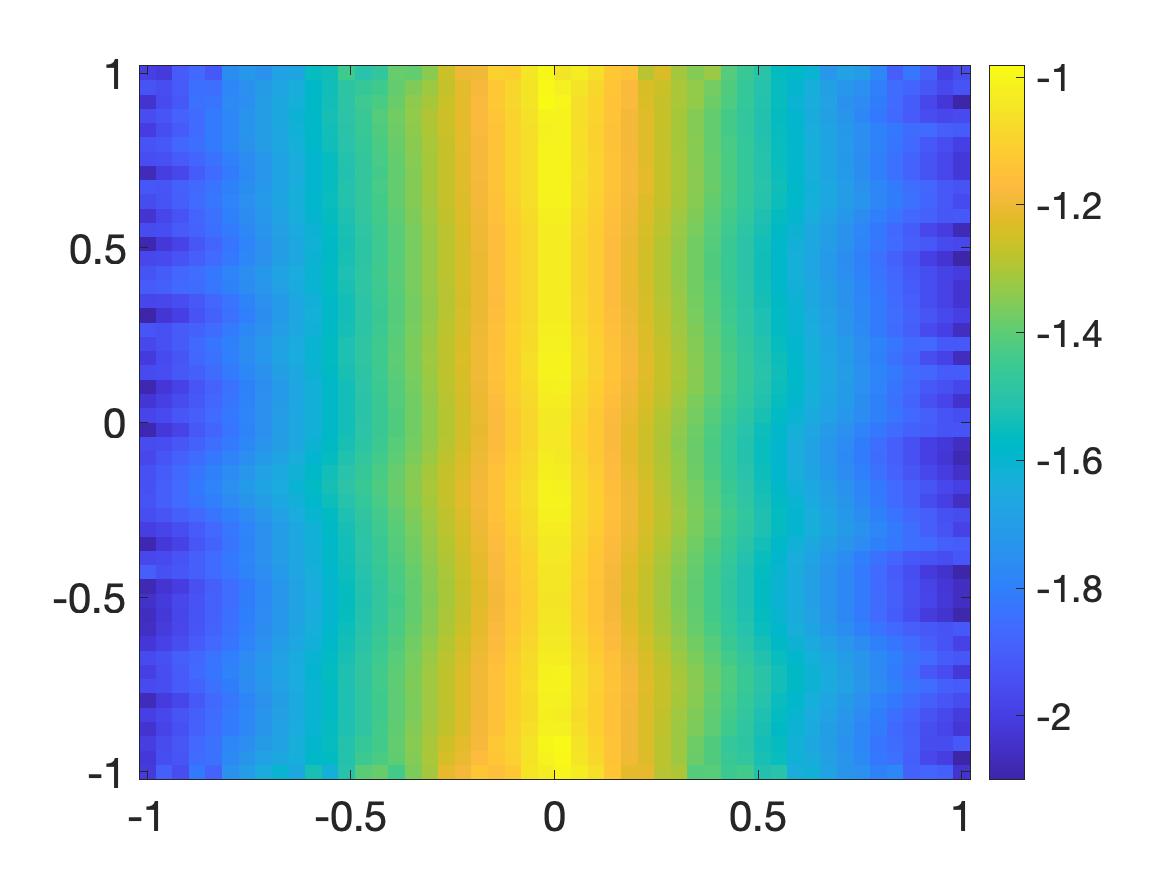}}
		\quad
		\subfloat[The solution $u_{\rm comp}$, computed from 10\% noisy boundary data.]{\includegraphics[width=.3\textwidth]{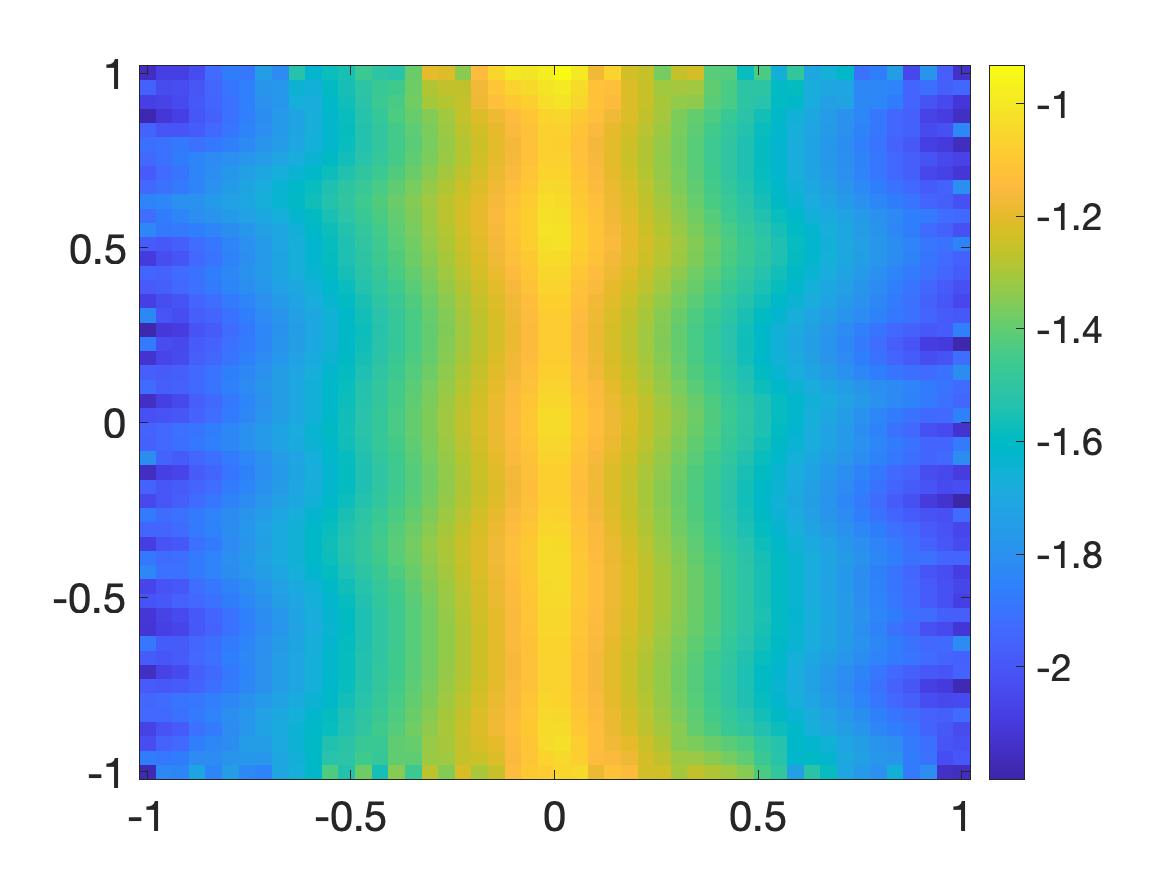}}

		\subfloat[{The relative error $\frac{|u_{\rm comp} - u_{\rm true}|}{\|u_{\rm true}\|_{L^{\infty}}}$, $\delta = 0\%$.}]{\includegraphics[width=.3\textwidth]{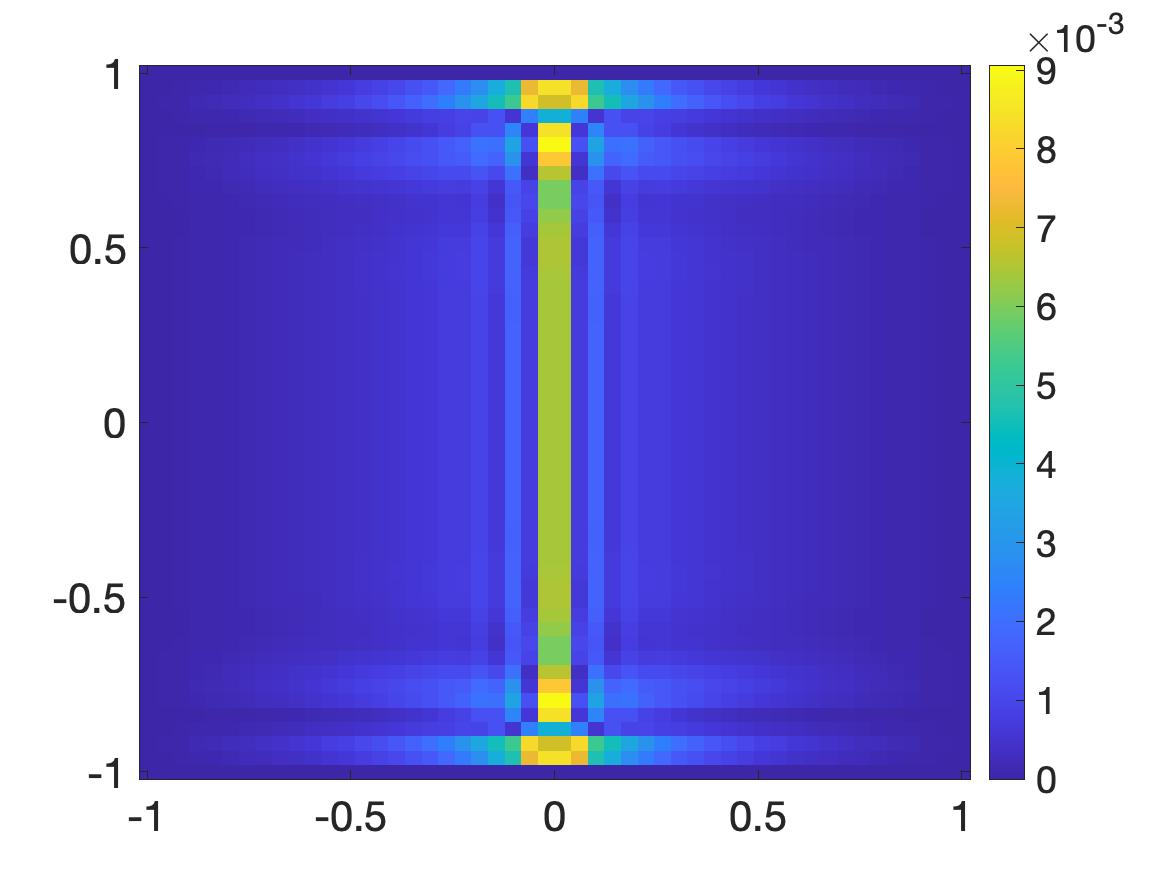}}
		\quad
		\subfloat[[{The relative error $\frac{|u_{\rm comp} - u_{\rm true}|}{\|u_{\rm true}\|_{L^{\infty}}}$, $\delta = 5\%$.}]{\includegraphics[width=.3\textwidth]{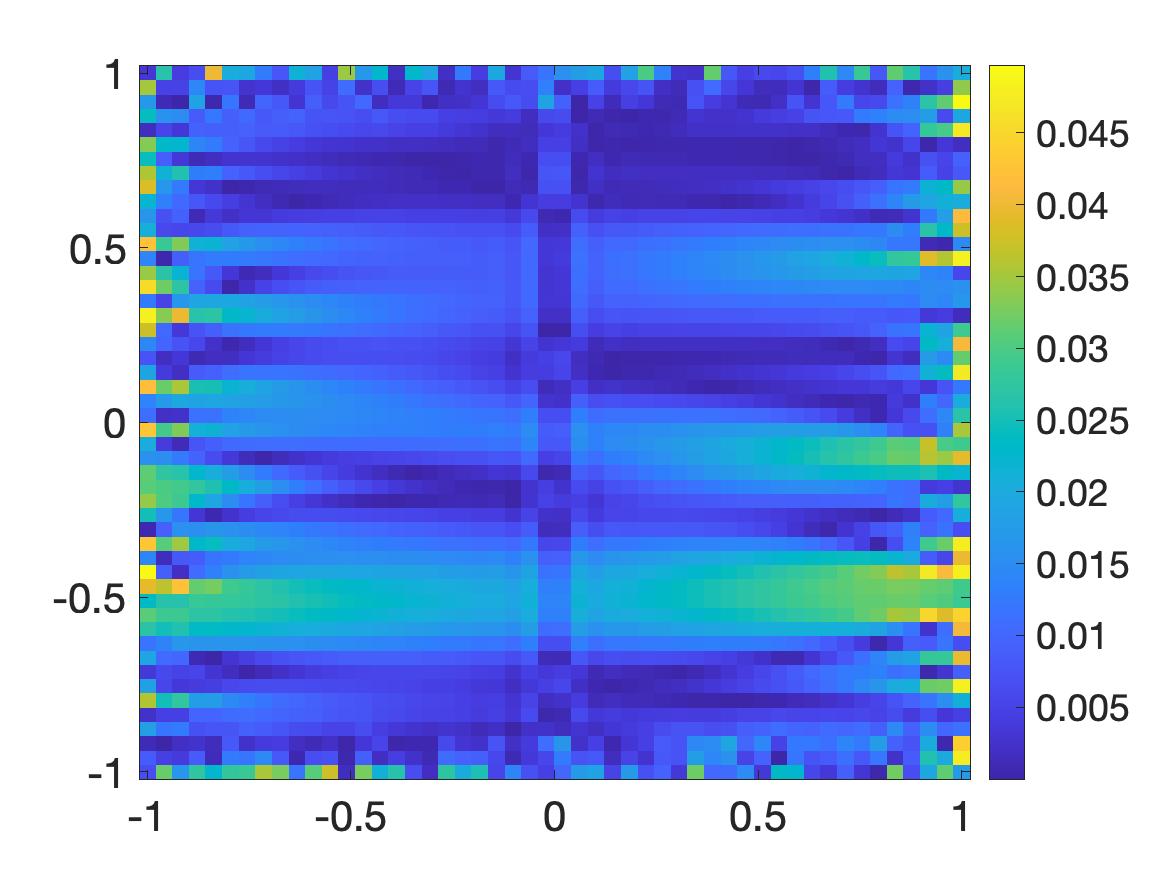}}
		\quad
		\subfloat[[{The relative error $\frac{|u_{\rm comp} - u_{\rm true}|}{\|u_{\rm true}\|_{L^{\infty}}}$, $\delta = 10\%$.}]{\includegraphics[width=.3\textwidth]{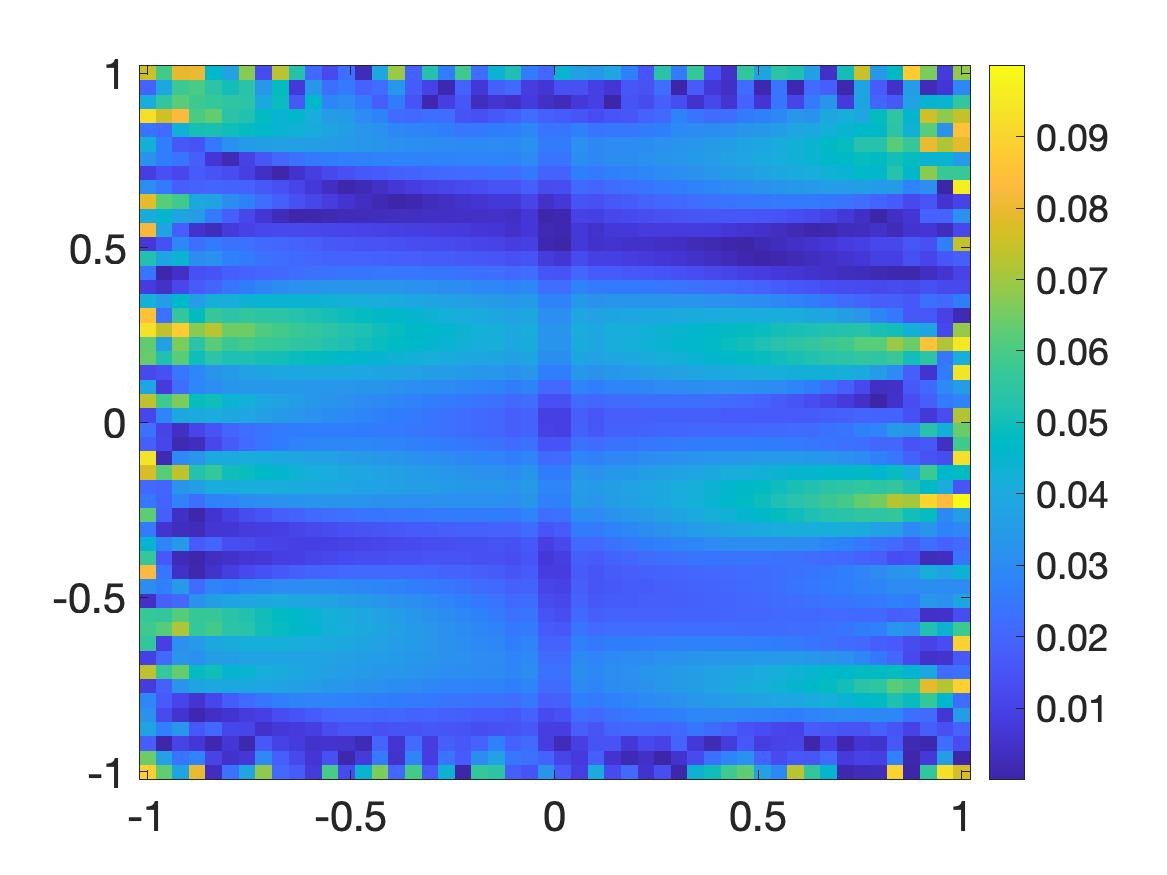}}

		\caption{\label{fig 5} Test 4. The true and computed viscosity solutions with $\delta$ is $0\%$,  5\% and 10\% noisy boundary Dirichlet data on $\partial \Omega$.}
\end{figure}

{\it Test 4.} We next consider the G-equation, which arises from instantaneous flame position. 
We solve \eqref{HJ}--\eqref{dir} when
\begin{equation}
	F(\x, s, \p) = s + |\p| - x p_1 
	\label{F4}
\end{equation}
for $\x \in \Omega, s \in \R, \p = (p_1, p_2)\in \R^2$
and the boundary data is given by
\begin{equation}
	f(\x) = -|x| - 1 \quad \x = (x, z) \in \partial \Omega.
	\label{f4}
\end{equation}
The function $u_{\rm true}(x, z) = -|x| - 1$ is the true viscosity solution for this test.
The verification of this is similar to that of Tests 2 and 3, and is hence omitted here.
 The numerical results  are given in Figure \ref{fig 5}. Relative errors are $0.91\%$, $4.97\%$ and $9.99\%$ when $\delta$ is 0\%, 5\% and $10\%$ respectively {(see Table \ref{tab4})}.
 The G-equation  is quite popular in the combustion science literature.
We refer the readers to \cite{CarNoSou, XinYu, LiuXinYu} for some recent important mathematical developments.

\begin{table}[h!]
\caption{\label{tab4}{Test 4. The performance of the convexification method. The computational time is the time for a Precisions Workstations T7810 with 24 cores to compute the solution $u_{\rm comp}$. 
The relative $L^\infty(\Omega)$ error is $\|u_{\rm comp} - u_{\rm true}\|_{L^\infty(\Omega)}/\|u_{\rm true}\|_{L^\infty(\Omega)}$.}} 
{
\centering
\begin{tabular}{c c c  c }
\hline\hline
Noise level &computational time&  number of iterations &  relative $ L^\infty(\Omega)$ error\\ 
\hline
0\%& 22.62 minutes & 280 &0.91\% \\
5\%& 24.99 minutes & 311 & 4.97\%\\
10\%&26.50 minutes & 322& 9.99\%\\
\hline
\end{tabular}
}
\end{table}

\medskip

{\it Test 5.} We finally consider a quite complicated form of the function $F$. 
For $\x = (x, z) \in \Omega$, let 
\[
	G(x, z) = \left\{
		\begin{array}{ll}
			1 + 2\pi \cos(2\pi(x + z)) &x > 0,\\
			-1 + 2\pi \cos(2\pi(x + z)) &x < 0.
		\end{array}
	\right.
\]

\begin{figure}[h!]
		\subfloat[The true solution $u_{\rm true} = -|x| + \sin(2\pi (x + z))$.]{\includegraphics[width=.3\textwidth]{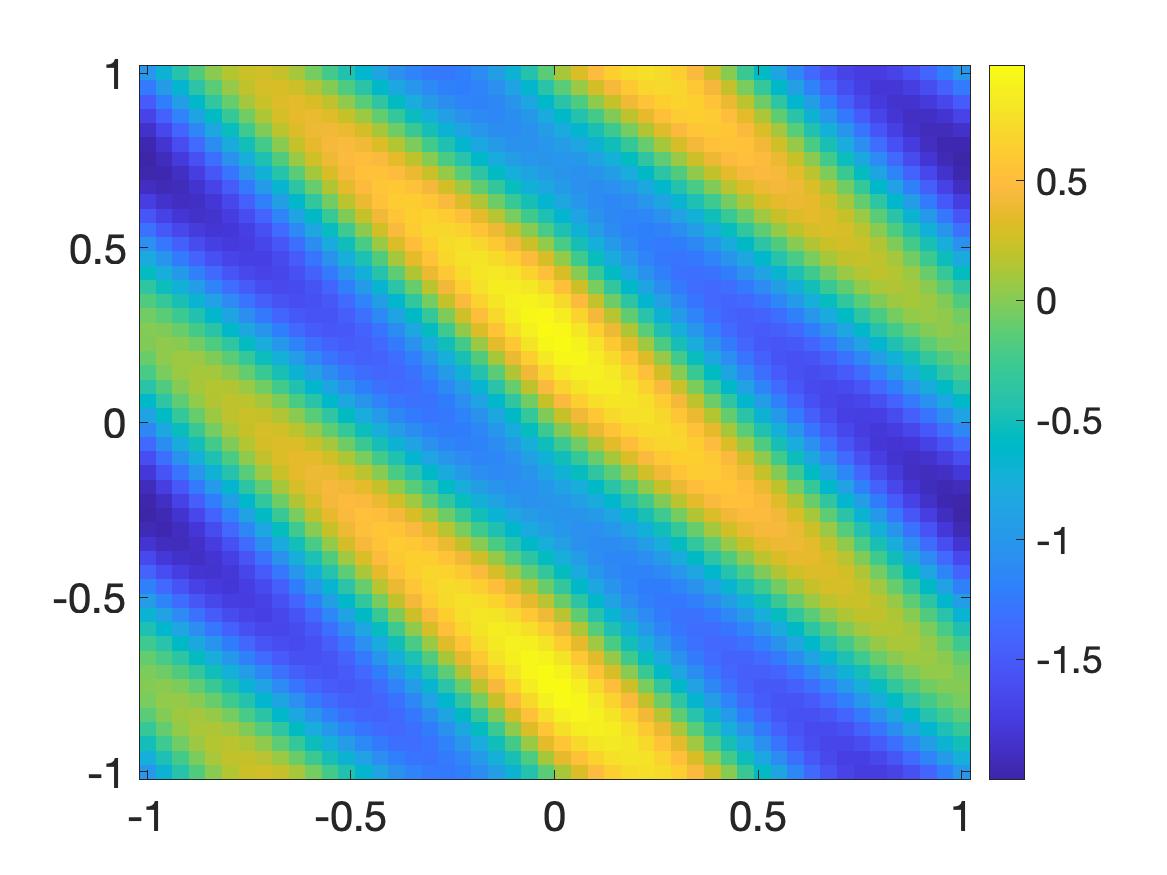}}
		
		\subfloat[The solution $u_{\rm comp}$, computed from noiseless boundary data.]{\includegraphics[width=.3\textwidth]{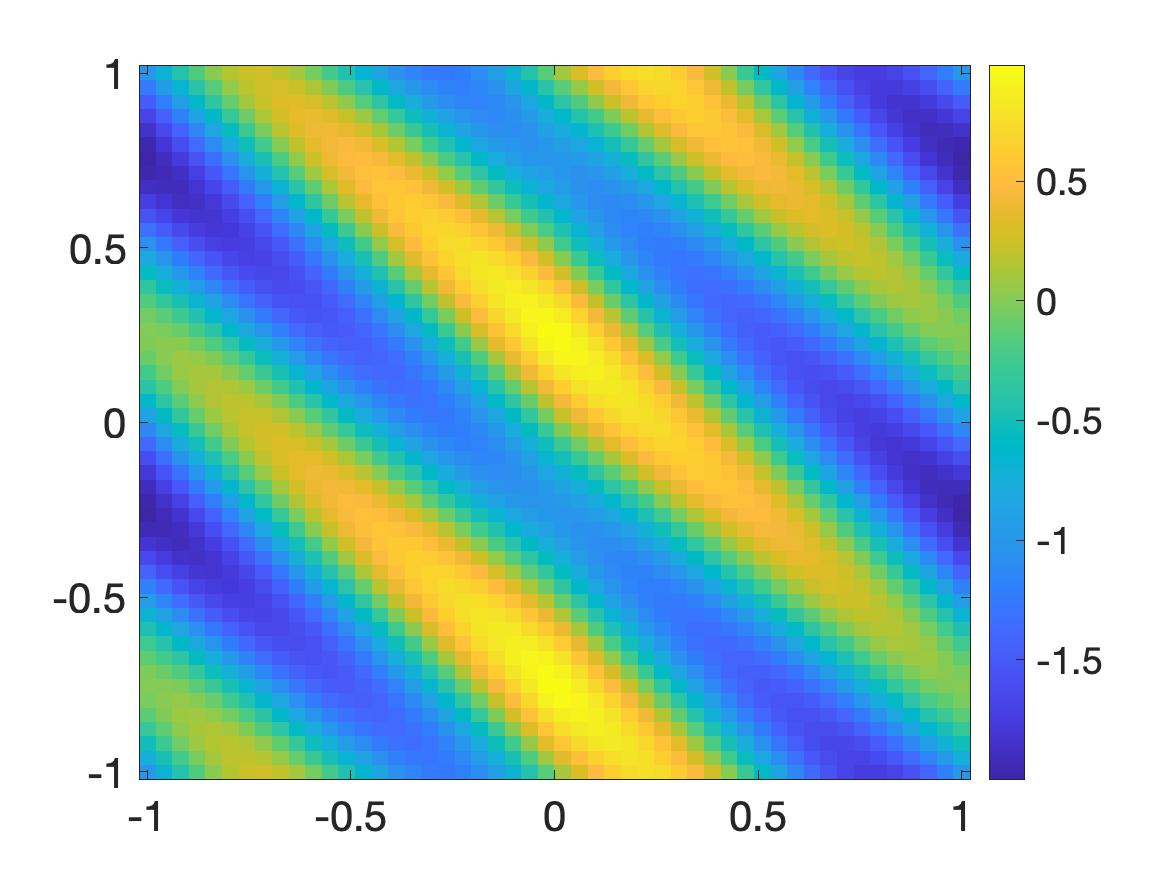}}
		\quad
		\subfloat[The solution $u_{\rm comp}$, computed from 5\% noisy boundary data.]{\includegraphics[width=.3\textwidth]{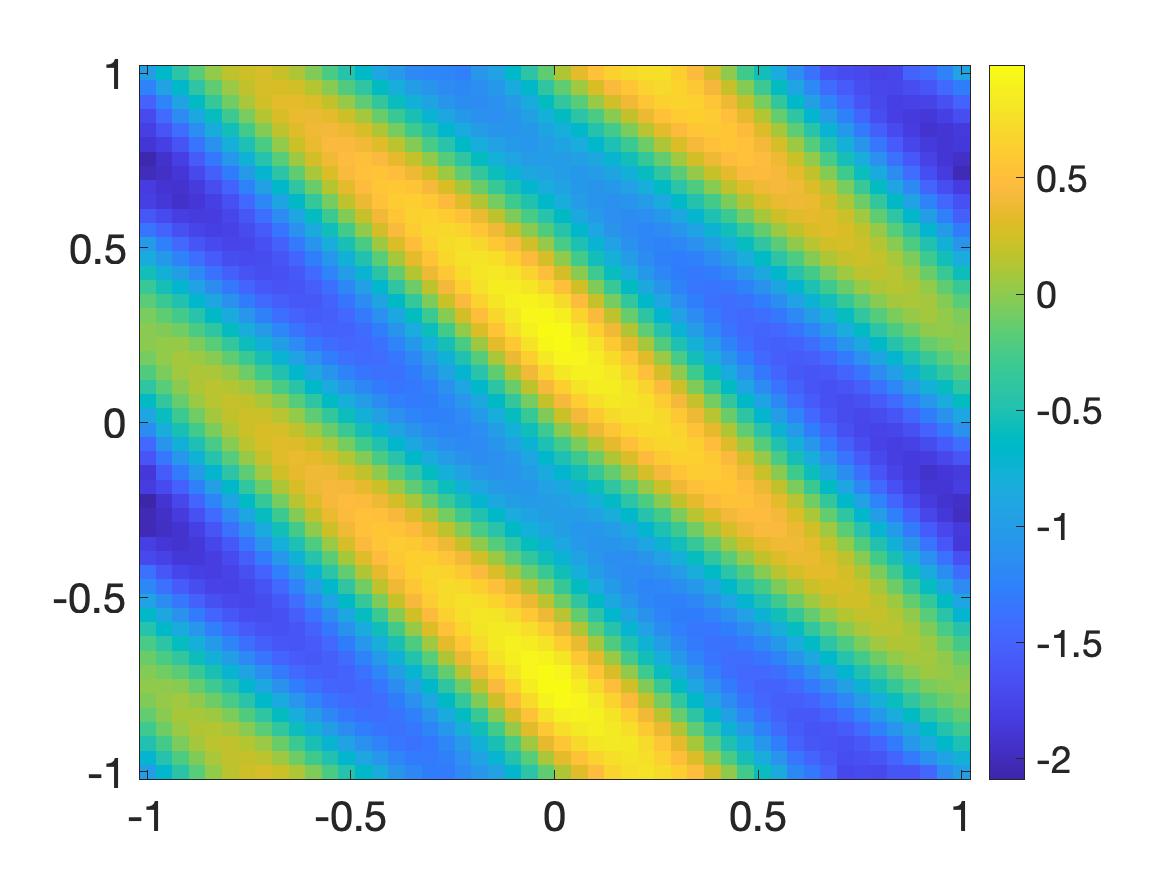}}
		\quad
		\subfloat[The solution $u_{\rm comp}$, computed from 10\% noisy boundary data.]{\includegraphics[width=.3\textwidth]{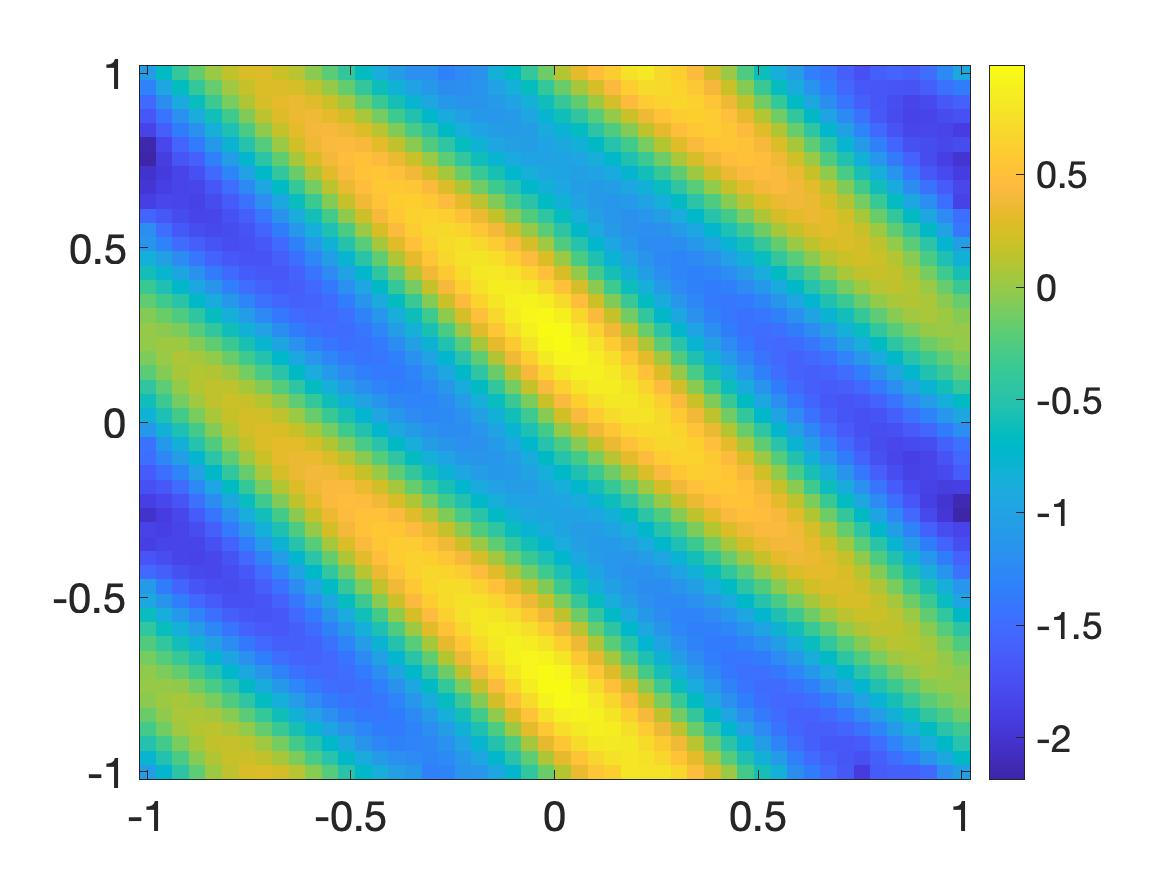}}
		
		\subfloat[{The relative error $\frac{|u_{\rm comp} - u_{\rm true}|}{\|u_{\rm true}\|_{L^{\infty}}}$, $\delta = 0\%$.}]{\includegraphics[width=.3\textwidth]{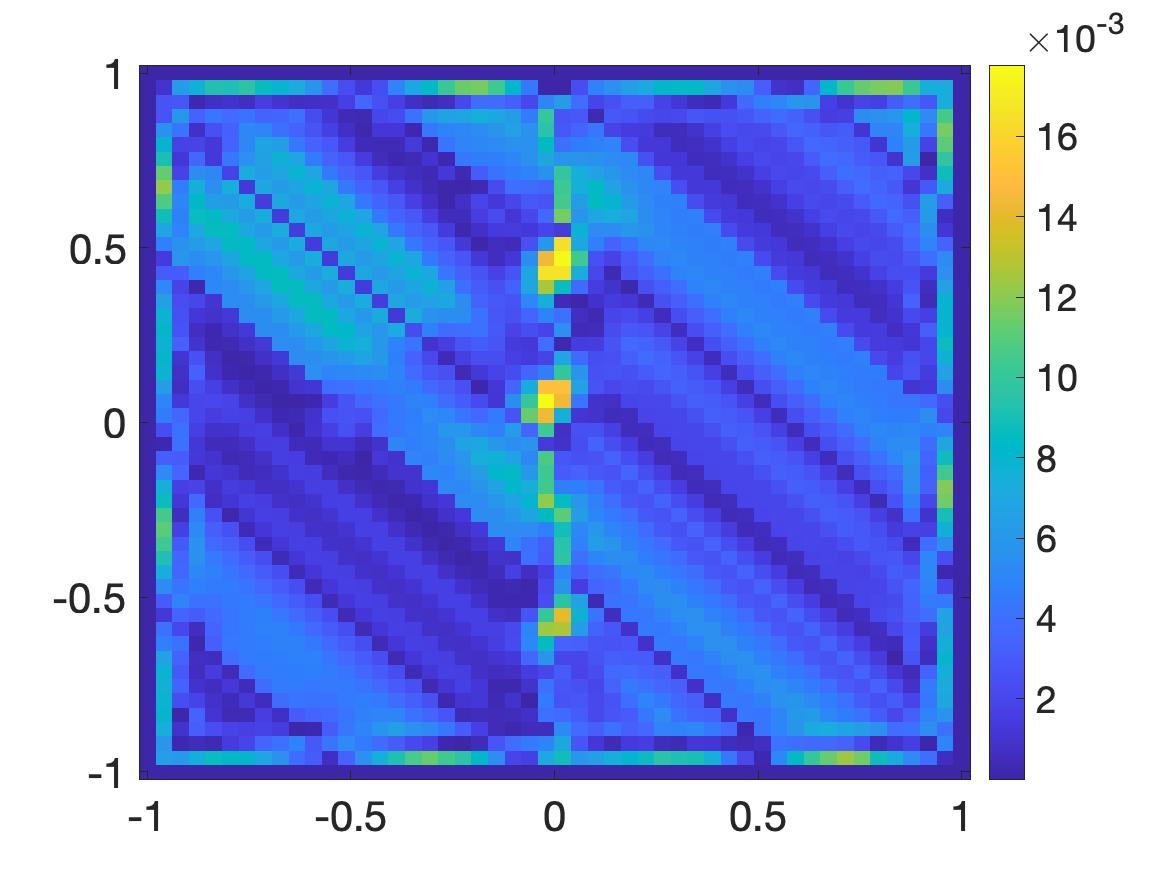}}
		\quad
		\subfloat[[{The relative error $\frac{|u_{\rm comp} - u_{\rm true}|}{\|u_{\rm true}\|_{L^{\infty}}}$, $\delta = 5\%$.}]{\includegraphics[width=.3\textwidth]{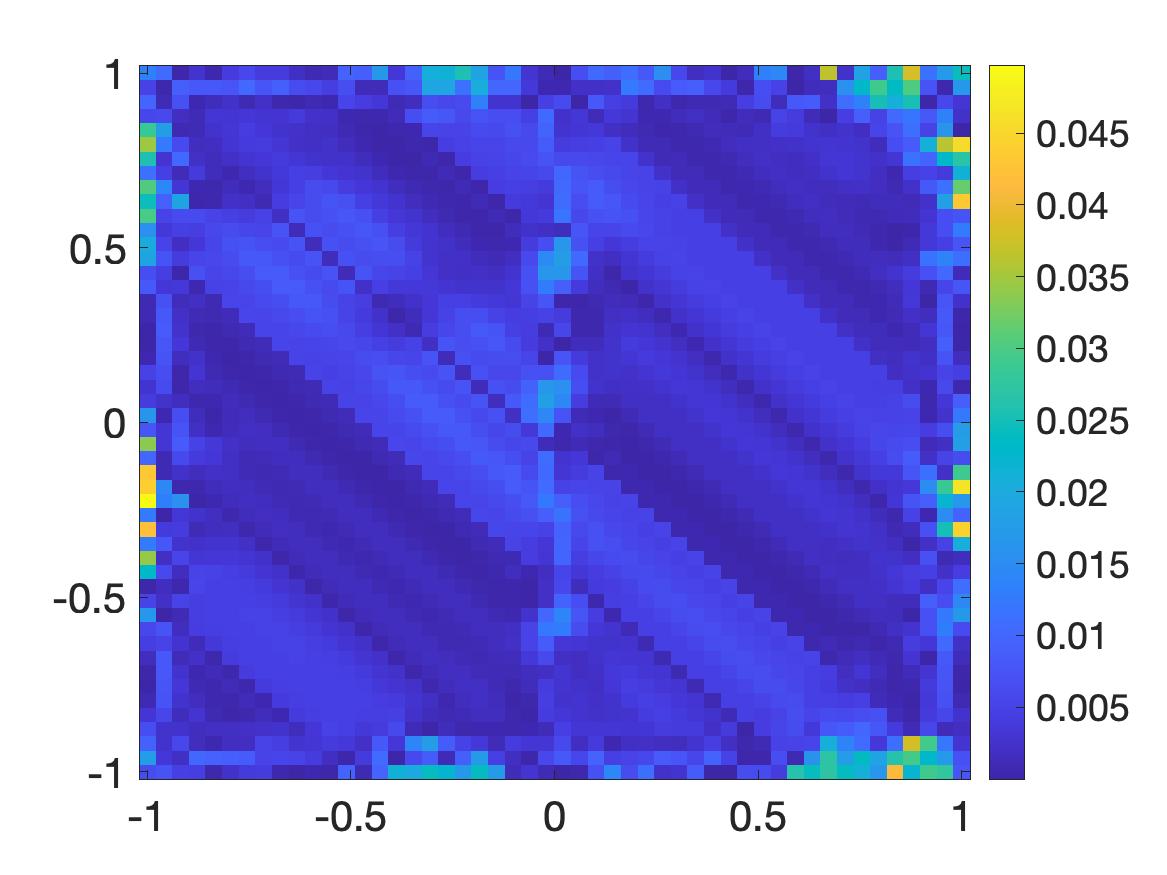}}
		\quad
		\subfloat[[{The relative error $\frac{|u_{\rm comp} - u_{\rm true}|}{\|u_{\rm true}\|_{L^{\infty}}}$, $\delta = 10\%$.}]{\includegraphics[width=.3\textwidth]{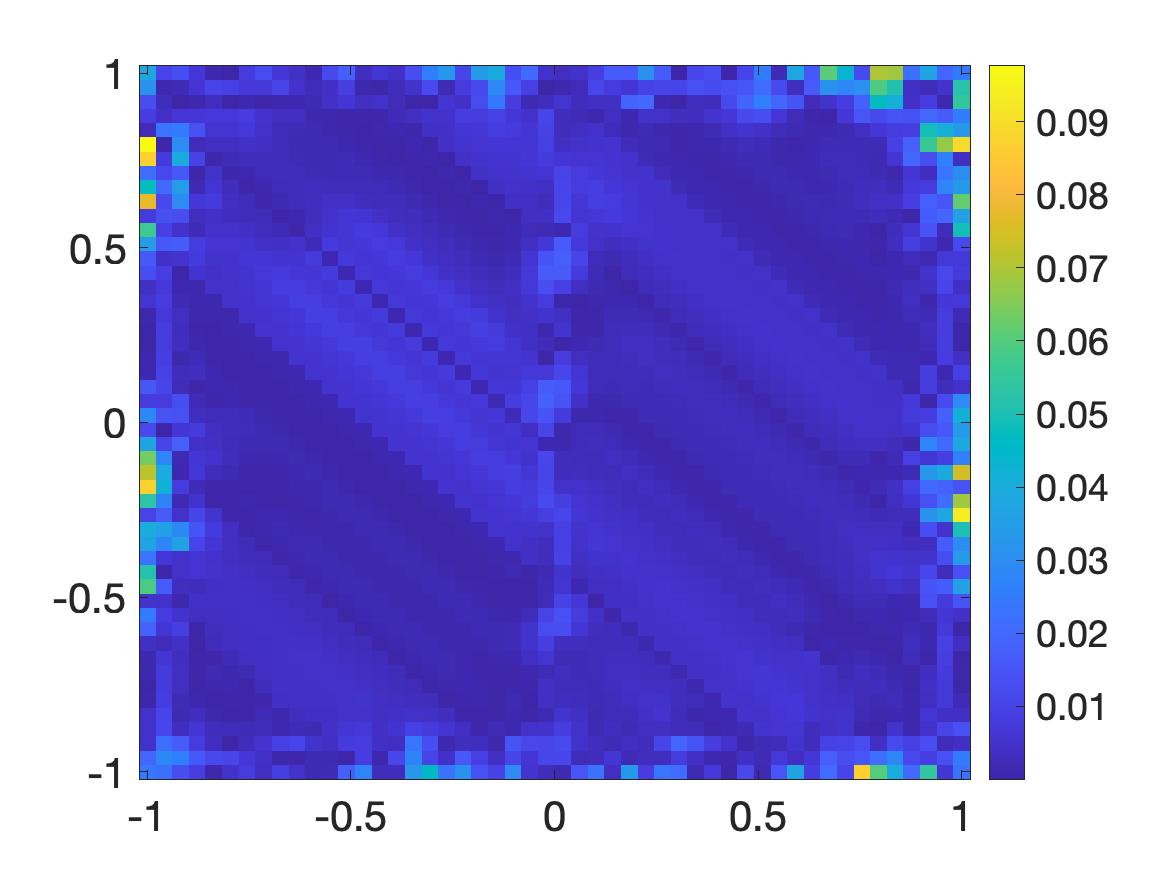}}

		\caption{\label{fig 6} Test 5. The true and computed viscosity solutions with $\delta$ is $0\%$,  5\%, and 10\% noisy boundary Dirichlet data on $\partial \Omega$. {The relative errors with these noise levels are given in Table \ref{tab5}.}}
\end{figure}

We solve \eqref{HJ}--\eqref{dir} when
\begin{multline}
	F(\x, s, \p) = 15s + \min \left\{ |\p|, ||\p|-10|+6 \right\}
	- \Big[15(-|x| + \sin(2\pi (x + z))) 
	\\
	+\min\left\{ |\sqrt{|G(x, z)|^2 + 4\pi^2\cos(2\pi (x + z))}|, \big|\sqrt{|G(x, z)|^2 + 4\pi^2\cos(2\pi (x + z))}-10\big|+6  \right\} \Big]
	\label{F6}
\end{multline}
for $\x \in \Omega, s \in \R, \p \in \R^2$,
and the Dirichlet boundary data is given by
\begin{equation}
	f(\x) = -|x| + \sin(2\pi (x + z)) \quad \text{ for all } \x = (x, z) \in \partial \Omega.
	\label{f6}
\end{equation}
It it worth mentioning that $F, G$ are not continuous in general at $x=0$.
The function $u_{\rm true}(x, z) = -|x| + \sin(2\pi (x + z))$ is the true viscosity solution for this test. 
Its graph and the graphs of the computed solutions are displayed in Figure \ref{fig 6}. 
Let us give a rigorous verification here.
If $x \neq 0$, then $u_{\rm true}$ is differentiable at $(x, z)$, and
\[
\nabla u_{\rm true}(x, z)  = \left(-\frac{x}{|x|}+2\pi \cos(2\pi(x+z)), 2\pi \cos(2\pi(x+z)) \right),
\]
which implies that $F(\x,u(\x),\nabla u_{\rm true}(\x))=0$.
If $x=0$, then $u_{\rm true}$ is not differentiable at $\x=(x,z)$.
We can only find smooth test functions that touch $u_{\rm true}$ from above at $\x$, and we cannot find smooth test functions that touch $u_{\rm true}$ from below at $\x$.
Let $\phi$ be a smooth test function that touches $u_{\rm true}$ from above at $\x$.
Then,
\[
\phi_x(\x) \in 2\pi \cos(2\pi z)+[-1,1], \ \phi_z(\x)=2\pi \cos(2\pi z),
\]
which yields that 
\[
|\nabla \phi(\x)|^2 \leq 1 + 4 \pi + 4 \pi^2 \quad \Rightarrow \quad |\nabla \phi(\x)| \leq 1 + 2\pi <8.
\]
Therefore,
\begin{align*}
F_*(\x,u_{\rm true}(\x),\nabla \phi(\x) ) &= 10 \sin(2\pi z) + |\nabla \phi(\x)|  - 10 \sin(2\pi z)  - \sqrt{|G|^*(0, z) + 4\pi^2\cos(2\pi  z)}\\
&=  |\nabla \phi(\x)|   - \sqrt{1 + 2\pi |\cos(2\pi  z)|+4\pi^2} \leq 0.
\end{align*}
Here, $F_*$ is the lower semicontinuous envelope of $F$, and $|G|^*$ is the upper semicontinuous envelope of $|G|$.
Thus, the subsolution test holds for $u_{\rm true}$ at $\x$.
We imply that $u_{\rm true}$ is a viscosity solution  to \eqref{HJ}--\eqref{dir}.

It is evidently clear that the convexification method delivers satisfactory numerical results.
{The relative errors are provided in Table \ref{tab5}}.
This kind of Hamiltonian was considered in \cite{QianTranYu} in the context of the periodic homogenization theory.

\begin{table}[h!]
\caption{\label{tab5}{Test 5. The performance of the convexification method. The computational time is the time for a Precisions Workstations T7810 with 24 cores to compute the solution $u_{\rm comp}$. In this table, 
the relative $L^\infty(\Omega)$ error is $\|u_{\rm comp} - u_{\rm true}\|_{L^\infty(\Omega)}/\|u_{\rm true}\|_{L^\infty(\Omega)}$.}} 
{
\centering
\begin{tabular}{c c c  c }
\hline\hline
Noise level &computational time&  number of iterations &  relative $ L^\infty(\Omega)$ error\\ 
\hline
0\%& 11.83 minutes &144 & 1.78\%\\
5\%& 10.06 minutes&121 &4.97\%\\
10\%& 11.54 minutes&143 &9.77\%\\
\hline
\end{tabular}
}
\end{table}

\begin{figure}[h!]
\begin{center}
	\subfloat[$\p \mapsto |p_1| - |p_2|$.]{\includegraphics[width=.35\textwidth]{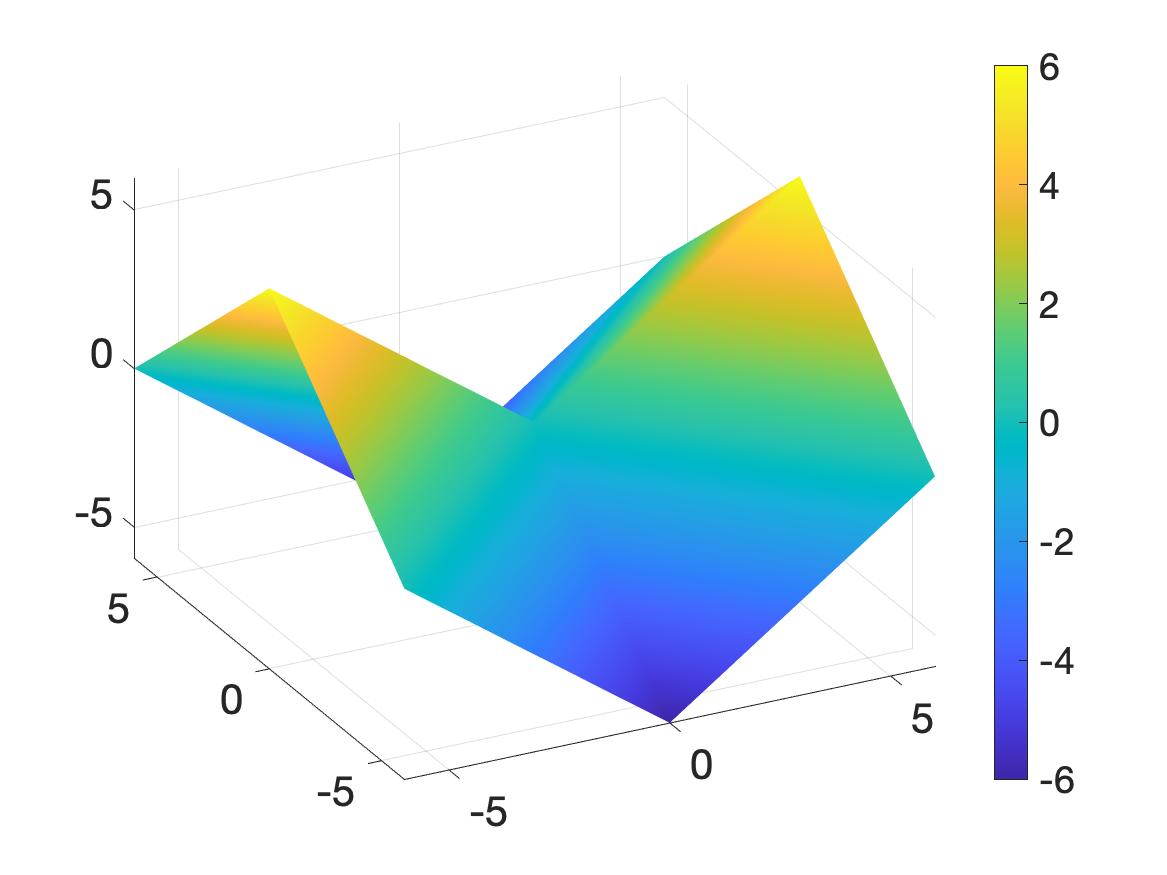}}
		\quad
		\subfloat[$\p \mapsto \min \left\{ |\p|, ||\p|-10|+6 \right\}$]{\includegraphics[width=.35\textwidth]{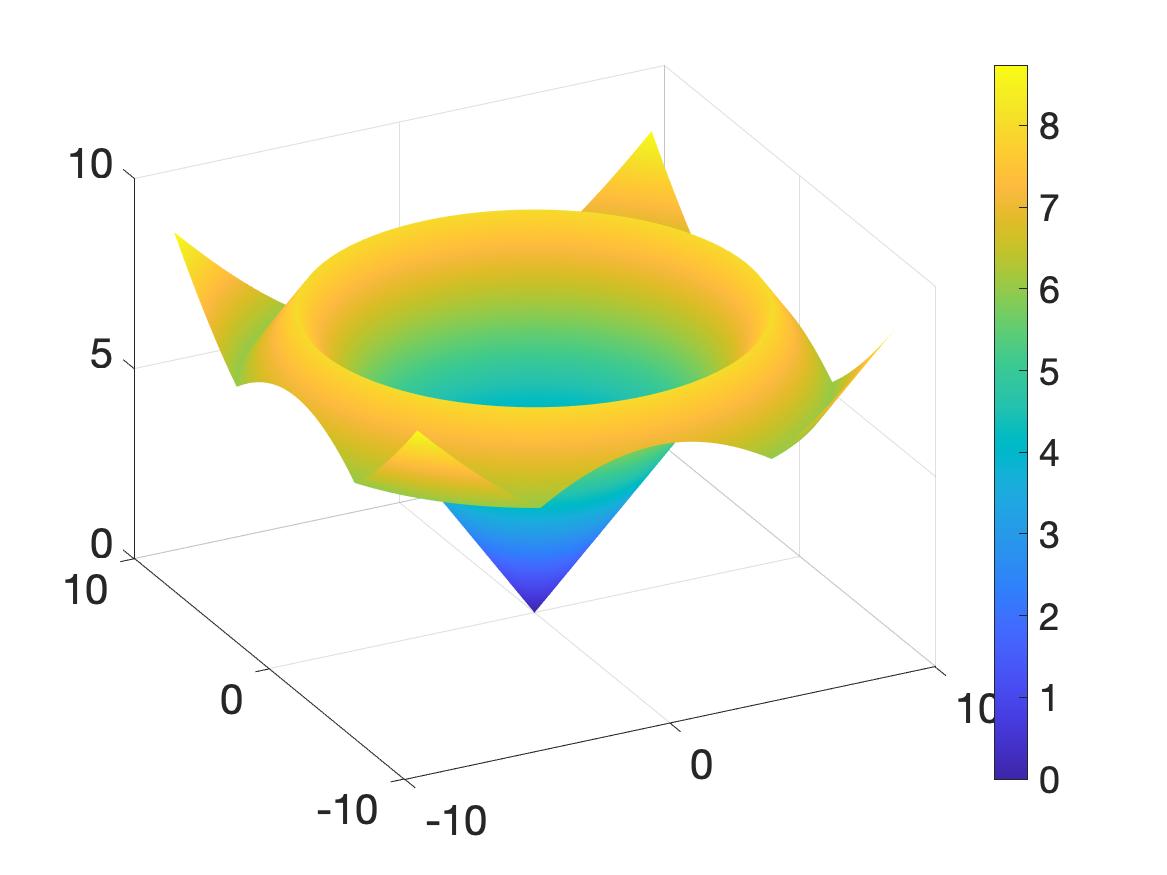}}
	\caption{\label{fig nonconvex} The dependence of the function $F(\x, s, \p)$ on the third variable $\p$ in  Test 3 (a) and  Test 5 (b). It is evident that in those tests, the function $F$ is nonconvex in $\p.$}
	\end{center}
\end{figure}

\begin{remark}
	It follows from tests 3 and 5  that the convexification method is effective in the interesting case when the function $F$ is nonconvex in $\p.$ 
	The non-convexity is illustrated in Figure \ref{fig nonconvex}. 
	This numerically confirms the strength of our method in solving Hamilton-Jacobi equations.
	We refer the readers to \cite{KaoOsherQian, LiQian} for some other examples dealing with non-convex, discontinuous Hamiltonians via the Lax-Friedrichs sweeping method.
\end{remark}

{
	\begin{remark}
		The relative errors in computation are displayed in the last rows of Figures \ref{fig 1}--\ref{fig 6}.
		It is evident from those figures that the errors in computation occur at $\partial \Omega$ where the noise is added and at the places where the true solution $u_{\rm true}$ is not differentiable.
		This again reflects the strength of the convexification method.
	\end{remark}
}

{
	It has been shown both analytically and numerically   that the convexification method is robust in solving  a large class of Hamilton-Jacobi equations. 
	The strengths of the convexification method involves the facts 
	\begin{enumerate}
		\item that it does not require any special structure of the Hamiltonian; especially, the convexity condition of the Hamiltonian with respect to the third variable is relaxed;
		\item that it yields the satisfactory numerical solutions even when the given boundary is noisy.
	\end{enumerate}
	However, the convexification method has a drawback. 
	It is time consuming in comparison to the well-known methods for nonconvex Hamiltonians; for e.g., the Lax--Friedrichs schemes (\cite{OsherShu, Abgrall, OsFe})  and the Lax--Friedrichs sweeping algorithm (\cite{KaoOsherQian,LiQian}). 
	In our tests, it takes  from 10 minutes to 90 minutes, depending on the forms of the given Hamiltonians, for a Workstations T7810 with 24 cores to compute the solutions (see Tables \ref{tab1}--\ref{tab5}).
	The slow computational time is acceptable in the sense that we consider the computational program as a ``proof of concept" to numerically confirm the analysis of the convexification method.
	The convexification method is the first generalization of the numerical method based on Carleman estimates to solve Hamilton-Jacobi equations. 
	We expect to improve the computational time in the next generation. 
	The next generation will be developed based on the fixed point iterative scheme similar to the ones in \cite{LeetalPreprint2021, LeNguyen:2020, NguyenKlibanov:preprint2021} for quasi-linear elliptic and hyperbolic equations. The rate of convergence in those papers is $O(\theta^n)$ as $n \to \infty$ for some $\theta \in (0, 1).$
	Hence, the success in reducing computational time is very promising. 
}

\section{Concluding remarks}

In this paper, we introduce a new method to solve highly nonlinear Hamilton-Jacobi equations in a rectangular domain. 
This method is called the convexification.
The key idea of the convexification is to involve a Carleman weight function into a cost functional defined directly from the equation under consideration. 
Using a Carleman estimate, we established some important theoretical results. 
The first theorem guarantees that this cost functional is strictly convex and has a unique solution.
Then, we proved in the second theorem that the gradient descent method with sufficiently small step size delivers a sequence converging to the unique minimizer.
Then, in the third theorem, we prove that the minimizer above converges to the solution in the vanishing viscosity process, a good approximation of the viscosity solution to the Hamilton-Jacobi equation,  as the noise contained in the boundary data tends to $0$. 
The rate of the convergence is Lipschitz.
All theoretical results are valid in the framework that we know the value of the solution on the boundary of the domain and its normal derivative in a part of this boundary.
We have pointed out that this framework is acceptable in some real-world applications.
We have shown some interesting numerical tests in 2D. 
These tests confirm the convergence of the convexification method even when the Hamiltonian is not convex or discontinuous. 
Moreover, these numerical results are out of expectation even when the solved equations are not in the framework above.

\medskip

As of now, the convexification method is more time consuming in comparison to the well-known methods as addressed by the end of Section \ref{sec:num}.
We expect to improve the computational time in the next generation. 
Besides, we also intend to study Hamilton-Jacobi equations with other types of boundary conditions in the near future.

\section*{Acknowledgement} The works of MVK and LHN were supported by US Army Research Laboratory and US Army Research
Office grant W911NF-19-1-0044. 
The work of HT is supported in part by NSF grant DMS-1664424 and NSF CAREER grant DMS-1843320.


\end{document}